\documentclass[11pt]{amsart}
\usepackage[utf8]{inputenc}
\usepackage{fullpage}
\usepackage{enumerate}
\usepackage[bbgreekl]{mathbbol}
\usepackage{amsmath,amssymb}
\usepackage{color}
\usepackage{marginnote}
\usepackage[linesnumbered,vlined,ruled]{algorithm2e}
\usepackage{graphicx}
\usepackage{tikz}
\usepackage{dsfont}
\usepackage{graphicx}
\usepackage{placeins}
\usepackage{xcolor}
\usepackage{enumerate}
\usepackage[shortlabels]{enumitem}

\usepackage{bbm}

\usepackage{qsymbols}
\usepackage{graphicx}
\usepackage{placeins}
\usepackage{float}
\usepackage{subfigure}
\usepackage{parskip}

\newtheorem{theorem}{Theorem}
\newtheorem{corollary}{Corollary}
\newtheorem{lemma}{Lemma}
\newtheorem{proposition}{Proposition}

\newtheorem{remark}{Remark}
\newtheorem{definition}{Definition}
\newtheorem{conjecture}{Conjecture}

\newcommand{\Pc}{{\mathcal P}}
\newcommand{\X}{{\mathcal X}}
\newcommand{\HR}{{\mathcal H}}
\newcommand{\N}{{\mathbb N}}
\newcommand{\R}{{\mathbb R}}
\newcommand{\Z}{{\mathbb Z}}
\newcommand{\V}{{\mathbb V}}

\newcommand{\Cc}{{\mathcal C}}
\newcommand{\Sc}{{\mathcal S}}
\newcommand{\supp}{\mathrm{supp}} 
\newcommand{\1}{\mathds{1}}
\newcommand{\B}{{\rho}}
\newcommand{\tz}{z'}
\newcommand{\ttz}{z''}

\begin{document}

\title[Criteria for entropic curvature]{Criteria for entropic curvature\\ on graph spaces}
\author{Martin Rapaport, Paul-Marie Samson}
\thanks{This research is funded in part, by the Agence nationale de la recherche (ANR), Grant ANR-23-CE40-0017, and by the B\'ezout Labex, reference ANR-10-LABX-58. A CC-BY public copyright license has been applied by the authors to the present document and will be applied to all subsequent versions up to the Author Accepted Manuscript arising from this submission, in accordance with the grant’s open access conditions.  The second author is supported by a grant of the Simone and Cino Del Duca Foundation.}
\address{Univ Gustave Eiffel, UPEM, Univ Paris Est Creteil, CNRS, F-77447 Marne-la-Vall\'ee, France}
\email{martin.rapaport@univ-eiffel.fr, paul-marie.samson@univ-eiffel.fr}
\date{\today}

\begin{abstract}
In this paper we  establish new simple  local geometric criteria for discrete entropic curvature introduced in \cite{Sam21} that are powerful enough to capture many geometric properties of complex models arising in mathematical physics. 
These results are robust in the sense that they apply to \emph{any} discrete graph equipped with a Markov reversible generator. Our definitions of entropic curvature differ from the one of the pioneering works of Erbar-Maas [19, 20] (which is already a discrete analog of the Lott-Sturm-Villani   entropic curvature in the continuous setting). Singularly, our results  provide refined concentration properties related to the celebrated convex-hull method by Talagrand \cite{Tal95,Tal96} for a large class of probability measures that cannot be captured from Erbar- Maas entropic definition of curvature. Our approach gives also a new insight of the convex hull method, without being related to induction arguments.

We illustrate the power of our results, as well as the general entropic strategy developed in this   paper, to tackle   challenging models studied in mathematical physics, including Gibbs measures with interaction potentials such as 
 Ising models on the discrete hypercube and    measures with interaction potential on the lattice $\Z^n$. 
For instance,  we significantly improve the constant of the refined convex concentration properties obtained in \cite{AKPS18} for Ising models.  Moreover,
when dealing with the antiferromagnetic Curie-Weiss model, we  improve the previously known bound for  entropic curvature  by a factor of $\sqrt{n}$. Our simple criteria also provides the expected right order of magnitude $C/\sqrt n$ for the lower-bound on the entropic curvature for the renowned Sherrington-Kirkpatrick model from the spin glass theory. This last result is consistant with the recent works \cite{BB19, EKZ22} on the modified logarithmic Sobolev and Poincaré inequalities for the Sherrington-Kirkpatrick model.

\end{abstract}
\maketitle

\tableofcontents
\section{Introduction and outline of the paper}

In this section, we present the goal and structure of the paper and announce the main results, which are provided in the subsequent sections. 

Let $G=(\X,d,L,m)$ be a {\it graph space}, i. e., $\X$ is the set of a  connected undirected graph equipped with its graph distance $d$ and a generator $L$ with reversible measure $m$ (see Section \ref{sectframework}). Since the pioneer works by Erbar-Maas \cite{EM12,EM14}, defining a \textit{good} notion of entropic curvature for such  general discrete  spaces, in comparison to the Lott-Sturm-Villani theory of curvature on geodesics spaces \cite{LV09,Stu06a,Vil09}, has been a challenging problem.  One of the main application of such notion is to capture geometric properties of the space such as measure concentration phenomena for the measure $m$, or functional inequalities related to estimates of the mixing time of the dynamics behind the generator~$L$.

In \cite{Sam21} the second named author, in  continuity of the work by \cite{Leo17},  introduced a displacement convexity property of the { relative entropy} along Schr\"odinger bridges at zero temperature, as  a guideline for other notions of entropic curvature on discrete spaces. Let $\Pc(\X)$ denote the set of probabilities on $\X$, and $\mathcal{P}_b(\X)$ the subset of probabilities  with bounded support. By definition,  the  {\it relative entropy} of a probability measure $q$ on a measurable space $\mathcal Y$ with respect to a probability  measure   $r\in \Pc(\mathcal Y)$  is  given by 
 \[\HR(q|r):= \int_{\mathcal Y} \log(dq/dr) \,dq\qquad \in [0, \infty], \]         
if $q$ is absolutely  continuous with respect to $r$ and    $\HR(q|r):=+\infty$ otherwise. As recalled in \cite{Sam21}, this definition extends to $\sigma$-finite non-negative measures $r$ adding  weak conditions on $q\in \Pc(\mathcal Y)$, and in that case $\HR(q|r)\in (-\infty, \infty]$. 
Schr\"odinger bridges at zero temperature are $W_1$-Wasserstein geodesics on the set of probability measures. The {\it $W_1$-Wasserstein distance}  between two probability measures $\nu_0$ and $\nu_1$ is defined as a minimal cost over all transference plans $\pi\in \Pi(\nu_0,\nu_1)$,  
\begin{eqnarray}\label{defW1}
W_1(\nu_0,\nu_1):=\inf_{\pi\in \Pi(\nu_0,\nu_1)}\iint d(x,y)\,d\pi(x,y),
\end{eqnarray}
 $\Pi(\nu_0,\nu_1)$ denotes the set of probability measures on the product space $\X\times \X$ with first marginal $\nu_0$ and second marginal $\nu_1$. 
We call {\it $W_1$-optimal coupling} of $\nu_0$ and $\nu_1$ any transference plan $\pi$  that achieves the infimum in \eqref{defW1}.
For the sake of consistency, the definition of Schr\"odinger bridge at zero temperature is postponed in the next Section \ref{sectframework}. As explained in \cite{Sam21}, these $W_1$-Wasserstein geodesics capture the geometry of the graph space, since the support of a geodesic between {two Dirac measures}  at $x\in\X$ and $y\in\X$  is the set of vertices which belong to any geodesic on the graph between the vertices $x$ and $y$. According to all these definitions, here is the displacement convexity property of entropy introduced in the seminal paper \cite{Sam21}.
\begin{definition}\label{defcourb}\cite{Sam21} On the graph space $(\X,d,m,L)$, one says that  the 
  relative entropy is $C$-displacement convex where $C=(C_t)_{t\in[0,1]}$,  if for any probability measures $\nu_0,\nu_1\in \mathcal{P}_b(\X)$, there exists a 
  Schr\"odinger bridge at zero temperature $(\widehat \nu_t)_{t\in [0,1]}$ whose structure is given by  \eqref{defhatnut}, and  such that 
   for any $t\in(0,1)$, 
 \begin{eqnarray}\label{deplacebis}
\HR(\widehat \nu_t|m)\leq (1-t) \HR(\nu_0|m)+t \,\HR(\nu_1|m)- \frac{t(1-t)}2C_t(\widehat \pi),
\end{eqnarray}
where   $\widehat \pi$ is the $W_1$-optimal coupling  between $\nu_0$ and $\nu_1$ that appears in the definition  \eqref{defhatnut} of $(\widehat \nu_t)_{t\in [0,1]}$.
\end{definition}
Observe that such a property on graphs was first proposed by  M. Erbar and J. Maas \cite{Maa11,EM12,EM14} where the cost $C_t(\widehat \pi)$ is replaced by $\kappa {\mathcal W}_2^2(\nu_0,\nu_1)$, with $\kappa\in \R$ and where ${\mathcal W}_2$ is an abstract Wasserstein distance on  $\Pc(\X)$.  In their definition, Schr\"odinger bridges at zero temperature are also replaced by a ${\mathcal W}_2$-geodesic, and the best constant $\kappa$ represents the so-called entropic curvature of the space.  Actually, the distance ${\mathcal W}_2$ is defined using a discrete type of Benamou-Brenier formula in order to  provide a Riemannian  structure for the probability space $\Pc(\X)$. This distance ${\mathcal W}_2$ is greater than $\sqrt 2 W_1$. Unfortunately ${\mathcal W}_2$ can not be  expressed as a minimum of a cost  among  transference plans $\pi$ as in the definition \eqref{defW1} of $W_1$. 

This paper shows that the main advantage of the Schr\"odinger method  is its ability to capture such types of costs greater than $W_1$ for almost all graph spaces, among them, the so called  {\it weak optimal transport costs} that  appear in the literature to describe refined concentration phenomena in discrete spaces (see \cite{GRST14bis}). 
Indeed, if $m(\X)<\infty$, then concentration properties for the renormalized probability measure $\mu=m/m(\X)$ are  a straightforward  consequence of the displacement convexity property  \eqref{deplacebis} when $C_t$ represents a positive cost. Since  property \eqref{deplacebis} is invariant under a scaling of the measure $m$, it also holds by replacing the measure $m$ by $\mu$. Therefore since  $\HR(\widehat \nu_t|\mu)\geq 0$, the $C$-displacement convexity property provides the following family of transport-entropy inequality, for any probability measure $\nu_0$ and $\nu_1$ on $\X$ and for any $t\in(0,1)$ \begin{equation}\label{transport0}
\frac12 \inf_{\pi\in \Pi(\nu_0,\nu_1)}C_t(\pi)\leq \frac1{t}\HR(\widehat \nu_0|\mu)+\frac1{1-t}\HR(\widehat \nu_0|\mu).
\end{equation}
Such types of transport-entropy inequalities are closely related to different types of concentration properties depending on the family of transport cost functions $C_t$, {including} among them  {\it weak optimal transport costs} that  appear in the literature to describe refined concentration phenomena in discrete spaces. This  has  been widely studied for different types of discrete measures with different types of costs, we refer to \cite{GRST14bis} for more references on this subject. 

The displacement convexity property \eqref{deplacebis} has been studied in \cite{Sam21} for  toy models of graphs with restricted classes of measures $m$. It was not clear what kind of generic cost functions $C_t$ could be considered in any graph spaces to measure the convexity property of the relative entropy along Schr\"odinger bridges. Moreover the paper \cite{Sam21} does not propose general statements linking  positive entropic curvature to functional inequalities (logarithmic Sobolev inequalities and Poincaré inequalities). 
This paper answer these two fundamental questions, with new applications and results  for complex graph models, validating the Schr\"odinger bridges approach. {To} emphasize the power of our approach, Theorem \ref{thmising}, at the end of this introduction, presents one of the relevant  applications for Ising models of the main general results of this paper.

First, in Section \ref{sectmainresults}, we consider several types of generic cost functions, $T_2$, $W_1$, and a weak transport cost  $\widetilde T$, {which}  define different types of entropic curvature. All these costs are defined with the graph distance $d$. As a main example in this paper, we call {\it entropic curvature} {the largest constant}  constant $\kappa$ for which  the displacement convexity property \eqref{deplacebis} holds with 
\begin{equation}\label{defT_2}
C_t(\pi):=\kappa \, T_2(\pi):= \kappa \iint d(x,y)\big(d(x,y)-1\big)\, d\pi(x,y).
\end{equation}

Similarly one defines the $W_1^2$-entropic curvature and the  $\widetilde T$-entropic curvature.
Theorem \ref{thmprinc} presents general simple geometric criteria to lower bound the  entropic curvature. These criteria are of local nature, they only depend on the geometric structure of the balls of radius two, and the generator $L$. Note that our criteria provide negative lower-bounds on entropic curvature  for trees, and more generally  for geodetic graphs as exposed in Section \ref{sectionneg}. In case of positive lower bounds, Theorem \ref{thmprincbis} gives refined criteria to  lower bound the $W_1^2$-entropic curvature and the  $\widetilde T$-entropic curvature. Then, Theorem \ref{Logsob} establishes that some modified  logarithmic Sobolev inequalities and Poincaré  inequalities  for the measure $\mu$ are consequences of positive $\widetilde T$-entropic curvature.

The results of Section \ref{sectmainresults} are very general since we do not  assume any particular geometric structure for the graph space. Section \ref{sectrefined} complements  Section \ref{sectmainresults} by considering  a large class of graph named {\it structured graphs} for which a refined weak cost functions named $\it \widetilde T_2$ can be considered in the convexity property \eqref{deplacebis}. These class includes the discrete hypercube, the lattice $\mathbb{Z}^n$ and many other Cayley graphs. Theorem \ref{thmstructure} establishes that structured graphs have non-negative entropic curvature and Theorem \ref{Thmstructure} gives other local criteria to lower bound the $\widetilde T_2$-entropic curvature. Finally Theorem~\ref{Logsobbis} shows that positive $\widetilde T_2$-entropic curvature also implies other types of modified  logarithmic Sobolev inequalities and Poincaré inequalities.  

We illustrate the power of our criteria {in handling} challenging Gibbs measures such as Ising models on the discrete hypercube, including the renowned Sherrington-Kirkpatrick model, or alternatively on the lattice  $\mathbb{Z}^{n}$ equipped with measures with interaction potentials (see Section \ref{sectionpertur} and Section \ref{section-interaction-pot}). The results of this paper, though general, offer new insights and particularly significant findings for such measures. 
A measure with interaction potential on a graph can be interpreted as a perturbation of the counting measure on the graph. In the first part of Section  \ref{sectionpertur}, we show how  potential perturbations affect the cost function or our main criteria for entropic curvature. We apply then these results in the next section on the discrete hypercube (see Proposition \ref{prophypercube}) and on the lattice $\Z^n$ (see Proposition \ref{propZn}).   In the second part of Section   \ref{sectionpertur}  we show that the convexity property \eqref{deplacebis} is stable by restriction of the measure $m$ to {\it convex subsets} (see definition  \eqref{defconvset}). 

Section \ref{studycrit} concerns the study of the main criteria of the paper. We first present tensorisation properties of the main criteria in Theorem \ref{ttens}. Then, assuming that the graph space is equipped with the counting measure, we discuss necessary or sufficient conditions for the criteria to get positive entropic curvature. We propose and discuss a necessary and sufficient condition (see Conjecture \ref{conj1}). Also, for the class of graphs for which there exists at most two midpoints between  two points at distance two, we show that our criteria is  related to the {\it Motzkin-Strauss Theorem} (see Theorem \ref{cliquethm} and Proposition \ref{graphMS})

In the last section, Section \ref{sectcomp},
we make simple  remarks to compare  our criteria with other notions of curvature. All the proofs of the results of this paper are postponed in Appendix A and B.

For the non\textcolor{blue}{-}specialist reader, and in order to emphasize the {robustness}  of our criteria, we end this presentation with concentration properties and modified logarithmic Sobolev inequalities, consequences of our results for Ising models. 
Let $\Lambda$ be a finite set of vertices of cardinality $n$ of a finite graph $G_{\Lambda}$ 
without  multiple edges and without loops.  One denotes $i\sim_{\!\Lambda} j$ if two vertices $i$ and $j$ of $\Lambda$ are neighbors. Let $m_0$ denote the counting measure on $\{-1,1\}^\Lambda$ and  $m_w=e^{-w}m_0$ be the measure on $\{-1,1\}^\Lambda$ with potential of interaction defined as
\begin{equation}\label{defwising}
w(z):=-\sum_{i\in \Lambda} T_i z_i -\frac{\beta}{2} \sum_{(i,j)\in \Lambda^2,i{\sim_{\!\Lambda}}j} W_{ij}\, z_iz_j, \qquad z\in\{-1,1\}^\Lambda,
\end{equation}
where  $\beta>0$ represents an inverse temperature parameter, and the $T_i$'s and $W_{ij}$'s are real interaction parameters with $W_{ij}=W_{ji}$. Let $\mu_w:=m_w/m_w\big(\{-1,1\}^\Lambda\big)$ be the associated renormalized probability measure.
Let $\lambda_{\max}(W)$ denotes the largest eigenvalue of the  $n$ by $n$ symmetric matrix $W$, with interaction coefficient $W_{ij}$ for $i{\sim_{\!\Lambda}}j$,  and 0 otherwise. Let $\lambda_{\max}|W|$ denotes the largest eigenvalue of the  symmetric matrix with coefficients $|W_{ij}|$, with interaction coefficient $W_{ij}$. Define also $|W|_{\max}=\max_{\{i,j\}\subset \Lambda} |W_{ij}|$. 
Let 
\begin{equation}\label{defrhoising}
\rho_\beta(W):=1- 2\beta{\lambda_{\max}(W)} - 2\beta{\lambda_{\max}|W|} \,k\big({2\beta|W|_{\max}}\big),
\end{equation}
where $k(s):=\frac1s(e^s-s-1)$, $s>0$. 
In Section \ref{sectionhypecube}, Proposition \ref{prophypercube} ensures that if $\rho_\beta(W)>0$ then the entropic curvatures defined in this paper are positive,  from which we  derive the following results.
\begin{theorem}\label{thmising} According to the above notations, assume that  $\rho_\beta=\rho_\beta(W)>0$, 
\begin{itemize}
    \item The probability measure $\mu_w$  satisfies the following concentration properties. 
    For any 1-Lipschitz function $f:\{-1,1\}^\Lambda\to \R$ with respect to the graph distance $d(x,y)=\sum_{i\in \Lambda} \1_{x_i\neq y_i}$, $x,y\in \{-1,1\}^\Lambda$, one has  for all $s\geq 0$, 
    \[\mu_w(f\geq \mu_w[f]+s) \leq e^{-{2\rho_\beta s^2}/{n}}.\]
    For any 1-Lipschitz convex function $f:[-1,1]^\Lambda\to \R$ with respect to the Euclidean metric, one has  for all $s\geq 0$,
    \[\mu_w(f\geq \mu_w[f]+s)\leq e^{-  {\rho_\beta s^2}/8} \quad\mbox{and}\quad \mu_w(f\leq \mu_w[f]-s)\leq e^{-{\rho_\beta s^2}/{8}}.\]
    \item The probability measure $\mu_w$ satisfies the following modified logarithmic Sobolev inequality, for any function $f:\{-1,1\}^\Lambda\to (0,+\infty)$,
    \begin{equation*}
{\rm Ent}_{\mu_w}(f)\leq \frac 1{\rho_\beta} \int \sum_{i\in[n]} {[ \log f(\sigma_i(x))- \log f(x)]_{_-} [ f(\sigma_i(x))-  f(x)]_{_-}} \,d\mu_w(x),
\end{equation*}
   where $[u]_{_-}=\max(0,-u),u\in \R$, ${\rm Ent}_{\mu_w}(f):=\mu_w[f\log f]-\mu_w[f]\log \mu_w[f]$, and given a configuration $x\in \{-1,1\}^\Lambda$ and $i\in \Lambda$,  $\sigma_i(x)$ denotes the configuration that differs from $x$ only at vertex $i$.  
    \end{itemize}
    \end{theorem} 


Let us provide a few comments that are further elaborated in Section \ref{sectionhypecube}.
As a first example, consider the so-called simple Curie-Weiss model, where $W$ is the adjacency matrix of a complete graph $G_{\Lambda}$, with $W_{ij}=1$ for $i \neq j$. In this case, $\lambda_{\max}(W)=n-1$, and $\beta < \frac{1}{2n}$ is a sufficient condition for $\rho_\beta(W)= 1 - (n-1)(e^{2\beta}-1) > 0$.
The critical value $\overline{\beta}_n$ for the Curie-Weiss model, beyond which some Poincaré inequality is known to fail, is $\overline{\beta}_n=\frac{1}{n}$ (see \cite{DLP09}). Achieving positive entropic curvature up to this critical value remains a challenging problem.

Considering the antiferromagnetic Curie-Weiss model where  $W_{ij}=-1 $ for $i\neq j$, the condition improved by a factor $\sqrt{n}$. Indeed, since in that case  $\lambda_{\max}(W)=1$ and  $\lambda_{\max}|W|=\lambda_{\max}(A)=n-1$, one gets 
\[\rho_\beta=1-2\beta-2(n-1)\beta\,k(2\beta)\geq 1-2n\big(e^\beta-1)^2.\] 
It follows that $\rho_\beta>0$ as soon as $\beta\leq \frac1{1+\sqrt{2n}}$. Clearly this improvement follows from the fact that the first order Taylor expansion in $\beta$ of $\rho_\beta$ only involves $\lambda_{\max}(M)$ which is much smaller than the quantity $|W|_{\max} \lambda_{\max}|M|$ appearing in the second order term, namely,
\[\rho_\beta=1-2\beta{\lambda_{\max}(W)}-2\beta^2 |W|_{\max} \lambda_{\max}|M| +o(\beta^2).\]
Note that we have not found references where similar results have been obtained for the antiferromagnetic Curie-Weiss model.
As explained in Section \ref{sectionhypecube}, many results in the literature are obtained under   stronger Dobrushin-type conditions for which $\beta |W|_{\max} \max_{i\in \Lambda}\sum_{j, j\sim_{\! \Lambda}i} |W_{ij}|$ needs to be bounded by a numerical constant. 
This condition is stronger since the Perron-Frobenius theorem gives $\lambda_{\max}(M)\leq \lambda_{\max}|M|$.

The fact that our condition mainly concerns the largest eigenvalue of $M$ also allow us to consider random interaction coefficients $W_{ij}$ when the order of $\lambda_{\max}(W)$ is well understood for large $n$. This is the case of the so-called Sherrington-Kirkpatrick model from spin glass theory \cite{Tal10} for which $\lambda_{\max}(W)$ is almost surely equivalent to $2\sqrt n$. As a consequence, we prove that $\rho_\beta>0$ as soon as $\beta=\beta_n\leq \frac{1-\varepsilon}{4\sqrt n}$, $\varepsilon >0$.

As explained with more details in Section \ref{sectionhypecube}, it follows from these remarks that the concentration results  of Theorem \ref{thmising} strongly improve the convex concentration results of \cite[Proposition 5.4]{AKPS18} with Dobrushin-type conditions, and that the results dealing with modified logarithmic Sobolev 
inequalities are comparable to the ones by Bauerschmidt-Bodineau \cite{BB19} and  Eldan-Koehler-Zeitouni \cite{EKZ22} for the Sherrington-Kirkpatrick model.

\section{Framework of the paper - Schr\"odinger bridges at zero temperature}\label{sectframework}
Let $\X$ be the set of vertices of a connected undirected graph $G=(\X,E)$ where $E$ denotes the set of edges, without multiple edges and without loops. Two vertices $x$ and $y$ are {\it neighbours} if $\{x,y\}$ is an edge of $E$, we  write $x\sim y$ in this case. Let $d$ denote the combinatorial graph distance, so that $d(x,y)=1$ if and only if  $x\sim y$.
The graph $G$ is supposed to be \textit{locally finite}, that is, the \textit{vertex degree} $\text{deg}(x):=\sum_{y\sim x} 1$ of any $x\in \mathcal{X}$ is finite. The \textit{maximal degree} of the graph is denoted by   $\Delta(G):=\sup_{x\in \mathcal{X}}\text{deg}(x)\in \N\cup\{+\infty\}$. A {\it discrete geodesic path} $\gamma$ joining two  vertices $x$ and $y$ is a sequence of  neighbours of minimal size $d=d(x,y)$:  $\gamma=(z_0,\ldots,z_d)$ with $z_0=x$ and $z_d=y$ and for any $i\in[d]:=\{1,\cdots ,d\}$, $\{z_{i-1},z_i\}\in E$. In the sequel,  $z\in\gamma$ means that there exists $i\in\{0,\ldots, d\}$ such that $z=z_i$, and  $(z,w)\in \gamma$ means that there exists $0\leq i<j\leq \ell$ such that $z=z_i$ and $w=z_j$. Let $G(x,y)$ be the set of all geodesic paths joining $x$ to $y$, and let   $[x,y]$ be the set of all vertices that belong to a geodesic from $x$ to $y$,
\[[x,y]:=\big\{z\in \X\,\big|\, z\in \gamma,\gamma\in G(x,y)\big\},\]  
$]x,y[:=[x,y]\setminus\{x,y\}$ and $[x,y[:=[x,y]\setminus \{y\}$. More generally, given two subsets  $A$ and $B$ of $\X$, one defines
\[[A,B]:=\bigcup_{x\in A,y\in B} [x,y]\quad \mbox{and} \quad ]A,B[:=\bigcup_{x\in A,y\in B} ]x,y[.\]

The set $\X$ is endowed with the $\sigma$-algebra generated by singletons. 
The subset of probability measures $\mu$ satisfying $\int d(x_0,y) d\mu(y)<\infty$ for some $x_0\in\X$, denoted by $\Pc_1(\X)$, can be  endowed with the $W_1$-Wasserstein distance defined by \eqref{defW1}.
For any non negative measure $M$ on a measurable space $\mathcal Y$, $\supp(M)$ denotes the support of this measure.  For further use, note that any measure $\pi\in\Pi(\nu_0,\nu_1)$ admits the two  following decompositions : for any $(x,y)\in \supp(\pi)$, 
 \[ \pi(x,y)=\nu_0(x)\,{\pi}_{_\rightarrow}(y|x)= \nu_1(y)\,{\pi}_{_\leftarrow}(x|y),\]
 defining thus two Markov kernels ${\pi}_{_\rightarrow}$ and ${\pi}_{_\leftarrow}$. 

On a discrete space $\X$, recall that any generator $L$ acting on functions from $\X$ to $\R$ is entirely given by the jump rates from $x\in \X$ to $y\in \X$ denoted by $L(x,y)$, $L(x,y):=L\delta_y(x)$ with $\delta_y(y)=1$ and $\delta_y(z)=0$ for $z\neq y$.

In this paper, we call \textit{graph space} any locally finite graph $G=(\X,E)$ as above endowed with a reference measure $m$ on $\X$ and a generator $L$ satisfying the following two conditions:
\begin{itemize}
\item[$\cdot$] The measure $m$ is reversible with respect to $L$, namely  for any $x,y\in\X$,
\[m(x)L(x,y)=m(y)L(y,x).\]
\item[$\cdot$] 
For any $x,y\in\X$, one has 
\begin{equation}\label{PropL}
L(x,y)>0 \quad\mbox{ if and only if }\quad  x\sim y,
\end{equation}
(and $ L(x,x):=- \sum_{y\in \X, y\neq x} L(x,y)$). 
\end{itemize}
For simplicity, one also denotes $G=(\X,d,L,m)$ in that case. Note that in this paper, we are able to relax strong technical assumptions given in the seminal paper \cite[Assumptions (12)-(15)]{Sam21}, allowing consideration of spaces such as $M/M/\infty$ processes on $\mathbb{N}$ with Poisson stationary measure, which were previously out of scope.

In this paper, given $m$, two specific generic generators $L_1$ and $L_2$   will be considered at different places defined by, for all $x\neq y$
\begin{equation}\label{defL1L2}
L_1(x,y)=\Big(\frac{m(x)}{m(y)}\Big)^{1/2}\1_{x\sim y}\quad\mbox{and} \quad L_2(x,y)=\frac 12 \Big(1+\frac{m(x)}{m(y)}\Big)\1_{x\sim y} .
\end{equation}
Note that  for the counting measure on $\X$, denoted by $m_0$ in this paper, these two generators are the same and one denotes in that case $L_0:=L_1=L_2$.

By definition, for $x,y\in \X$ and $k\in \N^*$, one denotes  \[L^k(x,y):=\sum_{z_1, z_2, \ldots,z_{k-1}\in\X} L(x, z_1)L(z_1,z_2)\cdots L(z_{k-1},y).\] 
The property \eqref{PropL} ensures that 
\[ L^{d(x,y)}(x,y)=\sum_{\alpha\in G(x,y)} L(\alpha), \quad \mbox{where} \quad L(\alpha):= L(z_0, z_1)\cdots L(z_{d(x,y)-1},z_{d(x,y)}),\]
for any $\alpha=(z_0,z_1,\ldots, z_{d(x,y)})\in G(x,y)$.

As defined in \cite{Sam21}, given  $\nu_0,\nu_1\in \mathcal{P}_b(\X)$, {\it a Schr\"odinger bridge at zero temperature}, denoted by $(\widehat \nu_t)_{t\in [0,1]}$ in the present paper, is a particular $W_1$ constant speed geodesic between $\nu_0$ and $\nu_1$ on $\Pc(\X)$, namely,  $ \widehat \nu_0= \nu_0$, $ \widehat \nu_1= \nu_1$, and for any $0\leq s\leq t\leq 1$,
 \[W_1\big(\widehat \nu_t,\widehat \nu_s\big)=(t-s) W_1(\nu_0,\nu_1).\]
Such a  path is obtained from a mixture of Schr\"odinger bridges,  by a slowing down procedure as a temperature term goes to zero due to C. Léonard \cite[Theorem 2.1]{Leo16} (see also \cite{Sam21}). These geodesic paths are mixture of $W_1$-constant speed geodesics $\nu_t^{x,y}$ between the Dirac measures $\delta_x$ at $x\in \supp(\nu_0)$ and $\delta_y$ at $y\in\supp(\nu_1)$, according to a coupling $\widehat \pi\in \Pi(\nu_0,\nu_1)$. Observe that given  bounded marginals $\nu_0$ and $\nu_1$, the L\'eonard slowing down procedure selects a single  coupling $\widehat \pi\in \Pi(\nu_0,\nu_1)$ if some conditions are satisfied on the underlying space \cite[Result 0.3]{Leo16}.  The main property of $\widehat \pi$ is to be a $W_1$-optimal coupling.  As explained in \cite{Sam21}, the structure of Schr\"odinger bridges at zero temperature that we also consider in this paper is the following: 
for any $z\in \X$
\begin{equation}\label{defhatnut}
\widehat \nu_t(z):=\sum_{x,y\in \X} \nu_t^{x,y}(z) \,\widehat \pi(x,y),
\end{equation}
with for any $x,y\in \X$,
\begin{equation}\label{pathdirac}
\nu_t^{x,y}(z):=
\1_{[x,y]}(z)\, r(x,z,z,y)\,\B_t^{d(x,y)}(d(x,z)),
\end{equation}
where for $x,z,v,y\in\X$,
\begin{equation}\label{defrpont}
 r(x,z,v,y):= \frac{L^{d(x,z)}(x,z) L^{d(v,y)}(v,y)}{L^{d(x,y)}(x,y)},
 \end{equation}
and $\B_t^d$ denotes the binomial law with parameter $t\in [0,1]$, $d\in \N$ :
\[ \B_t^d(k):= \binom{d}{k}\, t^k(1-t)^{d-k},\quad k\in\{0,\ldots,d\},\]
 with the binomial coefficient $\binom{d}{k}:=\frac{d!}{k!(d-k)!}$.
 All along the paper one omits the dependence in $\nu_0$ and $\nu_1$ of $(\widehat \nu_t)_{t\in [0,1]}$ and $\widehat \pi$ to lighten the notations.

\section{Main results for any graph space}\label{sectmainresults}
In this part, we focus on the convexity property \eqref{deplacebis} for very general graph spaces, which means without a particular geometric structure. One of the main result of this paper, Theorem \ref{thmprinc} below,  introduces a uniform local criteria on the graph space  $(\X,d,m,L)$ under which the  cost $C_t(\widehat \pi)$ can be replaced by $\kappa\, T_2(\widehat \pi)$  where the definition of the cost $T_2$ is given by \eqref{defT_2}.
Recall that by definition, we call  {\it entropic curvature} of the graph space $(\X,d,m,L)$ denoted by $\kappa$ the supremum of $k\in \R$ such that the $C$-displacement convexity property $\eqref{deplacebis}$ holds with $C_t=k \,T_2$. Observe that if $\kappa=+\infty$ then $\eqref{deplacebis}$ ensures that $T_2(\widehat \pi)=0$ since $H(\widehat \nu_t|m)<+\infty$ for any $\nu_0,\nu_1\in \Pc_b(\X)$. 
As a convention $\kappa T_2(\widehat \pi)=0$ if $\kappa=+\infty$. 
When $\kappa>0$ (respectively $\kappa\geq 0$), one says that the space $(\X,d,m,L)$ has positive entropic curvature (respectively non-negative entropic curvature).  
More generally, given a family of  cost functions $c=(c_t)_{t\in (0,1)}$, $c_t:\N\to \R$, we call  {\it $T_c$-entropic curvature} of the graph space $(\X,d,m,L)$ the best constant $\kappa_c\in \R\cup \{+\infty\}$ such that the $C$-displacement convexity property $\eqref{deplacebis}$ holds for all $t\in(0,1)$ with $C_t=\kappa_c \,T_{c_t}$ with 
\[ T_{c_t}(\widehat \pi):=  \iint  c_t(d(x,y))\, d {\widehat \pi}(x,y).\]

Similarly, let us also introduce a definition of {\it ${\widetilde T}$-entropic curvature}  as the best constant $\widetilde{\kappa}\in \R$ such that \eqref{deplacebis} holds with
$C_t= \widetilde{\kappa}\,\widetilde{T}$ and $\widetilde{T}=\max\big(\widetilde{T}_{_\rightarrow},\widetilde{T}_{_\leftarrow}\big),$
where
\[\widetilde{T}_{_\rightarrow}(\widehat \pi):=  \int  \left(\int d(w,x) \,d {\widehat \pi}_{_\rightarrow}(w|x)\right)^2 d\nu_0(x)\]
and \[ \widetilde{T}_{_\leftarrow}(\widehat \pi):= \int \left(\int d(w,y) \,d {\widehat \pi}_{_\leftarrow}(w|y)\right)^2 d\nu_1(y).\]
Let us note that as soon as $\X $ is not reduced to a singleton, there always exist $\nu_0$ and $\nu_1\in 
{\mathcal P}_b(\X)$ such that $\widetilde{T}_{_\rightarrow}(\widehat \pi)>0$ or $\widetilde{T}_{_\leftarrow}(\widehat \pi)>0$ in such a way that $\widetilde \kappa<+\infty$. Observe  that due to the convexity property in $p$ of the cost function $c:(w,p)\in \X\times \mathcal P(\X)\mapsto \left(\int d(w,y) \,d p(y)\right)^2$ in the definition of $\widetilde{T}_{_\rightarrow}$ or $\widetilde{T}_{_\leftarrow}$, these costs are part of the family of the so-called {\it weak transport costs} introduced in \cite{GRST14}, in the continuity of Marton's and Talagrand's works.   

Analogously,  we also call {\it $W_1^2$-entropic curvature} of the graph space $(\X,d,m,L)$ the best constant $\kappa_1\in \R$ such that \eqref{deplacebis} holds with
$C_t = \kappa_1 \,W_1^2.$
By the Cauchy-Schwarz inequality and since $\widehat \pi$ is a $W_1$-optimal coupling, one has $\widetilde{T}_{_\rightarrow}(\widehat \pi)\geq W_1^2(\nu_0,\nu_1)$ and therefore, if $\widetilde{\kappa}\geq 0$ then one has  
\[\kappa_1\geq \widetilde{\kappa}.\]
Further in this introduction, one gives examples of graphs for which this inequality is strict. A second main result of this paper, Theorem \ref{thmprincbis} below, presents refined criteria for non-negative  {\it $W_1^2$-entropic curvature}, non-negative {\it ${\widetilde T}$-entropic curvature} and also some non-negative
 {\it $T_c$-entropic curvature} for a specific family of costs $c=(c_t)_{t\in (0,1)}$, related to the cost $T_2$. 

Let us now introduce the key quantities that allow to define the criteria of  Theorem \ref{thmprinc} and Theorem \ref{thmprincbis}.
For $z\in \X$, let $B_1(z):=\{w\in \X\,|\, d(z,w)\leq 1\}$ denotes \textit{the ball of radius one} centered at $z$, and  for $k=1$ or $k=2$  let \textit{the combinatorial sphere} $S_{k}(z)$ denotes the set of vertices at distance $k$ from $z$
\[S_k(z):=\Big\{w\in \X\,\Big|\, d(z,w)=k\Big\}.\]
Given a vertex $z\in\X$ and a 
 subset
$W\subset S_2(z)$, 
the next  non-negative key quantity  will be used for local lower bound on entropic curvature 
\begin{multline}\label{defR_2}
K_L(z,W):=\sup  \Biggl\{  \sum_{\ttz\in W} L^2(z,\ttz) \prod_{\tz\in  ]z,\ttz[}\left(\frac{\alpha(\tz)}{L(z,\tz)}\right)^{ 2\ell(z,\tz,\ttz) }\Biggl|\, {\alpha}:]z,W[\to \mathbb{R}_{+},\sum_{v \in ]z,W[} \alpha(v)=1 \Biggr\} ,
\end{multline}
with $\ell(z,\tz,\ttz):=\frac{L(z,\tz)L(\tz,\ttz)}{L^2(z,\ttz)}$.
To simplify the notations, one  omits the dependence in $L$ and notices $K=K_L$ when there is no possible confusions.  
In this definition as in the all paper,  we use the convention that a sum indexed by an empty set is $0$. Therefore,      $K(z,W)=0$ holds if and only if $W=\emptyset$. 
One may easily check that given $z$ the quantity $K(z,W)$ is increasing in $W$. Namely,  for  $W\subset W'\subset S_2(z)$, it holds 
\begin{equation}\label{enfin}
K(z,W) \leq K(z,W')\leq  K\big(z,S_2(z)\big).
\end{equation}

For more comprehension about this quantity, consider the special case where $m=m_0$ is the counting measure and $L=L_0$. 
Let $|A|$ denote the cardinal of any finite set $A\subset \X$.  Then \eqref{defR_2} becomes 
\begin{equation}\label{defRbis}
K_0(z,W):=K_{L_0}(z,W)=\sup_{\alpha} \Biggl\{ \sum_{\ttz\in W} \big|]z,\ttz[\big| \Big(\prod_{\tz\in  ]z,\ttz[} {\alpha(\tz)}\Big)^{\frac{2}{|]z,\ttz[|}}\Biggr\},
\end{equation}
where the supremum runs over all $\alpha:]z,W[\to \mathbb{R}_{+}$ such that $\sum_{v \in ]z,W[} \alpha(v)=1$.
Note that given a vertex $z\in \X $, the existence of edges between vertices exclusively within  $S_ {1}(z)$ or $S_{2}(z)$ respectively does not change the value of $K\big(z, S_{2}(z)\big)$. 

In order to introduce a first general criterion that provides lower bounds on entropic curvature, one needs to recall few  notions dealing with optimal transport.
According to  the theory of optimal transport, the support of any $W_1$-optimal coupling is {\it $d$-cyclically monotone} (see \cite[Theorem 5.10]{Vil09}). Recall that by definition, a subset $S\subset \X\times \X$ is $d$-cyclically monotone if for any family $(x_1,y_1),\ldots,(x_N,y_N)$ of points in $S$,
 \[\sum_{i=1}^N d(x_i,y_i)\leq \sum_{i=1}^N d(x_i,y_{i+1}),\]
 with the convention $y_{N+1}=y_1$. As a trivial example, for any $z\in \X$ and  $V\subset \X$, the set $\{z\}\times V$ is a $d$-cyclically monotone subset of $\X\times \X$.
Given  such a subset $S$, recall that due to the $d$-cyclically monotonicity property, if $\{z,z'\}$ is an edge shared by two geodesics $\gamma$ and $\gamma'$, the first from $x$ to $y$ with $(x,y)\in S$ and the second from $x'$ to $y'$ with $(x',y')\in S$, then  either $(z,z')\in \gamma\cap \gamma'$ or $(z',z)\in \gamma\cap \gamma'$ (see \cite[Lemma 4.3]{Sam21}).  It follows that the $d$-cyclically monotone subset $S$ induces a direction on each edge shared by geodesics  from $x$ to $y$ for $(x,y)\in S$. 
One defines $Z(S)$ as the set of vertices which belong to a geodesic from $x$ to $y$ for $(x,y)\in S$,
\[Z(S):=\Big\{z\in \X\,\Big|\,\exists (x,y)\in S, z\in [x,y]\Big\}.\]
Let also define
\[C_{_\rightarrow}(S):=\Big\{(z,w)\in\X\times \X\,\Big|\,z\neq w, \exists (x,y)\in S, \exists \gamma\in G(x,y), (z,w)\in \gamma\Big\},\]
  \[C_{_\leftarrow}(S):=\Big\{(z,w)\in\X\times \X\,\Big| (w,z)\in C_{_\rightarrow}(S)\Big\}.\]
  and for $z\in Z(S)$,
  \[ V^S_{_\rightarrow}(z):=\Big\{\tz\in S_1(z)\,\Big|\, (z,\tz)\in C_{_\rightarrow}(S)\Big\},\quad\quad  V^S_{_\leftarrow}(z):=\Big\{\tz\in S_1(z)\,\Big|\, (z,\tz)\in C_{_\leftarrow}(S)\Big\} . \] 
  $V^S_{_\rightarrow}(z)$ can be interpreted as neighbours of $z$ in the direction given by the $d$-cyclicaly monotone subset $S$ and $V^S_{_\leftarrow}(z)$ as neighbours of $z$ in the opposite direction.  Let us note that there may be vertices in $S_{1}(z)$ which do not belong to $V^S_{_\rightarrow}(z)$ nor to $V^S_{_\leftarrow}(z)$. 
  Setting  $\ell(z,\tz,\ttz):=\frac{L(z,\tz)L(\tz,\ttz)}{L^2(z,\ttz)}$ and $\mathcal L^2(z,\ttz):=\sqrt{L^2(z,\ttz)L^2(\ttz,z)}$, let 
\[K_L(S)=K(S):=\sup_{\alpha,\beta}\Biggl\{\sum_{(z,\ttz)\in C_{_\rightarrow}(S), d(z,\ttz)=2} \mathcal L^2(z,\ttz) \prod_{\tz\in  ]z,\ttz[}\left(\frac{\beta(z,\tz)}{L(z,\tz)}\frac{\alpha(\ttz,\tz)}{L(\ttz,\tz)}\right)^{\ell(z,\tz,\ttz)}\Biggl\}, \]
where  the supremum runs over all non negative functions $\alpha$ and $\beta$ on $\X\times\X$ satisfying  
\[\sum_{z\in Z(S)}\Big(\sum_{\tz\in V^S_{_\rightarrow}(z)} \beta(z,\tz)\Big)^2=1 \quad\mbox{and}\quad   \sum_{\ttz\in Z(S)}\Big(\sum_{\tz\in V^S_{_\leftarrow}(\ttz)} \alpha(\ttz,\tz)\Big)^2=1.\]
Observe that if there exists $(x,y)\in S$ such that $d(x,y)\geq 2$ then $(x,y)\in\big\{(z,\ttz)\in C_{_\rightarrow}(S), d(z,\ttz)=2\big\}\neq \emptyset$ and therefore $K(S)>0$.

For the uniform measure $m_0$ on $\X$ with $L=L_0$, one has  
\begin{equation}\label{defK_0(S)}
K_0(S):=K_{L_0}(S)=\sup_{\alpha,\beta}\Biggl\{\sum_{(z,\ttz)\in C_{_\rightarrow}(S), d(z,\ttz)=2} \big|]z,\ttz[\big|\prod_{\tz\in  ]z,\ttz[}\big({\beta(z,\tz)}{\alpha(\ttz,\tz)}\big)^{\frac1{|]z,\ttz[|}}\Biggl\}
\end{equation}

The quantity $K_L(S)$ is upper-bounded by  $\sup_{z\in \X} K(z,S_2(z))$. Indeed, the reversibility property implies $\ell(z,\tz,\ttz)=\ell(\ttz,\tz,z)$ and one has $\sum_{\tz\in]z,\ttz[}\ell(z,\tz,\ttz)=1 $. Therefore by setting 
\[\mathbb V^S_{_\rightarrow}(z):=\Big\{\ttz\in S_2(z)\,\Big|\, (z,\ttz)\in C_{_\rightarrow}(S)\Big\}\quad \mbox{and}\quad \mathbb V^S_{_\leftarrow}(z):=\Big\{\ttz\in S_2(z)\,\Big|\, (z,\ttz)\in C_{_\leftarrow}(S)\Big\}, \]  
for $z\in Z(S)$, the Cauchy-Schwarz inequality provides 
\begin{align}\label{compK}
 K(S)^2& \leq \sup_{\beta}\Biggl\{\sum_{z\in Z(S)} \sum_{\ttz\in \V^S_{_\rightarrow}(z)} L^2(z,\ttz) \prod_{\tz\in  ]z,\ttz[}\left(\frac{\beta(z,\tz)}{L(z,\tz)}\right)^{2\ell(z,\tz,\ttz)}\Biggl\}\nonumber\\
 &\qquad\qquad\qquad \cdot\sup_{\alpha}\Biggl\{\sum_{\ttz\in Z(S)} \sum_{z\in \V^S_{_\leftarrow}(\ttz)} L^2(\ttz,z) \prod_{\tz\in  ]z,\ttz[}\left(\frac{\alpha(\ttz,\tz)}{L(\ttz,\tz)}\right)^{2\ell(\ttz,\tz,z)}\Biggl\}\nonumber \\
 &\leq \sup_{z\in Z(S)} K\big(z,\V^S_{_\rightarrow}(z)\big)\,  \sup_{\ttz\in Z(S)} K\big(\ttz,\V^S_{_\leftarrow}(\ttz)\big) \leq  \sup_{z\in \X} K(z,S_2(z))^2,
\end{align}
where the last inequality is a consequence of the monotonicity property \eqref{enfin}.

Here is one of the main result of this paper whose proof is given in Appendix B. We emphasize the local-global aspect of this result, which presents a local criterion that enables the derivation of entropic global lower bounds.  Remarkably, it provides a straightforward criterion for efficiently establishing a lower bound on the entropic curvature of any graph space $(\X, d, m, L)$. The cases of the hypercube and the lattice $\Z^n$ are deeply studied in Section \ref{section-interaction-pot}.  Some examples of graphs are not included in this paper, such as the Bernoulli Laplace model, the hard-core model, the multibinomial law, or the transposition model among others and are part of current research that hold potential promising results.

\begin{theorem}\label{thmprinc} Let $(\X,d,m,L)$ be a graph space. Let  \[K=K_L:=\sup_{S} K_L(S) ,\]
where the supremum is over all $d$-cyclically monotone subsets of $\X\times\X$.
Then  the entropic curvature $\kappa$ of the space $(\X,d,m,L)$  is bounded from below by $r=r^L:=-2\log K$ if $K>0$, and $\kappa=+\infty$ if $K=0$. 
\end{theorem}

{\bf Comments:}
\begin{enumerate}[label=(\roman*)]
\item Observe that $K=0$   if and only if $G=(\X,E)$ is a complete graph. Indeed, if $K=0$ then for any singleton $\{(x,y)\}$ of $\X\times \X$, one has $K(\{(x,y)\})=0$ and therefore $d(x,y)\leq 1$.
\item 
For most of graphs,  computing $K$ is not easy since we need first to  characterize the set of  $d$-cyclically monotone subsets. However one may use the estimate given by inequality \eqref{compK} which is more tractable. More precisely, for $z\in \X$, let us define 
\begin{equation}\label{defr}
   r(z)=r^L(z):=-2 \log K_{L}\big(z,S_2(z)\big). 
\end{equation}
Then, according to \eqref{compK}, Theorem \ref{thmprinc} ensures that $\kappa\geq  r\geq\inf_{z\in\X} r(z)$.
Actually, the quantity
$r(z)$ can be interpreted as a local lower  bound on entropic curvature at vertex $z$.
This quantity only depends on the value of the jump rates on the ball of radius 2 centered at $z$. If the graph space is equipped with the counting measure ($m=m_0$ and $L=L_0$), then  $r_0(z):=r^{L_0}(z)$ only depends on the structure of the ball $B_{2}(z)$. Therefore a  lower bound on this local quantity can be  interpreted as a geometric property of the balls of radius 2. The lower bound  on entropic curvature $\inf_{z\in\X} r(z)$ is really useful since it can be estimated on a wide range of graphs. Morevover, as shown in Section \ref{examplesZcube}, this lower bound  provides new  results for complex measures with interaction potentials, like for Ising models. Therefore it should give promissing results for many other specific graph spaces which are not considered in this paper. 

\item\label{distance1} For more comprehension, let us present a simple necessary condition for $r>0$, as $m=m_0$ is the counting measure. According to \eqref{defK_0(S)}, if for some $z_0\in \X$, there exists $z''_0\in S_2(z_0)$ such that $\big|]z_0,z''_0[\big|=1$, or equivalently $]z_0,z''_0[=\{z'_0\}$, then for  $S_0:=\{(z_0,z''_0)\}$ one has $K(S_0)=1$ and  therefore $r\leq 0$.
 As a consequence,  if the criteria of Theorem \ref{thmprinc} provides positive entropic curvature for the space $(\X,d,m_0,L_0)$, 
 then there is  at least 2 midpoints between any two vertices at distance 2. 
 \end{enumerate}

Theorem \ref{thmprinc} ensures that if $K\leq 1$  then the graph space $(\X,d,m,L)$ has non-negative  entropic curvature, i.e. $\kappa\geq  0$. Section \ref{sectionneg} and Section \ref{sectionpos} are devoted to the study of the upper bound $\sup_{z\in \X}K_0(z,S_2(z))\geq K$ in order to give necessary or sufficient conditions for non-negative or positive entropic curvature for the space $(\X,d,m_0,L_0)$.  The next  result states that $K\leq 1$  also implies non-negative $\widetilde{T}$ and $W_1$-entropic curvature. In order to get refined  lower bound on the $W_1$-entropic curvature one introduces the following quantity: for any  $d$-cyclically monotone subset $S$, let 
\begin{align*}
    R_1&(S)=R_1^L(S):=\inf_{\alpha,\beta,\nu}\Biggl\{\sum_{z\in Z(S)} \frac{\Big(\sum_{\tz\in V^S_{_\rightarrow}(z)} \beta(z,\tz)\Big)^2}{\nu(z)} + \sum_{\ttz\in Z(S)} \frac{\Big(\sum_{\tz\in V^S_{_\leftarrow}(\ttz)} \alpha(\ttz,\tz)\Big)^2}{\nu(\ttz)}\\
    & -2\sum_{(z,\ttz)\in C_{_\rightarrow}(S), d(z,\ttz)=2} \mathcal L^2(z,\ttz) \prod_{\tz\in  ]z,\ttz[}\left(\frac{\beta(z,\tz)}{L(z,\tz)\sqrt{\nu(z)}}\frac{\alpha(\ttz,\tz)}{L(\ttz,\tz)\sqrt{\nu(\ttz)}}\right)^{\ell(z,\tz,\ttz)}\Biggl\}, 
\end{align*}
where  the supremum runs over all non negative functions $\alpha$ and $\beta$ on $\X\times\X$ and $\nu$ on $\X$ satisfying  
\begin{equation}\label{condalphabetanu}
\sum_{z\in Z(S)}\sum_{\tz\in V^S_{_\rightarrow}(z)} \beta(z,\tz)=1, \sum_{\ttz\in Z(S)}\sum_{\tz\in V^S_{_\leftarrow}(z)} \alpha(\ttz,\tz)=1 \quad\mbox{and}\quad  \sum_{z\in Z(S)} \nu(z)=1.
\end{equation}
Since by the Cauchy-Schwarz inequality $B^2:=\sum_{z\in Z(S)} \frac{1}{\nu_(z)}\left(\sum_{\tz\in V^S_{_\rightarrow}(z)} \beta(z,\tz)\right)^2\geq 1$ and \\$A^2:=\sum_{\ttz\in Z(S)} \frac1{\nu(\ttz)}\left(\sum_{\tz\in V^S_{_\leftarrow}(\ttz)} \alpha(\ttz,\tz)\right)^2\geq 1$,
the quantity $R_1(S)$ is controlled by  $K(S)$ as follows 
\begin{multline}\label{vac}
R_1(S)\geq \inf_{A,B\geq 1} \big\{B^2+A^2-2 AB\, K(S)\big\}\geq \inf_{A,B\geq 1} \big\{(A^2+B^2)(1- K(S))\big\}=2(1- K(S)).
\end{multline}
As a consequence, setting 
\begin{eqnarray}\label{defr_1}
r_1=r_1^L:=\inf_SR_1(S),
\end{eqnarray}
where the infimum is over all $d$-cyclically monotone subsets $S$, the assumption $K\leq 1$ implies $r_1\geq 0$.

\begin{theorem}\label{thmprincbis} Let  $(\X,d,m,L)$ be a graph space. 
\begin{enumerate}[label=(\roman*)]
\item Assume that  $K\leq 1$. Then the $\widetilde{T}$-entropic curvature $\widetilde{\kappa}$ of  $(\X,d,m,L)$ is   bounded from below by \[\widetilde{r}= \widetilde{r}^L:=1-K^2=1-e^{-r}.\]
For the choice $\widetilde{T}:=\widetilde T_{_\rightarrow}+\widetilde T_{_\leftarrow}$, instead of $\widetilde{T}:=\max(\widetilde T_{_\rightarrow},\widetilde T_{_\leftarrow})$, the $\widetilde{T}$-entropic curvature  of  $(\X,d,m,L)$ is   bounded from below by $1-K$.

\item \label{(ii)thmprincbis} If $K\leq 1$ then the $W_1$-entropic curvature $\kappa_1$  of  $(\X,d,m,L)$ is   bounded from below by $r_1$ defined by \eqref{defr_1}.
Moreover if for any $z\in \X$, $K(z, S_2(z))<1$ then for any $d$-cyclically monotone subset $S$, one has
\begin{equation}\label{R_1R_1K}
R_1(S)\geq 4 \,\bigg(\sup_{z\in Z(S)}  \Big\{ \frac{\1_{V^S_{_\rightarrow}(z)\neq \emptyset}   }{1-K(z,\V^S_{_\rightarrow}(z))}+\frac{\1_{V^S_{_\leftarrow}(z)\neq \emptyset}}{1-K(z,\V^S_{_\leftarrow}(z))} \Big\}\bigg)^{-1}.
\end{equation}

\item For any $t\in(0,1)$ and any integer $d$ let
\begin{equation}\label{defvt}
  u_t(d)
:=\frac{d(d-1)}2\Big[\1_{d=2}+\1_{d=3}\Big]+\frac{\1_{d\geq 4}}2 \Big[d(d-1) +
\sum_{k=2}^{d-2} 
\sqrt{k(k-1)(d-k)(d-k-1)}\frac{\rho_t^d(k)}{t(1-t)}\Big] 
\end{equation}
and 
\begin{equation}\label{defcout2}
\overline{c}_t(d):=\int_0^1 u_s(d)\,q_t(s)\,ds,
\end{equation}
where $q_t$ is the kernel on $[0,1]$ defined by \[q_t(s)=\frac{2s}t \1_{[0,t]}(s)+ \frac{2(1-s)}{1-t} \1_{[t,1]}(s),\qquad s\in[0,1].\] 
Assume that   $K(z, S_2(z))<1$ or equivalently $r(z)>0$ for any $z\in \X$. Then, for $\overline{c}=(\overline{c}_t)_{t\in(0,1)}$ the $T_{\overline{c}}$-entropic curvature $\kappa_{\overline{c}}$ of  $(\X,d,m,L)$ is lower bounded by $ \overline{r}$, where
 $\overline{r}=\overline{r}^L:=\inf_{z\in \X} \overline{r}(z)$ 
 with, 
\begin{equation*}
\overline{r}(z):= 4 \bigg(\sup_{W_+,W_-} \Big\{ \frac{\1_{W_+\neq \emptyset}   }{-\log K(z,W_+)}+\frac{\1_{W_-\neq \emptyset}}{- \log K(z,W_-)} \Big\}\bigg)^{-1},
\end{equation*}
where the supremum runs over all  subsets $W_-,W_+$ of $S_2(z)$ such that for all $(w_-,w_+)\in W_-\times W_+$, $d(w_-,w_+)=4$.
\end{enumerate}
\end{theorem}
The proof of this result is postponed in Appendix B.

{\bf Comments:}
\begin{enumerate}[label=(\roman*)]
\item Since $\kappa_1\geq \widetilde{\kappa}$, the first result of this theorem implies $\kappa_1\geq 1-K^2$. Actually $r_1$ is a better lower bound for $\kappa_1$. Indeed, the inequality \eqref{vac} implies that for $K\leq 1$,    
\[r_1\geq 2(1-K)\geq 1-K^2.\]

As exposed in the next section, for the discrete hypercube $\{0,1\}^n$, 
$r_1$ is an asymptotically optimal lower bound in $n$ for the $W_1$-entropic curvature.
\item\label{obsr_1} Note that for any $d$-cyclically monotone subset $S$ and any $z\in Z(S)$, the subsets $V^S_{_\leftarrow}(z)\times V^S_{_\rightarrow}(z)$ and  $\V^S_{_\leftarrow}(z)\times\V^S_{_\rightarrow}(z)$ are also $d$-cyclically monotone. As a consequence, the inequality \eqref{R_1R_1K} ensures that 
\[\kappa_1\geq r_1\geq\inf_{z\in \X} r_1(z),\]
with
\[r_1(z):=4 \,\bigg(\sup_{V_-,V_+,W_-, W_+}  \Big\{ \frac{\1_{V_+\neq \emptyset}   }{1- K(z,W_+)}+\frac{\1_{V_-\neq \emptyset}}{1- K(z,W_-)} \Big\}\bigg)^{-1},\]
where the supremum runs over all  subsets $W_-, W_+$ of $S_2(z)$ and all subsets $V_-, V_+$ of $S_1(z)$, with $]z,W_-[\subset V_-$, $]z,W_+[\subset V_+$, and for all $(v_-,v_+)\in V_-\times V_+$ and all $(w_-,w_+)\in W_-\times W_+$
\begin{equation}\label{condVW}
d(v_-,v_+)=2 \quad\mbox{and}\quad d(w_-,w_+)=4.
\end{equation}
The quantity $r_1(z)$  only depends on the values of the jump rates on the ball of radius 2 centered at $z$, it  can be interpreted as a local lower-bound on the $W_1$-entropic curvature. 
Up to constant, $r_1(z)$ is comparable to the local lower bound on entropic curvature $r(z)$. Namely, the monotonicity  property \eqref{enfin} gives  
\begin{multline*}2(1-e^{-r(z)/2})=2\big(1-K(z,S_2(z))\big) =4 \left(\frac1{1-K(z,S_2(z))}+ \frac1{1-K(z,S_2(z))}\right)^{-1}\\
\leq  r_1(z)\leq 4\big(1-K(z,S_2(z))\big)=4(1-e^{-r(z)/2})\leq 2 r(z).
\end{multline*}
For example, on the discrete hypercube $\{0,1\}^n$ equipped with the counting measure $m_0$, we show in Section \ref{sectionhypecube} that for any $z\in \{0,1\}^n$, $r_1(z)=4\big(1-K(z,S_2(z))\big)=4/n$ whereas $r(z)=-2\log(1-1/n)\sim_{n\to +\infty} 2/n$,  which shows that the last inequality is optimal.  
In practice, the lower bound $\inf_{z\in \X} r_1(z)$ on the $W_1$-entropic curvature is more easy to handle than $r_1$ since we don't need to specify the structure of  $d$-cyclically monotone subsets. 
\item About the last result of Theorem \ref{thmprincbis}, let us first give estimates of the family of cost functions $\overline{c}=(\overline{c}_t)_{t\in(0,1)}$. 
Since
\[\sum_{k=2}^{d-2} 
{k(d-k)}\frac{\rho_t^d(k)}{t(1-t)}= \sum_{k=0}^{d} 
{k(d-k)}\frac{\rho_t^d(k)}{t(1-t)}-(d-1)\big(t^{d-2}+(1-t)^{d-2}\big)\leq d(d-1),\]
one easily checks that for any integer $d$,
\[u_t(d)\leq d(d-1)\quad \mbox{and therefore}\quad \overline{c}_t(d)\leq d(d-1).\]
Easy computations also give 
\[\sum_{k=2}^{d-2} 
{(k-1)(d-k-1)}\frac{\rho_t^d(k)}{t(1-t)}=d(d-1)-(d-1)\frac {1-t^{d}-(1-t)^{d}}{t(1-t)},\]
that implies 
\[u_t(d)\geq \frac{d(d-1)}2\Big[\1_{d=2}+\1_{d=3}\Big]+\1_{d\geq 4} \Big[d(d-1)-(d-1)\,\frac {1-t^{d}-(1-t)^{d}}{2t(1-t)} \Big]\geq \frac{d(d-1)}2.\]  
As a consequence for $d\geq 4$, one has
\[\overline{c}_t(d)\geq d(d-1)-\frac{d-1}2\int_0^1\gamma_s(d)\, q_t(s) ds,\]
with for $d\geq 2$
\[\gamma_s(d):=\frac {1-s^{d}-(1-s)^{d}}{s(1-s)}=\sum_{k=0}^{d-2} \big((1-s)^k+s^k\big).\]
Since 
\begin{align}\label{subtile}
\nonumber\int_0^1\gamma_s(d)\, q_t(s)\, ds&=2\sum_{k=0}^{d-2}\Big(\frac1{k+1}-\frac1d\Big)\big(t^k+(1-t)^k\big)\\
&\nonumber\leq 2\sum_{k=0}^{d-2} \frac1{k+1}=2\1_{d=2}+ 3\1_{d=3} + \1_{d\geq 4} \Big(3+2\sum_{k=2}^{d-2} \frac1{k+1}\Big)\\
&\leq 2\1_{d=2}+ 3\1_{d=3} + \1_{d\geq 4} \Big(3+2\log\frac{d-2}2\Big),
\end{align}
we get for any $t\in(0,1)$, $\overline{c}_t\geq \overline{c}_*$ where the cost function $\overline{c}_*$ is given by
\begin{equation*}
  \overline{c}_*(d):=\frac{d(d-1)}2\1_{d<4}+\1_{d\geq 4} \Big[d(d-1)-(d-1)\Big(\frac{3}{2} +\log\frac{d-2}{2}\Big)\Big], \quad d\in \N 
\end{equation*}
As a main property for the cost $\overline{c}_*$, for large values of $d$, $\overline{c}_*(d)$ is of order $d^2$, and for any $d\in\N$, $\overline{c}_*(d)\geq \frac{d(d-1)}2$ (which follows from the inequality $\log u\leq u-1, u>0$).

\item 
The quantity $\overline{r}(z)$   can be interpreted as a local lower-bound on the $T_{\overline c}$-entropic curvature. 
As in the last remark, the monotonicity property \eqref{enfin} also gives 
\[0<r(z) \leq  \overline{r}(z)\leq 2r(z).\]
Actually if $K\big(z,S_2(z)\big)$ is close to 1, $r_1(z)$ and $\overline{r}(z)$ are of the same order since  $\log u\sim 1-u$ for $u$ close to 1. Namely, the inequality  $\log u\leq 1-u$ for $u\geq 0$ implies $\overline r(z)\geq r_1(z)$. This inequality can be improved by using the concavity of the function $g:x\in]1,+\infty]\mapsto \Big[-\log(1-1/x)\Big]^{-1}$:
\[\overline{r}(z)\geq -2\log\Big(1-\frac{r_1(z)}2\Big)\geq r_1(z).\] 
 
 For example  on the discrete hypercube $\{0,1\}^n$ equipped with the counting measure, one has for any $z\in\{0,1\}^n$, $\overline r(z)\geq 4/n=r_1(z)$. We know from \cite{Sam21} that this lower-bound $4/n$ on the $T_{\overline c}$-entropic curvature is asymptotically optimal in $n$. Indeed, one may recover the optimal $T_2$-transport entropy inequality for the standard Gaussian measure on $\R$ from the transport entropy inequality with cost $T_{\overline{c}}$ derived from this entropic lower bound (see  \cite[Lemma 4.1]{Sam21}). 
\end{enumerate}
 
As a consequence of  Theorem 2.1 of \cite{Sam21}, a first straight forward application of Theorem \ref{thmprinc} and Theorem \ref{thmprincbis} is  the following curved Prékopa-Leindler type of inequality on discrete spaces. This functional inequality can be interpreted as the dual functional expression of the displacement convexity property \eqref{deplacebis}.  
\begin{theorem}\label{PL} Let $(\mathcal{X},d,m,L)$ be a graph space. Given a family of  cost functions $c=(c_t)_{t\in (0,1)}$, $c_t:\N\to \R$, assume that the  $T_c$-entropic curvature $\kappa_c$ of the discrete space $(\X,d,m,L)$ is bounded from below, $\kappa_c>-\infty$.  
 If $f,g,h$ are real functions on $\mathcal{X}$ satisfying for all $x,y\in \X$,
\begin{equation*}
(1-t)f(x)+tg(y)\leq \int h\, d \nu_t^{x,y}+\frac{\kappa_c}2\,t(1-t)\,c_t\big(d(x,y)\big),
\end{equation*}
then 
\begin{equation*}
\Big( \int e^{f} dm\Big)^{1-t} \Big(\int e^{g}dm \Big)^{t}\leq \int e^{h} dm \hspace{0.1cm} . 
\end{equation*}
\end{theorem}
\noindent According to Theorem \ref{thmprinc}, this result applies replacing  $\kappa_{c}$ by $r=-2\log K$, and with $c_t(d)=d(d-1), d\in \N$ for any $t\in(0,1)$.
Let us note that $\kappa_{c}$ does not need to be positive. According to Theorem \ref{thmprincbis}, it also applies if $r(z)> 0$ for all $z\in \X$ replacing $\kappa_{c}$ by $\overline{r}$ and with the family of cost functions $\overline{c}=(\overline{c}_t)_{t\in (0,1)}$ given by \eqref{defcout2}.  
 
 As for other notions of discrete curvature such as the coarse Ricci curvature \cite[Proposition 23]{Oll09}, Lin-Lu-Yau curvature \cite[Theorem 4.1]{LLY11}, the Bakry-\'Emery curvature-dimension conditions
 \cite[Theorem 6,3]{FS18}-\cite[Theorem 2.1]{LMP18}, and the entropic curvature by Erbar-Maas \cite[Theorem 1.3]{Kam20}, we easily prove a Bonnet-Myers type of theorem given next (the proof is given in Appendix B). It ensures that if the graph space  has positive entropic curvature and finite maximal degree, then its set of vertices $\X$ is finite under  bounded assumptions  on the measure $m$. 

\begin{theorem}\label{BM} Let $G=(\X,d,m,L)$ be a graph space. We assume that $\Delta(G)\neq +\infty$  and that the measure $m$ is bounded and bounded away from 0 :
\begin{equation*}\label{condm}
\sup_{x\in \X} m(x)<\infty,\qquad \inf_{x\in \X} m(x)>0.
\end{equation*}
If the entropic curvature $\kappa$ of  $(\mathcal{X},d,m,L)$ is positive  then the diameter of the space $\X$ is bounded and therefore $\X$ is finite. More precisely one has  
\begin{equation*}
\textnormal{Diam}(\mathcal{X})\leq\frac{8\log \Big(\Delta(G)\, \frac{\sup_{x\in\X} m(x)}{\inf_{x\in \X}m(x)}\Big)}{\kappa}+1,
\end{equation*}
where $\textnormal{Diam}(\mathcal{X}):=\sup_{x,y} d(x,y)$.
Same type of results hold if the $\widetilde{T}$-entropic curvature or the $W_1$-entropic curvature of the space is positive.
\end{theorem}

Assume that entropic curvature is positive ($\kappa,\kappa_1,\widetilde\kappa,\kappa_{\overline{c}}>0$) and  that $m(\X)<\infty$. Let $\mu:=m/m(\X)$ be the renormalized probability measure and let us  define the optimal transport costs
\[\widetilde{T}(\nu_0,\nu_1):=\inf_{\pi\in \Pi(\nu_0,\nu_1)} \widetilde{T}(\pi),\]
and  similarly $T_2(\nu_0,\nu_1)$ and $T_{\overline{c}_*}(\nu_0,\nu_1)$.
The displacement convexity property \eqref{deplacebis} holds with  $C_t(\widehat\pi)\geq \max\big(\kappa_1\,W_1^2(\nu_0,\nu_1), \widetilde\kappa \,\widetilde{T}(\widehat\pi) , \kappa\, {T}_2(\widehat\pi),\kappa_{\overline{c}}\, {T}_{\overline{c}_t}(\widehat\pi)\big)$ and recall that $\overline{c}_t\geq \overline{c}_*$ for any $t\in(0,1)$. Therefore optimizing over all $t\in(0,1)$ in the right-hand side the inequality \eqref{transport0} provides the following result.
\begin{corollary}\label{Transport} Let   $(\X,d,m,L)$ be a graph space with positive curvature, namely  $\kappa>0$ or $ \kappa_1>0$ or $\widetilde\kappa>0$ or $\kappa_{\overline{c}}>0$. If $m(\X)<\infty$ 
then the probability measure $\mu=m/m(\X)$ satisfies the following transport-entropy inequality, for any probability measures $\nu_0,\nu_1\in \mathcal{P}_b(\X)$  
\begin{equation*}
\frac{1}{2} \max\Big(\kappa_1\,W_1^2(\nu_0,\nu_1),\widetilde\kappa\,\widetilde{T}(\nu_0,\nu_1)  , \kappa\,T_2(\nu_0,\nu_1),\kappa_{\overline{c}}\, T_{\overline{c}_*}(\nu_0,\nu_1)   \Big) \leq \Big(\sqrt{\HR(\nu_0|\mu)}+\sqrt{\HR(\nu_1|\mu)}\Big)^2.
\end{equation*}
\end{corollary}

\begin{remark} The above transport-entropy inequalities  provides bounds on the diameter $\textnormal{Diam}(\mathcal{X})$ of the space $\X$ when the probability measure $\mu$ is bounded away from 0.
Choosing Dirac measures for $\nu_{0}$ and $\nu_{1}$, one gets 
\[\textnormal{Diam}(\mathcal{X})\leq \sqrt{-\,\frac{8\log \inf_{x\in \X} \mu(x) }{\kappa_{1}}}.\]
In some cases this upper bound is very accurate. For example, for the $n$-dimensional hypercube $\{0,1\}^{n}$ for which $\textnormal{Diam}(\mathcal{X})=n$, endowed with the uniform probability measure $\mu(x)=1/2^n$ for all $x\in\X$, the  right hand side of this inequality is $n \sqrt{2\ln(2)}$ since $\kappa_{1}\geq\frac{4}{n}$ (as we show in Section \ref{sectionhypecube}).
\end{remark}

As the space has positive $\widetilde{T}$-entropic curvature $\widetilde{\kappa}$, following the ideas of the seminal work \cite{GRST14}, the probability measure $\mu$ also satisfies a modified logarithmic-Sobolev inequality, and therefore a discrete Poincaré inequality. 
Recall that for any positive function $f$ on $\X$,  
 the entropy of $f$ with respect to $\mu$ is given by
\[{\rm Ent}_\mu(f):=\mu( f\log f )-\mu(f)\log \mu(f)=\HR(\nu|\mu),\]
where $\mu(f):=\int f d\mu$ and $\nu$ is the probability measure on $\X$ with density $f/\mu(f)$ with respect to $\mu$.  The variance with respect to $\mu$ of any function  $g:\X\to \R$  is 
\[{\rm Var}_\mu(g)=\mu(g^2)-\mu(g)^2.\]
\begin{theorem}\label{Logsob}
Let $(\X,d,m,L)$ be a graph space, with positive ${\widetilde T}$-entropic curvature $\widetilde{\kappa}$ and such that $m(\X)<\infty$.  
Then the probability measure $\mu:=m/m(\X)$ satisfies the following modified logarithmic-Sobolev inequality, for any bounded function $f:\X\to [0,+\infty)$,
\begin{equation}\label{logsob}
{\rm Ent}_\mu(f)\leq \frac{1}{2\widetilde{\kappa}} \int \sup_{x', x'\sim x}\left[\log f(x)-\log f(x')\right]_+^2 f(x)\,d\mu(x),
\end{equation}
where $[a]_+=\max(0,a)$, $a\in\R$.
It follows that  the probability measure $\mu$ also satisfies  the following Poincaré type of inequality, for any real bounded  function $g:\X\to \R$, 
\[{\rm Var}_\mu(g)\leq \frac{1}{\widetilde{\kappa}} 
\int \sup_{x', x'\sim x}\left[g(x)-g(x')\right]_+^2 d\mu(x)
.\]
\end{theorem}
\noindent The proof of this result is given in Appendix B.      According to Theorem \ref{thmprincbis}, this result applies as soon as the space $(\X,d,m,L)$ satisfies $r>0$ since $\widetilde{\kappa}\geq \widetilde r=1-e^{-r}>0$.

{\bf Comments:} Following the work \cite{BHT00}, let $\lambda_\infty$ and $\lambda_2$ be the optimal constants in the following Poincaré-type of inequalities 
\[\lambda_{\infty} \,{\rm Var}_\mu(g)\leq 
\int \sup_{x', x'\sim x}[g(x)-g(x')]^2 d\mu(x)
,\]
\begin{equation}\label{deflambda2}
\lambda_2 \,{\rm Var}_\mu(g)\leq 
\int \sum_{x', x'\sim x}[g(x)-g(x')]^2 d\mu(x),
\end{equation}
where  $g:\X\to \R$ is arbitrary. 
Note that $\lambda_2$ corresponds to the second eigenvalue of the operator $\mathcal{L}=-2 L_2$ on the $\mathbb{L}_2(\mu)$ space (recall definition  \eqref{defL1L2}), since 0 is the  first smallest  eigenvalue with eigenspace  the set of constant function on $\X$. Indeed by symmetrization, one has 
\[\int \sum_{x', x'\sim x}[g(x)-g(x')]^2 d\mu(x)= \int \sum_{x', x'\sim x}[g(x)-g(x')]^2 L_2(x,x')\,d\mu(x)=\int g \,\mathcal{L} g \,d\mu. \]

The inequality $\lambda_\infty\geq \lambda_2/\Delta(G)$ is obvious. Applying the Poincaré inequality of Theorem \ref{Logsob} with $g$ or $-g$ also provides  
$\lambda_{\infty}\geq\widetilde\kappa$. As example, for the discrete hypercube $\X=\{0,1\}^n$ equipped with any product probability measure, we prove in Section \ref{sectionhypecube} that $\widetilde\kappa\geq \frac2n\Big(1-\frac1{2n}\Big)$. It is well known that the $\lambda_2$ constant that holds for any product probability measure on $\{0,1\}^n$ is $\lambda_2=2$, and from the tensorization properties of $\lambda_\infty$ (see \cite[Introduction]{BHT00}), $\lambda_\infty=2/n$. Therefore the lower bound $\lambda_{\infty}\geq\widetilde\kappa$ is asymptotically optimal as $n$ grows to infinity in that case.

Note that $\lambda_\infty$ is related to the  following Cheeger constants (also called    \textit{inner and outer vertex expansion} of the graph if $\mu=\mu_0$), namely 
\begin{equation*}
    g_{in}:=\inf_{A\subset \X, \mu(A)\leq 1/2} \frac{\mu(\delta_{in} A)}{\mu(A)}\quad\mbox{and}\quad g_{out}:=\inf_{A\subset \X, \mu(A)\leq 1/2} \frac{\mu(\delta_{out} A)}{\mu(A)},
\end{equation*}
where $\delta_{in} A$ 
and $\delta_{out} A$ 
denotes respectively {\it the inner and the outer  vertex boundary} of the subset $A$ defined as 
\begin{equation*}
\delta_{in} A:=\big\{x\in  A\,\big|\, \exists x'\in \X\setminus A, x'\sim x\big\}\quad\mbox{and}\quad \delta_{out} A:=\big\{x\in  \X\setminus A\,\big|\, \exists x'\in  A, x'\sim x\big\}.
\end{equation*}
Choosing $g=\1_A$, the Poincaré inequality of Theorem \ref{Logsob} provides $\widetilde{\kappa}\leq 2 g_{in} $. Recall that according to  Theorem 1 in  \cite{BHT00},
one has $\lambda_{\infty} \geq g_{in}^2/4$  and $\lambda_{\infty}\geq \frac12 \big(\sqrt{1+g_{out}}-1\big)^2$.

\section{Refined results  for structured graph spaces}\label{sectrefined}

As mentioned in the introduction, in this part we complement the results of the last part  when the graph has particular geometric structure. For that purpose we introduce a class of {\it structured} graphs, that contains as basic examples,  the discrete hypercube and the lattice $\Z^n$.
\begin{definition}
We say that a graph $(\X,E)$ is {\it structured} if there exists a finite set $\Sc$ of maps $\sigma:\X\to  \X$, with the following properties
\begin{enumerate}[label=(\roman*)]
    \item For any $z\in \X$ and any $\sigma\in S$, $d\big(z,\sigma(z)\big)\leq 1$.
    \item If  $z$ and $\tz$ are  neighbours in $\X$, then there exists a single $\sigma\in S$ such that $\tz=\sigma(z)$. One defines 
    \[ \Sc_z:=\big\{\sigma\in \Sc\,\big|\, \sigma(z)\sim z\big\},\]
    so that  
    \[S_1(z):=\big\{\sigma(z)\,\big|\,\sigma\in \Sc_z\big\} \quad\mbox{and}\quad|S_1(z)|=|\Sc_z|.\]
    
    \item For any $z\in \X$ and for any  $\tau\in \Sc$, setting 
    \[\Sc_z^{\tau\rightarrow \cdot}:= \big\{\sigma\in \Sc_{\tau(z)}\,\big|\,d\big(z,\sigma\tau(z)\big)=2\big\}\quad \mbox{and}\quad \Sc_z^{\cdot\rightarrow \tau}:= \big\{\sigma\in \Sc_{z}\,\big|\,d\big(z,\tau\sigma(z)\big)=2\big\},\]
    $\Sc_z^{\tau\rightarrow \cdot}$ is empty if and only if $\Sc_z^{\cdot\rightarrow \tau}$ is empty.
    \item For any $z\in \X$ and for any  $\tau\in \Sc$, if $\Sc_z^{\cdot\rightarrow \tau}\neq \emptyset$ then there exists a one to one map $\psi:\Sc_z^{\cdot\rightarrow \tau}\to \Sc_z^{\tau\rightarrow \cdot} $ such that for all $\sigma\in \Sc_z^{\cdot\rightarrow \tau}$,
    \[ \tau\sigma(z)=\psi(\sigma)\tau(z).\]
\end{enumerate}
\end{definition}
We call $\Sc$ {\it the set of moves} of the structured graph $(\X,E)$. For a better understanding, let us introduce some paradigmatic examples of structured graphs.\

{\it Example 1 (Cayley graph).}\label{ex1} 
 Let $({\mathcal G},\ast)$ be a finite group and let $\Sc$ be a subset of generators of the group (which does not contain the neutral element of the group denoted by $e$ in order to avoid loops in the consequent graph).  The group ${\mathcal G}$ and a subset $\Sc$ determine a \textit{Cayley graph} $(\X,E)$ as follows: $\X={\mathcal G}$ and $x\sim y$ for $x,y\in \X$ if and only if $y=x\ast s$ for some $s\in \Sc$. Let us consider Cayley graphs that satisfy certain conditions: $\Sc=\Sc_{e}=\Sc_{g}$ for all $g\in {\mathcal G}$ and $\Sc$ is conjugacy stable which means that  for all $s,h\in \Sc$, $shs^{-1}\in \Sc$. Then $(\X,E)$ is a structured graph with set of moves $\Sc$. Indeed, the first three axioms are obviously satisfied and $\psi(h):=s\ast h\ast s^{-1}\in \Sc$ is a one to one map and satisfies for any $l\in {\mathcal G}$, $s\ast h\ast l=\psi(h)\ast s\ast l$. The next three examples can be seen as Cayley graphs that satisfy the above hypotheses.\

{\it Example 2.}\label{ex2} 
The \textit{ discrete hypercube} $\X=\{0,1\}^n$  is a structured graph 
with set of moves $\Sc:=\big\{\sigma_i\,\big|\,i\in [n]\big\}$,
where for any $i\in [n]$, $\sigma_i(z)$ is defined by flipping the $i$'s coordinate of $z\in \{0,1\}^n$. The discrete hypercube will be endowed with the \textit{Hamming distance}  : $d(x,y)=\sum_{i=1}^{n} \1_{x_{i}\neq y_{i}}$ for $x,y\in \{0,1\}^n$.\

{\it Example 3.}\label{ex3}
The \textit{lattice} $\X=\Z^n$ is a structured graph with set of moves 
$\Sc:=\{\sigma_{i+},\sigma_{i-}\,|\,i\in [n]\}$ where for any $z\in \Z^n$
$\sigma_{i+}(z)=z+e_i$, $\sigma_{i-}(z)=z-e_i$, and $(e_1,\ldots,e_n)$ is the canonical basis of $\R^n$.
The graph distance is given by 
$d(x,y):=\sum_{i=1}^{n} |x_{i}-y_{i}|$ for $x,y\in \Z^n$.

{\it Example 4 (Transposition model).}\label{ex4}
Let $S_{n}$ be the \textit{symmetric group} consisting of all bijective maps $\sigma: [n]\rightarrow [n]$. For any $z\in S_n$ and $\{i,j\} \subset [n]$, $i\neq j$, let $\sigma_{ij}(z)$ be  the neighbour of $z$  that differs from $z$ by a \textit{transposition} $(ij), \hspace{0.1cm} \sigma_{ij}(z):=z(ij).$ The graph distance between two elements of $x$ and $y$ of $S_{n}$, is the minimal number of transpositions $ \tau_1,..., \tau_k$ such that $x \tau_1\cdots \tau_k=y.$ The transposition model on $S_{n}$ is a structured graph with $\Sc:=\{\sigma_{ij}\hspace{0.1cm}|\hspace{0.1cm} \{i,j\}\subset [n]\}$.\

{\it Example 5 (Bernoulli-Laplace model).} Let $\X=\X_{m}$ be the slice of the discrete hypercube $\{0,1\}^n$ ,
$ \X_{m}:=\left\{x\in \{0,1\}^{n}\,\big| \,x_1+\ldots +x_n=m\right\}.$
For any $\{i,j\}\subset [n]$, let $\sigma_{ij}:\X_{m}\to\X_{m}$ denote the one to one functions that exchanges the value of coordinate $i$ with the one of coordinate $j$, namely for any $z=(z_1,\ldots,z_n)\in \X_{m} $ 
\begin{equation*}
\big(\sigma_{ij}(z)\big)_j:=z_i,\qquad \big(\sigma_{ij}(z)\big)_i:=z_j,
\end{equation*}
and for any $k\in[n]\setminus\{i,j\}$,
$\big(\sigma_{ij}(z)\big)_k:=z_k.$ Two vertices in $\X_{m}$ are declared neighbours if they differ by exactly two coordinates and $d(x,y):=\frac12\sum_{i=1}^{n} \1_{x_{i}\neq y_{i}}$ for $x,y\in \X_m$.
The Bernoulli-Laplace model is a structured graph with 
$\Sc:=\{\sigma_{ij}\hspace{0.1cm}|\hspace{0.1cm} \{i,j\}\subset [n]\}$.\

In the literature, there exists a notion related to structured graphs, the so-called \textit{Ricci flat graphs}. The concept of Ricci flat graphs was first introduced by Chung and Yau for the study of logarithmic Harnack inequalities on graphs in \cite{CY96} and recently revisited in \cite{CKKLP21}. These graphs generalize the Cayley graphs of Abelian groups. 
\begin{definition}[\cite{CY96}]
Let $(\X,E)$ be a $D$-regular graph. We say that $z\in \X$ is \textit{Ricci-flat} if there exist some maps $\sigma:B_{1}(z)\rightarrow \X$,  $1\leq i\leq D$, with the following properties
\begin{enumerate}[label=(\roman*)]
    \item  $\sigma_{i}(z^{\prime})\sim z^{\prime} \hspace{0.2cm} \text{for all } z^{\prime}\in B_{1}(z)$ ,
    \item $\sigma_{i}(z)\neq \sigma_{j}(z) \hspace{0.1cm} \text{if} \hspace{0.1cm} i\neq j$ ,
    \item $S_1(\sigma_{i}(z))= \sigma_{i}(S_1(z))$ for any  $i\in [D]$.
\end{enumerate}
A graph $(\X,E)$ is said to be \textit{Ricci flat} if it is Ricci flat for every $z\in \X$.
\end{definition}

In the case that the maps $\sigma_{i}:B_{1}(z)\rightarrow \X$ do not depend on a chosen vertex $z\in \X$, Ricci flat graphs are examples of structured graphs. Indeed under this condition, a Ricci flat graph is a structured graph with $\Sc:=\{\sigma_i\,|\, i\in[D]\}$. Given $\sigma_i\in \Sc$,
by the third property of Ricci flat graphs it follows that $\sigma_{i}\sigma_{j}(z)\in S_{1}(\sigma_{i}(z))\cap S_{1}(\sigma_{j}(z))$ and thus it is immediate that there exists a one to one map $\phi:[D]\to [D]$ such that $\sigma_i\sigma_j(z)=\sigma_{\phi(j)}\sigma_i(z)$. Since $d\big(z,\sigma_i\sigma_j(z)\big)=2$ if and only if $d\big(z,\sigma_{\phi(j)}\sigma_i(z)\big)=2$, it follows that the map $\psi:\sigma_j\to \sigma_{\phi_j}$ is one to one from $\Sc_z^{\cdot\rightarrow \sigma_i}$ to $\Sc_z^{\sigma_i\rightarrow \cdot}$, and for any $\sigma_j\in \Sc_z^{\cdot\rightarrow \sigma_i}$, $\sigma_i\sigma_j(z)=\psi(\sigma_j)\sigma_i(z)$.

Ricci flat graphs have non negative 
\textit{Bakry-\'Emery} curvature as well as non negative \textit{Ollivier curvature} or \textit{Lin-Lu-Yau} curvature \cite{KKRT16,CKKLP21}.  Structured graphs for which moves commute also have non negative Bakry-\'Emery curvature (see Proposition \ref{BE} in Appendix A, whose proof is given for completeness and which is an easy adaptation of \cite{CY96,YS10} revisited in \cite{CKKLP21}). As regards to the Erbar-Maas entropic curvature, note that similar conditions on the graph structure are given in \cite[Proposition 5.4]{EM12} to ensure non-negative Erbar-Maas entropic curvature. The next theorem asserts that structured graphs also have non-negative  entropic curvature (as defined in this paper).

\begin{theorem}\label{thmstructure}
Let $(\X,E)$ be a structured graph associated with  finite set of moves $\Sc$. The lower bound $r=r_0$ of the entropic curvature $\kappa$ of the space $(\X,d,m_0,L_0)$ given by Theorem \ref{thmprinc} is non-negative. 

Moreover, given $z\in \X$, if for any $\sigma\in \Sc$, $d\big(z,\sigma\sigma(z)\big)\leq 1$ then 
\[r_0(z)=-2\log K_0\big(z,S_2(z)\big)\geq -2\log\big(1-1/|S_1(z)|\big)>\frac 2{|S_1(z)|},\]
and therefore  $\kappa\geq \min_{z\in\X} r_0(z)\geq \frac 2{\max_{z\in \X} |S_1(z)|}\geq \frac 1{|\Sc|} $.
 \end{theorem}
The proof of this general result is postponed in Appendix B.
\newpage

\textit{Example 1  (Cayley graph).}
Since the Cayley graph of the  group $({\mathcal G},\ast)$, with conjugacy stability of the sets moves $\Sc$, is a structured graph, it has non negative entropic curvature. Moreover, if $s=s^{-1} $ for all $s\in \Sc$ then $d\big(z,s\ast s \ast z\big)=0$ for all $z\in {\mathcal G}$ and therefore  
\begin{equation*}
r_0(z)>\frac{2}{|S_{1}(z)|}=\frac{2}{|\Sc|} \hspace{0.1cm} .
\end{equation*}

{\it Example 4 (Transposition model).}
 For all $z\in S_{n}$ let us note that 
\[S_1(z)=\Big\{\sigma_{ij}(z)\,\Big |\, \{i,j\}\in I\Big \}\quad \mbox{with}\quad I=\Big\{\{i,j\}\,\Big|\,1\leq i<j\leq n\Big\}.\] 
Also, $d\big(z,\sigma_{ij}\sigma_{ij}(z)\big)=0$ for all $\{i,j\}\in I$ . Thus,
\begin{equation*}
r_0(z)>\frac{2}{|S_{1}(z)|}=\frac{4}{n(n-1)} \hspace{0.1cm}.
\end{equation*}

{\it Example 5 (Bernoulli-Laplace model).} 
For all $z\in \X_{m}$, denoting $J_0(z):=\{i\in [n]\,|\, z_i=0\}$ and $J_1(z):=\{i\in [n]\,|\, z_i=1\}$, one has 
$S_1(z)=\Big\{\sigma_{ij}(z)\,\Big |\, i\in J_0(z), j\in J_1(z)\Big \}$. Moreover, $d\big(z,\sigma_{ij}\sigma_{ij}(z)\big)=0$ for all $\{i,j\}\in [n]$ . Thus,
\begin{equation*}
r_0(z)>\frac{2}{|S_{1}(z)|}=\frac{2}{m(n-m)} \hspace{0.1cm} .
\end{equation*}

For structured graphs, one  introduces another type of transportation cost $\widetilde T_2$  comparable to $\widetilde{T}$,  related to refined modified logarithmic Sobolev inequalities, as for the cost $\widetilde{T}$ in Theorem~\ref{Logsob}.

Given two probability measures $\nu_0,\nu_1\in \Pc(\X)$, for any coupling measure $\pi\in \Pi(\nu_0,\nu_1)$, let us define 
\[\widetilde T_2( \pi):=\int\sum_{\sigma\in \Sc_x} \Pi^\sigma_\rightarrow(x)^2d\nu_0(x)+\int\sum_{\sigma\in \Sc_y} \Pi^\sigma_\leftarrow(y)^2d\nu_1(y),
\]
with
\begin{align*}
 \Pi^\sigma_\rightarrow(x)&:=\int \1_{\sigma(x)\in ]x,y]} \,d(x,y)\, r(x,\sigma(x),\sigma(x),y)\,d\pi_{_\rightarrow}(y|x),\\ 
 \Pi^\sigma_\leftarrow(y)&:=\int \1_{\sigma(y)\in ]y,x]} \,d(x,y)\, r(y,\sigma(y),\sigma(y),x)\,d\pi_{_\leftarrow}(x|y).
\end{align*} 

As an example, on the discrete hypercube $\X=\{0,1\}^n$, since $L_{0}^{d(x,y)}(x,y)=d(x,y)!$ and $\sigma_i(x)\in ]x,y]$ if and only if $x_i\neq y_i$ for any $x,y\in \{0,1\}^n$,  a simple expression holds for the cost  $\widetilde T_2(\pi)$ in that case, namely 
\begin{equation}\label{Pisigma}
\Pi^{\sigma_{i}}_{\rightarrow}(x)=\int \1_{x_{i}\neq y_{i}}d\pi_{_\rightarrow}(y|x)  \quad \mbox{and}\quad \Pi^{\sigma_{i}}_{\leftarrow}(x)=\int \1_{x_{i}\neq y_{i}}d\pi_{_\leftarrow}(x|y). 
\end{equation} 
Such a type of weak transport cost has been first introduced by Marton \cite{Mar96} to get refined concentration properties on bounded spaces related to the one reached by Talagrand with the so-called convex-hull method (see \cite[Section 4]{Tal95}). Actually, the definition of $\widetilde T_2$ on any structure graph can be interpreted as an extension  of the transportation costs introduced by Marton and  Talagrand on the hypercube. These costs belong to a larger class of costs named {\it weak transport costs} introduced in the paper \cite{GRST14}.

Observing that $\sum_{\sigma\in S_x}r(x,\sigma(x),\sigma(x),y)=1$, by the Cauchy-Schwarz inequality, one has 
\begin{align*}
   2 \widetilde{T}( \pi)&\geq \widetilde{T}_{_\rightarrow}(\widehat \pi)+\widetilde{T}_{_\leftarrow}(\widehat \pi)\geq  \widetilde T_2( \widehat\pi)\\&\geq \int\frac1{|\Sc_x|}\left( \int d(x,y)\, \sum_{\sigma\in \Sc_x}r(x,\sigma(x),\sigma(x),y)\,d\pi_{_\rightarrow}(y|x)\right)^2d\nu_0(x)\\
&\qquad+\int\frac1{|\Sc_y|}  \left(\int d(x,y)\,\sum_{\sigma\in \Sc_y} r(y,\sigma(y),\sigma(y),x)\,d\pi_{_\leftarrow}(x|y)\right)^2d\nu_1(y)
\\&\geq \frac{\widetilde{T}_{_\rightarrow}(\widehat \pi)+\widetilde{T}_{_\leftarrow}(\widehat \pi)}{\sup_{x\in \X} |\Sc_x|}\geq \frac{\widetilde{T}(\widehat \pi)}{\sup_{x\in \X} |\Sc_x|}\geq \frac{\widetilde{T}(\widehat \pi)}{|\Sc|}.
\end{align*}
By definition, let us call {\it $\widetilde T_2$-entropic curvature} of the discrete space $(\X,d,m,L)$ the best constant $\widetilde\kappa_2\in \R$ such that \eqref{deplacebis} holds with
$C_t = \widetilde\kappa_2 \,\widetilde T_2.$
As a consequence of the last inequality, if $\widetilde \kappa_2\geq 0 $ or $\widetilde \kappa\geq 0 $ then $2 \widetilde \kappa_2\geq \widetilde \kappa\geq \widetilde \kappa_2/\sup_x |\Sc_x|$.

For a better lower-estimate of $\widetilde \kappa_2$, 
one introduces a new  quantity denoted by  $\widetilde{K}(z,W)$ defined for any $z\in \X$ and $W\subset S_2(z)$. Namely let $\widetilde{K}(z,\emptyset):=0$ and for $W\neq \emptyset$ , let 
\begin{align}\label{defKtilde}
\nonumber\widetilde{K}(z,W)&=\widetilde{K}_L(z,W):=\sup  \Biggl\{  \sum_{\ttz\in W} L^2(z,\ttz) \prod_{\tz \in ]z,\ttz[}\left(\frac{\beta(\tz)}{\big(L(z,\tz)\big)^2}\right)^{\ell(z,\tz,\ttz)}
\\&
\qquad-\sum_{(\tz,w')\in ]z,W[^2, \tz\neq w'} \sqrt{\beta(\tz)}\sqrt{\beta(w')}\, \Bigg|\,{\beta}:]z,W[\to \mathbb{R}_{+},\sum_{\tz\in ]z,W[} \beta(\tz)=1 \Biggr\}.
\end{align}
For a structured graph, this quantity can also be  expressed as follows,
\begin{align*}
\widetilde{K}_L(z,W)&=\widetilde{K}(z,W):=\sup  \Biggl\{  \sum_{\ttz\in W} L^2(z,\ttz) \prod_{\sigma\in \Sc_{]z,\ttz[}}\left(\frac{\beta(\sigma)}{\big(L(z,\sigma(z))\big)^2}\right)^{\ell(z,\tz,\ttz)}\!\!\!
\\&\nonumber\qquad
-\sum_{(\sigma,\tau)\in \Sc_{]z,W[}^2, \sigma\neq \tau} \sqrt{\beta(\sigma)}\sqrt{\beta(\tau)}\, \Bigg|\,{\beta}: \Sc_{]z,W[}\to \mathbb{R}_{+},\sum_{\sigma\in \Sc_{]z,W[}} \beta(\sigma)=1 \Biggr\},
\end{align*}
where for any subset $W\subset S_2(z)$, $\Sc_{]z,W[}:=\big\{\sigma\in S\,\big|\, \sigma(z)\in ]z,W[\big\}$. If $m=m_0$ and $L=L_0$, then we write $ \widetilde K_{0}(z,W):=\widetilde K_{L_0}(z,W)$.
\begin{theorem}\label{Thmstructure}
Let $(\X,d,m,L)$ be a graph space such that $(\X,E)$ is a structured graph  with set of moves $\Sc$.  For any $z\in \X$, let us define $\widetilde r_2=\widetilde r_2^L:=\inf_{z\in \X}\widetilde r_2(z) $, with
\begin{equation}\label{defr_3}
\widetilde r_2(z)=\widetilde r_2^L(z):=1-\widetilde{K}_L(z), \quad \mbox{and}\quad  \widetilde{K}_L(z):=\sup_{W\in S_2(z)} \widetilde{K}_L(z,W).
\end{equation}
For any $z\in \X$, one has
\begin{eqnarray}\label{rr_3} 1-K(z,S_2(z))\leq  \widetilde r_2(z)\leq |S_1(z)|\,\big(1-K(z,S_2(z))\big).
\end{eqnarray}

\begin{enumerate}[label=(\roman*)]
\item If the generator $L$ satisfies for any $z\in \X$ and any $\sigma,\tau \in \Sc$   with $d\big(z,\tau\sigma(z)\big)=2$ 
\begin{equation}\label{condL}
L\big(z,\sigma(z)\big)L\big(\sigma(z),\tau\sigma(z)\big)=L\big(z,\tau(z)\big)L\big(\tau(z),\tau\sigma(z)\big),
\end{equation}
then the
$\widetilde{T}_2$-entropic curvature $\widetilde{\kappa}_2$ of  $(\X,d,m,L)$ is   bounded from below by $\widetilde r_2\geq 0$.
\item Assume that  \eqref{condL} holds and moreover that for any $z\in \X$ and any  $\sigma\in \Sc$,
\begin{equation}\label{restr}
d\big(z,\sigma\sigma(z)\big)\leq 1
\end{equation}
then the above result can be improved replacing the curvature cost $C_t(\widehat \pi)=\widetilde r_2 \widetilde T_2(\widehat\pi)$ in the $C$-displacement convexity property of entropy \eqref{deplacebis} by  the cost $C_t(\widehat \pi)=\widetilde r_2\widetilde C_t^1(\widehat \pi)$, $t\in (0,1)$, where for any $D\geq1$, the cost  $\widetilde C_t^D(\widehat \pi)$ is given by  
\begin{equation}\label{defC_t^D}
\widetilde C_t^D(\widehat \pi):=\int \sum_{\sigma\in \Sc} D^2 h_t\left(\frac{\Pi^\sigma_\rightarrow(x)}{D}\right) d\nu_0(x)+\int \sum_{\sigma\in \Sc} D^2 h_{1-t}\left(\frac{\Pi^\sigma_\leftarrow(y)}D\right) d\nu_1(y),
\end{equation}
where  for any $u\geq 0$, 
\[
h_t(u):=\frac{ th(u)-h(tu)}{t(1-t)}\quad\mbox{with} \quad h(u):=
\left\{\begin{array}{ll}
2\left[(1-u)\log(1-u)+u\right] &\mbox{ for } \;0\leq u\leq 1,
\\
+\infty &\mbox{ for } \;u>1.
\end{array}\right.
\]
Assume that $D:=\textnormal{Diam}(\mathcal{X})<\infty$. If condition \eqref{restr} is not satisfied, then  the $C$-displacement convexity property of entropy \eqref{deplacebis} also holds with the cost $C_t(\widehat \pi)=\widetilde r_2\widetilde C_t^D(\widehat \pi)$, $t\in (0,1)$.   
\end{enumerate}
\end{theorem}

The proof of this Theorem is given in Appendix B.

{\bf Comments:}
\begin{enumerate}[label=(\roman*)]

\item Condition \eqref{condL} makes sense since if $d\big(z,\tau\sigma(z)\big)=2$ then $\sigma(z)\in ]z,\tau\sigma(z)[$, and according to the definition of structured graph there exists $\psi(\sigma)\in \Sc_{\tau(z)}$ such that  
$\tau\sigma(z)=\psi(\sigma)\tau(z)$, and therefore  $\tau(z)\in ]z,\tau\sigma(z)[.$ 
 Condition \eqref{condL} actually provides needed properties for the proof of Therorem \ref{Thmstructure} which are collected in  Lemma \ref{lemL} in Appendix A.
\item The second part of the theorem improves its first part
since for any $t\in (0,1), u\geq 0, h_t(u)\geq u^2$, and therefore $\widetilde C_t^D(\widehat \pi)\geq \widetilde T_2(\widehat \pi)$. This improvement is useful in particular when considering the derived modified logarithmic Sobolev inequalities. It allows to reach smaller discrete Dirichlet forms in the right-hand side of the modified Sobolev inequality (see the comments of Theorem \ref{Logsobbis} below).
\item Condition \eqref{restr} implies that any discrete geodesic $(z_0,\ldots,z_d)$  from $z_0=z$ to any vertex $z_d\in \X$  ($d=d(z, z_d)$) does not use any move $\sigma\in \Sc$ more than one time. Indeed, if $(z_0,\ldots, z_d) $ is such that for some $0\leq k<\ell\leq d-1$, $z_{k+1}=\sigma(z_k)$ and $z_{\ell+1}=\sigma(z_\ell)$, then  Lemma \ref{lemL} implies that $(z_0,\sigma(z_0), \sigma\sigma(z_0), \ldots,\sigma  \sigma(z_\ell), z_{\ell+2},\ldots,z_d)$ is also a geodesic. Therefore $d\big(z_0, \sigma\sigma(z_0)\big)=2$ which is a contradiction with condition \eqref{restr}. 
\end{enumerate}
As for $\widetilde{T}$-entropic curvature, positive $\widetilde T_2$-entropic curvature also provides transport entropy inequalities and also  modified logarithmic-Sobolev  and Poincaré inequalities.
For any $\sigma\in \Sc$, and $g:\X\to \R$, let 
\[\partial_\sigma g(z):=g(\sigma(z))-g(z),\qquad z\in \X.\]

\begin{theorem}\label{Logsobbis}
Let $(\X,d,m,L)$ be a graph space with $m(\X)<+\infty$ and such that $(\X,E)$ is a structured graph with set of moves $\Sc$. Let $\mu:=m/m(\X)$.
\begin{enumerate}[label=(\roman*)]
 \item If the $\widetilde T_2$-entropic curvature $\widetilde \kappa_2$ of the space $(\X,d,m,L)$ is positive, then  $\mu$ satisfies the following transport-entropy inequality, for any probability measures $\nu_0$ and $\nu_1$ on $\mathcal{P}_b(\X)$
\begin{equation*}
\frac{\widetilde\kappa_2}{2} \,\widetilde{T}_2(\nu_0,\nu_1)  \leq \Big(\sqrt{\HR(\nu_0|\mu)}+\sqrt{\HR(\nu_1|\mu)}\Big)^2,
\end{equation*}
 with $\widetilde{T}_2(\nu_0,\nu_1):=\inf_{\pi\in \Pi(\nu_0,\nu_1)} \widetilde{T}_2(\pi)$. If moreover condition \eqref{restr} holds, we also have for any $\nu\in \mathcal{P}_b(\X)$,
 \begin{equation}\label{transportentC_0^D}
 \frac{\widetilde \kappa_2}{2}\inf_{\pi\in \Pi(\mu, \nu)} \widetilde C_0^D(\widehat \pi)\leq \HR(\nu|\mu),
 \end{equation}
 where the cost $\widetilde C_0^D$ is defined like in \eqref{defC_t^D}
 with   $h_0:=h$ and 
 \[h_1(u):=
\left\{\begin{array}{ll}
-2\log(1-u)-2u &\mbox{ for } \;0\leq u< 1,
\\
+\infty &\mbox{ for } \;u\geq1.
\end{array}\right.\]

 \item 
 If the $\widetilde T_2$-entropic curvature $\widetilde \kappa_2$ of the space $(\X,d,m,L)$ is positive, then  $\mu$ satisfies the following modified logarithmic-Sobolev inequality, for any bounded function $f:\X\to [0,+\infty)$, 
 \begin{equation}\label{logsobT3}
{\rm Ent}_\mu(f)\leq \frac{1}{2\widetilde \kappa_2} \int \sum_{\sigma\in \Sc}  [\partial_\sigma \log f]_-^2  f\,d\mu \hspace{0.2cm} .
\end{equation}
 \item
 If the $C$-displacement convexity property of entropy \eqref{deplacebis} holds with  the cost $C_t(\widehat \pi)= \widetilde{\kappa}_2\widetilde C_t^D(\widehat \pi)$  given in Theorem \ref{Thmstructure} for some $D\geq 1$
then for any bounded function $f:\X\to [0,+\infty)$,
\begin{equation}\label{logsobbis}
{\rm Ent}_\mu(f)\leq  \int \sum_{\sigma\in \Sc} \frac{\widetilde \kappa_2 D^2}2 \, h^*\left(\frac{2}{D\widetilde\kappa_2} [\partial_\sigma \log f]_-\right)  f\,d\mu,
\end{equation}
where $[a]_-=\max(0,-a)$, $a\in\R$ and $\frac{1}2 h^*(2v)=e^{-v}+v-1$, $v\geq 0$.
 \end{enumerate}
In any case, it follows that   $\mu$ also satisfies  the following Poincaré  inequalities, 
\begin{equation}\label{poinT3bis}
{\rm Var}_\mu(g)\leq \frac{1}{\widetilde{\kappa}_2} 
\int \sum_{\sigma\in \Sc}[\partial_\sigma g]_-^2 d\mu
,
\end{equation}
and therefore
 \begin{equation}\label{poinT3}
{\rm Var}_\mu(g)\leq \frac{1}{2\widetilde{\kappa}_2} 
\int \sum_{\sigma\in \Sc}(\partial_\sigma g)^2 d\mu
,
\end{equation}
for any real bounded function $g:\X\to \R$.
\end{theorem}
The proof of the transport entropy inequality is identical to the one of Corollary \ref{Transport}. Proofs of  modified logarithmic Sobolev inequalities  and the Poincaré inequality are given  together with the one of Theorem \ref{Logsob} in Appendix B.

{\bf Comments:}
\begin{enumerate}[label=(\roman*)]

\item According to the definition \eqref{deflambda2} of the Poincaré constant $\lambda_2$, the Poincar{\'e} inequality \eqref{poinT3} ensures that $\lambda_{2}\geq  2 \widetilde{\kappa}_{2}$. For example on the discrete hypercube equipped with the uniform probability measure $\mu=\mu_0$, we proves in Section \ref{sectionhypecube} that $\widetilde \kappa_2\geq \widetilde r_2\geq 1$. As a consequence since it is known that $\lambda_2=2$ for $\mu_0$ on the discrete hypercube,  $\kappa_2= \widetilde r_2= 1$ on this space.

For completeness, recall also that in discrete setting, when $\X$ is finite and $\mu=\mu_0=m_0/|\X|$, the Poincaré constant $\lambda_2$ is also related to the \textit{Cheeger constant} $h_G$ of the graph  defined by
\begin{equation*}
    h_G:=\min_{A\subset \X, |A|\leq |\X|/2} \frac{|\partial A|}{|A|}\,,
\end{equation*}
where  $\partial A$ denotes {\it the edge boundary} of the subset $A$ defined as 
\begin{equation*}
 \partial A:=\big\{(x,x')\,\big|\, x\in A, x'\in \X\setminus A, x'\sim x\big\}.\]
 It is a well known fact that $2h_G\geq \lambda_2\geq h_G^2/2 $ 
 (see for example \cite{Chu07}). Applying the Poincaré inequality 
 \eqref{poinT3} to the function $g=\1_A$ also provides 
 $h_G\geq 4\widetilde \kappa_2$. 
 For general probability measure $\mu$ one may introduce the  conductance $\Phi$ that generalize the above Cheeger constant, 
 \[\Phi_\mu:=\inf_{A\subset \X,\mu(A)\leq 1} \frac{\sum_{(z,z')\in \partial A } \mu(z)}{\mu(A)}.\]
 Similar connections are proved between $\lambda_2$ and $\Phi$ in \cite[section 3]{MT06}.
 The Poincaré inequality \eqref{poinT3bis} applied with $g=-\1_A$ provides 
 $\Phi_\mu\geq 2\widetilde\kappa_2$. 
 
\item Let us give few  inequalities which are useful to compare the discrete Dirichlet form on the right-hand side of \eqref{logsobT3} and \eqref{logsobbis} with other  discrete Dirichlet forms. 
Let $\alpha,a,b$ be positive real numbers with $a\geq b$, one easily proves that
\[\frac a2 h^*\big(2\alpha[\log a-\log b]\big)\leq  \alpha^2 \frac{a}2 [\log a-\log b]^2. 
\]
If $\alpha\leq 1$, then the  convexity property of the function $h^*$ implies
\[
  \frac a2 \,h^*\big(2\alpha[\log a-\log b]\big)\leq \alpha\, \frac a 2\, h^*\big(2[\log a-\log b]\big) 
  =\alpha\left( a[\log a-\log b]-[a-b]\right)
\]
and if  $\alpha\geq 1$, then  the decreasing monotonicity property of the function $u\in (0,+\infty)\to \frac1{2u^2} h^*(2u)$ gives 
\[\frac a2 h^*\big(2\alpha[\log a-\log b]\big) \leq\alpha^2\, \frac{a}2\,  h^*(2[\log a-\log b]).\]
As a consequence, for any $\alpha>0$,
\begin{align*}
\frac a 2 h^*\big(2\alpha[\log a-\log b]\big)&\leq \frac{\alpha\max(1,\alpha)}2 \,h^*(2[\log a-\log b])a \\
&\leq \alpha\max(1,\alpha) \min\left\{[\log a-\log b][a-b],\frac{[a-b]^2}{2b} \right\}.
\end{align*}
Applying these inequalities  with $\alpha=1/(D\widetilde\kappa_2)$, $a=f(x)$ and $b=f(\sigma(x))$, one gets
the following comparisons
\begin{align*}\int \sum_{\sigma\in \Sc} \frac{\widetilde \kappa_2 D^2}2 \, h^*\left(\frac{2}{D\widetilde\kappa_2} [\partial_\sigma  \log f]_-\right)  f\,d\mu
&\leq \min\left\{\frac1{2\widetilde \kappa_2} \int \sum_{\sigma\in \Sc} [\partial_\sigma \log f ]_-^2 f\,d\mu,\right.\\
&\qquad \frac 12 \max\left(D, \frac 1{\widetilde \kappa_2}\right) \int \sum_{\sigma\in \Sc} \frac{[\partial_\sigma f(x)]_-^2}{f(\sigma(x))} \,d\mu(x),\\
&\qquad \left.\max\Big(D, \frac 1{\widetilde \kappa_2}\Big) \int \sum_{\sigma\in \Sc} {[\partial_\sigma \log f]_-[\partial_\sigma f]_-} \,d\mu\right\}
\end{align*}
In particular, it follows that \eqref{logsobbis} is a refinement of \eqref{logsobT3}.
\item For $f,g: \X\to \R$ let $\mathcal{E}(f.g)$ denote the Dirichlet form defined by 
\[\mathcal{E}_L(g,f)=-\int g\,Lf\, d\mu=\frac12  \int \sum_{z', z'\sim w} (g(z')-g(z))(f(z')-f(z)) L(z,z')\, d\mu(z).\]
It is a well known fact (see \cite[section 2]{MT06}) that the Poincaré inequality 
\begin{equation}\label{PEM}
c\, {\rm Var}_\mu(f)\leq \mathcal{E}_L(f,f), \quad\mbox{for all } f:\R\to \R,
\end{equation}
and the modified logarithmic-Sobolev inequality 
\begin{equation}\label{MLSIEM}
c\, {\rm Ent}_\mu(f)\leq \mathcal{E}_L(f,\log f), \quad\mbox{for all } f:\R\to \R,
\end{equation}
for some $c>0$, respectively implies exponential decay of the variance and the entropy of $P_tf={e^{tL}}f$, namely for all $t\geq 0$
\[{\rm Var}_\mu(P_tf)\leq e^{-ct} {\rm Var}_\mu(f) \quad\mbox{and}\quad{\rm Ent}_\mu(P_tf)\leq e^{-ct} {\rm Ent}_\mu(f) .\]
Bounds for mixing times then follows for the continuous time Markov chain associated to the generator $L$ (see \cite[Corollary 2.6]{MT06}). In the entropic curvature approach by Erbar-Maas \cite{EM12}, positive entropic curvature provides modified logarithmic-Sobolev inequality of type \eqref{MLSIEM}. For us, positive $\widetilde T_2$-entropic curvature $\widetilde \kappa_2$ of the space $(\X,d,m,L)$  implies the Poincaré inequality \eqref{poinT3} that corresponds to \eqref{PEM} with generator $L=L_2$ given by \eqref{defL1L2} and $c=\widetilde \kappa_2$, since by symmetrisation 
\begin{equation}\label{grrr}
\int \sum_{\sigma\in \Sc} (\partial_\sigma f)(\partial_\sigma g)\,  d\mu = 2\,\mathcal{E}_{ L_2}(f,g).
\end{equation}

Similarly,  if  \eqref{logsobbis} holds, then by using the above Dirichlet forms comparisons,   the modified logarithmic Sobolev inequality \eqref{MLSIEM} holds with the generarator $L=L_2$ and $c:=\min(1/D, \widetilde \kappa_2)/2$ since 
\[\int \sum_{\sigma\in \Sc} [\partial_\sigma \log f]_-[\partial_\sigma f]_- \,d\mu\leq\int \sum_{\sigma\in \Sc} (\partial_\sigma \log f)(\partial_\sigma f) \,d\mu= 2\,\mathcal{E}_{ L_2}(f,\log f).\]  
However, it remains a challenge  to introduce another  $C$-displacemnent convexity property \eqref{deplacebis} along Schr\"odinger bridges at zero temperature from which one could derive Poincaré or modified  logarithmic Sobolev inequalities with any generator $L$, instead of $L_2$.

A careful reading of the proof of Theorem  \ref{Logsobbis} shows that  the modified logarithmic Sobolev inequality  \eqref{logsobT3} may actually be improved by substracting on the right hand side the quantity 
\[ \frac{\widetilde  \kappa_2}2 \mu(f) \int \sum_{\sigma\in \Sc} \left(\Pi^\sigma_\leftarrow(y)\right)^2 d\mu(y),\]
where  $\Pi^\sigma_\leftarrow(y)$ is defined with a $W_1$-optimal coupling $\widehat \pi$ with first marginal $\nu_0=\frac{f}{\mu(f)} \mu$ and second marginal $\nu_1=\mu$.
However, we do not know how to get ride of this improvement. 
\end{enumerate}

\section{Perturbation results}
\label{sectionpertur}

\subsection{Perturbation with a potential}\label{sectionpertu}
Let $(\X,d,m,L)$ be a graph space satisfying a $C$-displacement convexity property of entropy \eqref{deplacebis}. Let $m_v$ denote the measure with density $e^{-v}$ with respect to $m$, where $v:\X\to \R$ is a potential. In this part, we  analyse the perturbations of the $C$-displacement convexity property of the relative entropy  along the Schr\"odinger bridges at zero temperature of the space $(\X,d,m,L)$ when $m$ is replaced by $m_v$.

Since for any probability measure $\nu\in \Pc(\X)$ absolutely continuous with respect to $m$,
\begin{equation}\label{Hm0mv}
H(\nu|m_v)=H(\nu|m)+\int v\, d\nu,
\end{equation}
convexity properties of $t\in(0,1)\to H(\widehat \nu_t|m_v)$ may follow from convexity properties of $t\in(0,1)\to H(\widehat \nu_t|m)$ and convexity properties of $\psi:t\in(0,1)\to \int v\, d\widehat \nu_t$.
According to Lemma \ref{deriveseconde} (see Appendix A), assuming $\nu_0$ and $\nu_1$ have bounded support, one has 
\begin{equation*}
\psi''(t)=\sum_{(x,y)\in\X^2}\left(\int v\, d \nu_t^{x,y}\right)^{''}\widehat \pi(x,y)\\
=\sum_{(x,y)\in\X^2}
d(x,y)\big(d(x,y)-1\big)\, D_t v(x,y)
\widehat \pi(x,y)
\end{equation*}
with 
\begin{equation}\label{mulh}
D_t v(x,y):=\sum_{(z,\ttz)\in [x,y], d(z,\ttz)=2} Dv(z,\ttz) \, L^2(z,\ttz)\, r(x,z,\ttz,y)\, \rho_t^{d(x,y)-2}(d(x,z)),
\end{equation}
and for $z,\ttz\in \X$ with $d(z,\ttz)=2$,
\[ Dv(z,\ttz):=\sum_{\tz\in  ]z,\ttz[} \left(v(\ttz)+v(z)-2v(\tz)\right)\frac{L(z,\tz)L(\tz,\ttz)}{L^2(z,\ttz)}.\]
Observe that $Dv(z,\ttz)=Dv(\ttz,z)$ can be interpreted as a local discrete laplacian of the potential $v$ at $(z,\ttz)$.

It follows that 
\begin{align*}
\psi(t)&=(1-t)\psi(0)+t\psi(1)- \frac{t(1-t)}2 \int_0^1\psi''(s)q_t(s)\, ds\\
&=(1-t)\psi(0)+t\psi(1)- \frac{t(1-t)}2 \iint c^v_t(x,y)\,d\widehat\pi(x,y),
\end{align*}
with 
\[c^v_t(x,y):=d(x,y)\big(d(x,y)-1\big)\, \int_0^1 D_s v(x,y) q_t(s)\, ds.\]
This together with \eqref{Hm0mv} gives the following result.
\begin{theorem}\label{entropicperturb}
Let $(\X,d,m,L)$ be a graph space. Assume that a   $C$-displacement convexity property of entropy \eqref{deplacebis} holds.
Given a  potential $v:\X\to \R$, let $m_v$ denote the measure with density $e^{-v}$ with respect to $m$. 
Then the relative entropy with respect to $m_v$, $\nu\in \Pc(\X)\mapsto H(\nu|m_v)$, satisfies the $C^v$-displacement convexity property \eqref{deplacebis} along the Schr\"odinger bridge at zero temperature of the space $(\X,d,m,L)$, with for any $t\in (0,1)$,
\[C^v_t(\widehat \pi)=C_t(\widehat \pi) +\iint c_t^v(x,y)\, d\widehat \pi(x,y) .\]
\end{theorem}

Another way to  get $C^v$-displacement convexity properties with the measure $m_v$ is to consider the generator $L_v$  defined by 
\[L_v(x,y)=e^{\frac12(v(x)-v(y))} L(x,y),\qquad x,y\in \X, x\neq y.\]
One easily checks that the measure $m_v$ is reversible with respect to $L_v$ and that the space $(\X,d,m_v,L_v)$ is a graph space.
Moreover, since for any $x,y\in \X$,
\[L_v^{d(x,y)}(x,y)=e^{\frac12(v(x)-v(y))} L^{d(x,y)}(x,y),\]
the Schr\"odinger briges at zero temperature of the space $(\X,d,m_v,L_v)$ are the same as the one of the space $(\X,d,m,L)$. Indeed for any $x,y\in \X$ and $z\in [x,y]$, the quantity
\[r(x,z,z,y)=\frac{L_v^{d(x,z)}(x,z) L_v^{d(z,y)}(z,y)}{L_v^{d(x,y)}(x,y)}=\frac{L^{d(x,z)}(x,z) L^{d(z,y)}(z,y)}{L^{d(x,y)}(x,y)},\]
does not depend on the potential $v$ and therefore the Schr\"odinger briges $\nu_t^{x,y}$  between Dirac measure on the space $(\X,d,m_v,L_v)$ are the same as the one of the space $(\X,d,m,L)$. Moreover, note that if on a structured graph the generator $L$ satisfies \eqref{condL} then the generator $L_v$ also satisfies  \eqref{condL}. As a consequence, 
any result we get on the graph space $(\X,d,m_v,L_v)$ on the lower bound on entropic curvature from Theorems \ref{thmprinc}, \ref{thmprincbis} and \ref{Thmstructure} can be interpreted as a perturbation result of the same result on  $(\X,d,m,L)$. Note that if one choose $m=m_0$ and $L=L_0$, the graph space $(\X,d,m_v,L_v)$ is exactly the space $(\X,d,m_v,L_1)$.

For further use, let us just simplify the definition of the key quantities $K^v(z,W):=K_{L_v}(z,W)$ and $\widetilde K^v(z,W):=\widetilde K_{L_v}(z,W)$ on the space  $(\X,d,m_v,L_v)$, for  $z\in\X$  and $W\subset S_2(z)$.
Observing that 
\[-Dv(z,\ttz)=2 \sum_{\tz\in ]z,\ttz[} \Big(\log \frac{L_v^2(z,\ttz)}{L^2(z,\ttz)}-2\log \frac{L_v(z,\tz)}{L(z,\tz)}\Big) \, \frac{L_v(z,\tz)L_v(\tz, \ttz)}{L_v^2(z,\ttz)},\]
according to \eqref{defR_2} and to \eqref{defKtilde}, one has 
\begin{align}\label{defKLv}
    K^v\big(z,W\big)&=\sup_\alpha \Biggl\{  \sum_{\ttz\in W} e^{-Dv(z,\ttz)/2} L^2(z,\ttz) \prod_{\tz\in  ]z,\ttz[}\left(\frac{\alpha(\tz)}{L(z,\tz)}\right)^{\frac{2L(z,\tz)L(\tz,\ttz)}{L^2(z,\ttz)}}\Biggl\},
\end{align}
where the supremum is over all ${\alpha}: ]z,W[\to \mathbb{R}_{+}$  such that $\sum_{v \in ]z,W[} \alpha(v)=1$, and $\widetilde K^v(z)=\sup_{W\subset S_2(z)} \widetilde{K}^v(z,W)$  with
\begin{multline}\label{defKvtilde}
\widetilde{K}^v(z,W):=\sup  \Biggl\{ \sum_{\ttz\in W}e^{-Dv(z,\ttz)/2} L^2(z,\ttz) \prod_{\sigma\in \Sc_{]z,\ttz[}}\left(\frac{\beta(\sigma)}{L(z,\sigma(z))^2}\right)^{\frac{L(z,\sigma(z))L(\sigma(z),\ttz)}{L^2(z,\ttz)}}\\-\sum_{(\sigma,\tau)\in \Sc_{]z,W[}^2, \sigma\neq \tau} \sqrt{\beta(\sigma)}\sqrt{\beta(\tau)}\Biggl\},
\end{multline}
where the supremum is over all ${\beta}: \Sc_{]z,W[}\to \mathbb{R}_{+}$ such that $\sum_{\sigma\in \Sc_{]z,W[}} \beta(\sigma)=1$.  

\subsection{Restriction to  convex subsets}\label{sectionrestconv}
This section concerns another type of perturbation result, when the measure $m$ is restricted to a 
{\it convex} subset $\Cc$ of $\X$. The convexity property of a subset is defined as follows in this paper.
\begin{definition}\label{defconvset} On a graph space $(\X,d,m,L)$, a subset $\Cc$ of $\X$ is convex if for any $x,y\in \Cc$,
$[x,y]\subset \Cc$.
\end{definition}
Let $(\X,d,m,L)$ be a graph space. Given a subset $\Cc$ of $\X$, let $(\Cc,d_\Cc, m_\Cc,L_\Cc)$ denotes the graph space restricted to $\Cc$ defined by :  $m_\Cc=\1_\Cc m$, , $d_\Cc(x,y):=1$ if $d(x,y)=1$,  and  
$L_\Cc(x,y):=L(x,y)$ for any $x,y\in \X_\Cc$, $x\neq y$. One easily checks that the space $(\Cc,d_\Cc, m_\Cc,L_\Cc)$ is also a graph space.

If $\Cc$ is a convex subset of $\X$, then the set of discrete geodesics on $(\Cc,d_\Cc, m_\Cc,L_\Cc)$ between two vertices $x$ and $y$ of $\Cc$ is the same as the one on $(\X,d,m,L)$. Since $L_\Cc(\gamma)=L(\gamma)$  for any discrete geodesic $\gamma$ between  $x\in \Cc$ and $y\in \Cc$, it follows that the Schr\"odinger bridge at zero temperature between the Dirac measures at $x$ and  $y$ is the same on the space $(\Cc,d_\Cc, m_\Cc,L_\Cc)$ as on the space $(\Cc,d_\Cc, m_\Cc,L_\Cc)$.  As a consequence,  this observation also holds  for any Schr\"odinger bridge at zero temperature between two probability measures  on $\Cc$. This remark implies the following result.
\begin{theorem}\label{thmrestconv}
Let $(\X,d,m,L)$ be a graph space and let $\Cc$ be a convex subset of $\X$. If the relative entropy with respect to $m$ satisfies a $C$-displacement convexity property \eqref{deplacebis} on the  space $(\X,d,m,L)$, then the same property holds for the relative entropy with respect to $m_\Cc$  on the space $(\Cc,d_\Cc, m_\Cc,L_\Cc)$. 
\end{theorem}

\section{Applications} \label{section-interaction-pot}
\label{examplesZcube} 
This part is devoted to applications of the perturbation results of the last Section \ref{sectionpertu} for two specific structured graphs, as a guideline for many other structure graphs which are not presented in this paper. We only focus on the discrete hypercube $\X=\{0,1\}^n$ and  the lattice $\X=\Z^n$ endowed with a measure $m_{v}$ with density $e^{-v}$ with respect to the counting measure $m_0$ on the set of vertices.  For $z\in \X$, we analyze the constants $r^v(z)=r^{L_v}(z),r_1^v(z)=r_1^{L_v}(z),\overline{r}^v(z)=\overline{r}^{L_v}(z),$ and $\widetilde r_2^v(z)=\widetilde r_2^{L_v}(z)$ that allow to bound from below the different types of entropic curvatures of the graph space $(\X,d,m_v,L_v)$, defined in Section \ref{sectionpertu}. These results show that our approach of entropic curvature is robust on discrete spaces. Indeed, by applying Theorem \ref{PL}, Corollary \ref{Transport}, Theorem \ref{Logsob} and Theorem \ref{Logsobbis}, one derives  functional inequalities for the measure $m_v$ or its associated normalized probability measure $\mu_v:=m_v/m_v(\X)$ under weak conditions on the potential $v$, involving eigenvalues of some Hessian type of matrices for the potential $v$. 

\subsection{Ising models on the discrete hypercube}\label{sectionhypecube}
As mentioned in the introduction, the discrete hypercube is a structured graph with set of moves $\Sc:=\{\sigma_i\,|\,i\in[n]\}$ where $\sigma_i(z)$ is defined by flipping the $i$'s coordinate of $z\in\X=\{0,1\}^n$.
In this part,  $\mu_0$ is the uniform probability measure on $\{0,1\}^n$. Given $v:\{0,1\}^n\to \R$,  $m_v=e^{-v}m_0$ is a perturbation of the counting measure $m_0$ on $\{0,1\}^n$.

\begin{remark}
For $z\in \{0,1\}^n$ and $W\subset S_2(z)$ the quantities $K_0(z,W)$ and $\widetilde K_0(z,W)$ computed in this section  for the hypercube are the same for any graph whose local structure is the one of the hypercube. Therefore, the lower bounds on entropic curvature reached from these two quantities are also the same. 
As pointed out in \cite{LMP17}, the hypercube is not determined by its local structure. Indeed, Laborde and Hebbare \cite{HL82} showed that the conjecture according to which every bipartite, regular graph satisfying that all balls of radius 2 are isomorphic to those of the hypercube is necessarily the hypercube is false. 
\end{remark}

Since for any $i\in[n]$ and $z\in \X$, $\sigma_i\sigma_i(z)=z$, Theorem \ref{thmstructure} gives $K_0\big(z,S_2(z)\big)\leq 1-1/n$ for any $z\in \X$. Actually, by choosing $\alpha(\tz)=1/n$ for any $\tz \in S_1(z)$ in the definition of $K_0\big(z,S_2(z)\big)$, one exactly gets $K_0\big(z,S_1(z)\big)=1-1/n$ and therefore according to  \eqref{defr}, the entropic curvature $\kappa$ of the graph space $\big(\{0,1\}^n,d,m_0,L_0\big)$ satisfies 
\[\kappa\geq r\geq \min_{z} r(z)=-2\log(1-1/n).\]
Theorem \ref{thmprincbis} also provides the lower bound for the $\widetilde{T}$-entropic curvature of this space \[\widetilde\kappa\geq {\widetilde{r}}= 1-K^2\geq 1-\Big(1-\frac1n\Big)^2=\frac2n\Big(1-\frac1{2n}\Big).\]

Let us now compute a lower bound on $r_1(z)$, $z\in \X$, to reach a lower bound on the $W_1$-entropic curvature.
Notice that 
\[S_2(z):=\Big\{\sigma_j\sigma_i(z)\,\Big|\, (i,j)\in I\Big\}\qquad \mbox{with}\quad I=\Big\{(i,j)\,\Big|\,1\leq i<j\leq n\Big\}.\]
Given $W\subset S_2(z)$,  $W=\Big\{\sigma_j\sigma_i(z)\,\Big|\, (i,j)\in A\Big\}$ for some $A\subset I$, and setting
\begin{equation}\label{hard} A^1:=\Big \{i\in[n]\,\Big|\,\exists j\in [n], (i,j)\in A \;\mbox{ or } \; (j,i)\in A\Big\},
\end{equation}
by using Cauchy-Schwarz inequality, the  expression \eqref{defRbis} provides 
\begin{equation}\label{automnebis}
K_0(z,W):=\sup_\alpha \sum_{(i,j)\in A} 2 \alpha_i\alpha_j\leq \sup_\alpha \sum_{(i,j)\in A^1\times A^1,i\neq j}  \alpha_i\alpha_j\leq 1-\frac{1}{|A^1|}=1-\frac1{|]z,W[|},
\end{equation}
where the supremum runs over all vectors $\alpha=(\alpha_1,\ldots,\alpha_n)$ with positive coordinates satisfying $\alpha_1+\cdots+\alpha_n=1$.
Let  $V_+,V_-\subset S_1(z)$ and  $W_+,W_-\subset S_2(z)$   with $V_+\supset]z,W_+[$ and $V_-\supset]z,W_-[$ and satisfying condition \eqref{condVW}.  Applying \eqref{automnebis} for $W=W_-$ and 
$W=W_+$ gives 
\[\frac{\1_{V_+\neq \emptyset}}{1-K_0(z,W_+)}+\frac{\1_{V_-\neq \emptyset}}{1-K_0(z,W_-)} \\
\leq |]z,W_+[|+|]z,W_-[|\leq  |V_+|+|V_-|\leq n,
\]
since $V_-$ and $V_+$ are disjoint.
If  $W_+=\emptyset$ and $W_-\neq \emptyset$ then since $V_+\cap V_-=\emptyset$ one has 
\[ \frac{\1_{V_+\neq \emptyset}}{1-K_0(z,W_+)}+\frac{\1_{V_-\neq \emptyset}}{1-K_0(z,W_-)}
=\1_{V_+\neq \emptyset} + \frac{1}{1-K_0(z,W_-)}\leq \1_{V_+\neq \emptyset}+|V_-|\leq n,\]
and if $(W_+,W_-)=(\emptyset,\emptyset)$ then 
\[\frac{\1_{V_+\neq \emptyset}}{1-K_0(z,W_+)}+\frac{\1_{V_-\neq \emptyset}}{1-K_0(z,W_-)} = \1_{V_+\neq \emptyset}+\1_{V_-\neq \emptyset}\leq n.\]
Thus one gets  $r_1(z)\geq 4/n$ and according to Theorem \ref{thmprincbis} the $W_1$-entropic curvature is bounded from below by $\min_z r_1\geq 4/n$. Asymptotically as $n$ goes to $+\infty$, this lower bound  is the best  one may expect (see \cite[Corollary 4.5]{GRST14}).

The estimate of the lower bound $\overline{r}=\min_{z\in \X} \overline{r}(z)$ on the $T_{\overline{c}}$-entropic curvature of the space is very similar ($\overline{c}=(\overline{c}_t)_{t\in (0,1)}$ with $\overline{c}_t$ defined by \eqref{defcout2}).  Let $W_+,W_-\subset S_2(z)$ with  $d(w_-,w_+)=4$ for all $(w_-,w_+)\in W_-\times W_+$. If $W_+\neq \emptyset$ and $W_-\neq \emptyset$ then inequality \eqref{automnebis} gives \begin{align*}\frac{\1_{W_+\neq \emptyset}}{-\log K_0(z,W_+)}+\frac{\1_{W_-\neq \emptyset}}{-\log K_0(z,W_-)}
&\leq \frac{\1_{W_+\neq \emptyset}}{1-K_0(z,W_+)}+\frac{\1_{W_-\neq \emptyset}}{1-K_0(z,W_-)} \\
&\leq  |]z,W_+[|+|]z,W_-[|\leq n,
\end{align*}
since $]z,W_+[\cap ]z,W_-[=\emptyset$.
If $W_+\neq \emptyset$ and $W_-=\emptyset$,  then 
\[\frac{\1_{W_+\neq \emptyset}}{1-K_0(z,W_+)}+\frac{\1_{W_-\neq \emptyset}}{1-K_0(z,W_-)}=\frac{1}{1-K_0(z,W_+)}\leq \frac{1}{1-K_0(z,S_2(z))}=n.\]
In any cases $\overline{r}(z)\geq 1/n$ and Theorem \ref{thmprincbis} ensures that the $T_{\overline{c}}$-entropic curvature of the hypercube is bounded from below by $\overline{r}\geq 4/n$. We know from \cite{Sam21} that this lower-bound is asymptotically optimal in $n$. Indeed, one may recover the optimal $T_2$-transport entropy inequality for the standard Gaussian measure on $\R$ from the transport entropy inequality with cost $T_{\overline{c}}$ derived from this entropic lower bound (see  \cite[Lemma 4.1]{Sam21}). 

Let us now compute $\widetilde r_2(z)=1-\sup_{W\subset S_2(z)} \widetilde{K}_0(z,W)$.
Using the above notations, for $z\in \X$, one has 
\[\widetilde{K}_0(z,W)=\widetilde{K}(z,W):=\sup_\beta  \Biggl\{  \sum_{(i,j)\in A} 2 \sqrt{\beta_i}\sqrt{\beta_j}
-\sum_{(i,j)\in A^1\times A^1, i\neq j} \sqrt{\beta_i}\sqrt{\beta_j}\Biggr\},\]
where the supremum runs over all $\beta=(\beta_i)_{i\in A^1}$ such that $\beta_i\geq 0$ and $\sum_{i\in A^1} \beta_i=1$. Obviously $\widetilde{K}_0(z,W)=0$ and therefore $\widetilde r_2(z)=1$. 
Theorem \ref{Thmstructure} indicates that the $\widetilde T_2$-entropic curvature $\widetilde \kappa_2$ of the discrete hypercube is bounded from below by $\widetilde r_2\geq 1$. We know that this lower bound is asymptotically optimal as $n$ goes to $\infty$. Indeed according to \cite[Corollary 5.5]{GRST14}, the modified logarithmic inequality \eqref{logsobT3} given by Theorem \ref{Logsobbis} with $\widetilde \kappa_2=1$ implies the well-known Gross logarithmic Sobolev inequality for the standard Gaussian measure with optimal constant.

Since  $\sigma_i\sigma_i(z)=z$ for any $i\in[n]$ and any $z\in \X$, each move $\sigma_i$ is used at most one time along any discrete geodesic. It follows that the $C$-displacement convexity property \eqref{deplacebis} holds with the cost $C_t(\widehat \pi)=\widetilde C_t^1(\widehat \pi)$, $t\in (0,1)$. 
As previously mentioned, the quantities $\Pi^{\sigma_i}_\rightarrow(x)$ and 
$\Pi^{\sigma_i}_\rightarrow(y)$ involved in the definition of the costs $\widetilde T_2$ and $\widetilde C_t^1(\widehat \pi)$ are given by \eqref{Pisigma}.
Thus, one exactly recovers the results of \cite[Theorem 2.5]{Sam21} for the uniform probability measure $\mu_0$ on the discrete hypercube. 
Applying Theorem \ref{Logsobbis}, gives the modified logarithmic Sobolev inequalities given in \cite[Comments (d) of Theorem 2.5]{Sam21} for $\mu_0$. 

We want now to go a step further by considering perturbation measures $m_v$ of $m_0$. In order to apply Theorem \ref{entropicperturb}, one needs to estimate the quantity $D_sv(x,y)$ given  by \eqref{mulh} for any $x,y\in \{0,1\}^n$. For that purpose let us introduce some kind of discrete Hessian matrix for the potential $v$.
For any $z\in\{0,1\}^n$ ($n\geq 2$), and $\{i,j\}\subset [n]$, 
one denotes
\[z_{\overline{ij}}:=(z_1,\ldots,z_{i-1},z_{i+1},\ldots, z_{j-1},z_{j+1},\ldots z_n)\in \{0,1\}^{n-2},\]
and one uses the notation $z_{\overline{ij}}z_iz_j=z$.
 Let $Hv(z)$ denote the symmetric matrix with off-diagonal entries and for $i\neq j$, \[(Hv(z))_{ij}:=\partial_{ij}^2v(z_{\overline{ij}}),\]
where 
\[\partial_{ij}^2v(z_{\overline{ij}}):=v(z_{\overline{ij}}11)+v(z_{\overline{ij}}00)-v(z_{\overline{ij}}01)-v(z_{\overline{ij}}10).\]
The minimum and maximum eigenvalues of the symmetric matrix $Hv(z)$ are denoted, respectively, by $\lambda_{\rm min}(Hv(z))$ and $\lambda_{\rm max}(Hv(z))$. 
We know that $\lambda_{\rm max}(Hv(z))\geq 0$ and $\lambda_{\rm min}(Hv(z))\leq 0$ since the matrix $Hv(z)$ has off diagonal. Let  also
\[\lambda_{\rm max}^\infty(Hv):=\max_{z\in\X}  \lambda_{\rm max}(Hv(z)) \quad \mbox{and}\quad \lambda_{\rm min}^\infty(Hv):=\min_{z\in\X}  \lambda_{\rm min}(Hv(z)).\] 
\begin{lemma}\label{Dvcube}
Let $v:\{0,1\}^n\to \R$ be a potential. If for any $z\in\{0,1\}^n$, the matrix $Hv(z)=V$ does not depend on $z$, then  $V_{ij}=\partial_{ij}^2v(z_{\overline{ij}})$,  and one has for any $x,y\in \X$ with $d(x,y)\geq 2$,
\[\int_0^1 D_sv(x,y)\,q_t(s)\, ds= \frac{2\sum_{\{i,j\}\subset [n] }(x_i-y_i)(x_j-y_j)V_{ij}} {d(x,y)(d(x,y)-1)}  \geq \frac{\lambda_{\min}(V)}{ d(x,y)-1}.\]
In any other cases we also have
\[\int_0^1 D_sv(x,y)\,q_t(s)\, ds\geq \frac{\lambda_{\rm min}^\infty(Hv)}{d(x,y)-1} \sum_{k=1}^{d(x,y)-1}\frac1k\geq \lambda_{\rm min}^\infty(Hv).\]
\end{lemma}
As an example, let $v$ be the potential  defined by 
\begin{equation}\label{vacances}
v(z)=\sum_{i\in [n]} T_i z_i +\frac{1}{2} \sum_{i,j\in [n], i\neq j} V_{ij}\, z_iz_j,
\end{equation}
where $T=(T_1, \ldots, T_n)\in \R^n$ and  $V=(V_{ij})_{i,j\in [n]}$ is a symmetric matrix of real coefficients  with off diagonal.  In that case $Hv(z)= V$ for any $z\in \X$, and  
Theorem \ref{entropicperturb} and Lemma \ref{Dvcube} imply that  the relative entropy with respect to $m_v$ satisfies the $C^v$-displacement convexity property \eqref{deplacebis}
 with for any $t\in (0,1)$,
 \begin{align}\label{mauvais}
 C^v_t(\widehat \pi)&\geq \iint \Big(\frac4n c_t(d(x,y))+{\lambda_{ \min}(V)}\,d(x,y)\1_{d(x,y)\geq 2} \Big)d\widehat\pi(x,y)\nonumber \\
 &\geq \iint d(x,y)\Big(\frac2n (d(x,y)-1) +\lambda_{\min}(V)\Big)\1_{d(x,y)\geq 2}\,d\widehat\pi(x,y) 
 \end{align} 
Observe that one needs $|\lambda_{\min}(V)|=-\lambda_{\min}(V)$ smaller than constant over $n$ to get positive curvature for any integer $n$ from the last estimates. This condition is very strong in high dimensions as regard to the condition we will present now  by applying directly Theorem \ref{thmprinc} on the space $(\X, d,L_v,m_v)$ as explained in Section \ref{sectionpertu}.

For that purpose, let us first observe  that for any  $z\in\{0,1\}^n$ and any  $i\neq j $, 
\begin{align}\label{Dvijcube}
Dv(z,\sigma_i\sigma_j(z))&= v(\sigma_i\sigma_j(z))+v(z)-v(\sigma_i(z))-v(\sigma_j(z))\nonumber\\
&=(2z_i-1)(2z_j-1)\,\partial_{ij}^2v(z_{\overline{ij}}).
\end{align}
Therefore, applying Theorem \ref{thmprinc}, the entropic curvature $\kappa^v$ of the space $(\X, d,L_v,m_v)$ is bounded from below by $r^v=-2\log\big(\max_{z\in\{0,1\}^n} K^v(z,S_2(z))\big)$ with, according to \eqref{defKLv}, 
\[K^v(z, S_2(z)):=\sup_\alpha\Big\{2\sum_{\{i,j\}\subset [n]} e^{-(2z_i-1)(2z_j-1)\,\partial_{ij}^2v(z_{\overline{ij}})/2} \alpha_i\alpha_j\Big\}.\]
In order to  estimate  this key quantity, let us introduce some notations.
For any $z\in\{0,1\}^n$, let $|Hv(z)|$ denotes the symmetric matrix with off diagonal and with coefficients $(|Hv(z)|)_{ij}:=|(Hv(z))_{ij}|$, $\{i,j\}\subset[n]$.  Setting 
\[|Hv(z)|_{\max}:=\max_{\{i,j\}\subset [n]} |(Hv(z))_{ij}|\quad \mbox{and}\quad |Hv|_{\max,\infty} :=\sup_{z\in \{0,1\}^n} |Hv(z)|_{\max},\]
since $Hv(z)$ is a symmetric matrix, one easily checks that 
\[\big|\lambda_{\min}(Hv(z))\big|\leq \lambda_{\max}|Hv(z)|,\qquad  \big|\lambda_{\min}^\infty(Hv)\big|\leq \lambda_{\max}^\infty|Hv|,\]
and since it has an off diagonal
\[|Hv(z)|_{\max}\leq \min\big\{\big|\lambda_{\min}|Hv(z)|\big|,\lambda_{\max}|Hv(z)|\big\},\]
and
\[ |Hv|_{\max,\infty}\leq\min\big\{\big|\lambda_{\min}^\infty|Hv|\big|,\lambda_{\max}^\infty(Hv)\big\}.\]

\begin{lemma}\label{lemcubeK_v} With the above notations, let 
\[\rho(v):= 1+ \frac{\lambda_{\min}^\infty(Hv)}2 - \frac{\lambda_{\max}^\infty|Hv|}2 \,k\Big(\frac{|Hv|_{\max,\infty}}{2}\Big) ,\]
with $k(s):=\frac1s(e^s-s-1)$, $s>0$.
For any $z\in \{0,1\}^n$ and for any $W\subset S_2(z)$, one has
\[ K^v(z, W)\leq 1-\inf_{\alpha}\Big\{\rho(v)\sum_{i\in A^1}\alpha_i^2\Big\}, \]
where the subset of indices $A^1\subset[n]$ is given by \eqref{hard}, and the infimum runs over all vectors $\alpha$ with positive coordinates $\alpha_i$ satisfying $\sum_{i\in A^1}\alpha_i=1$. Moreover if $\rho(v)>0$, then it holds for any $z\in\{0,1\}^n$ and for any $W\subset S_2(z)$,
\[-\frac{\lambda_{\min}^\infty(Hv)}2\leq \widetilde K^v(z,W)\leq -\frac{\lambda_{\min}^\infty(Hv)}2 + k\Big(\frac{|Hv|_{\max,\infty}}{2}\Big)= 1- \rho(v).\]
\end{lemma}
The proof of this lemma is postponed in Appendix B. Since $A^1=[n]$ for $W=S_2(z)$, the upper estimate of $K^v(z, S_2(z))$  of this lemma and Theorem \ref{thmprinc} give :  for $\rho(v)\leq 0$, $K(z,S_2(z))\leq 1-\rho(v)$ and therefore  the  entropic curvature $\kappa^v$ of the space $(\X, d,L_v,m_v)$ is bounded from below by
\[ \kappa^v\geq \min_z r^v(z)\geq -2\log(1-\rho(v))\geq 2\rho(v),\]
and for $\rho(v)> 0$, $K(z,S_2(z))\leq 1-\frac{\rho(v)}n$ and therefore
\[\kappa^v\geq -2\log\Big(1-\frac{\rho(v)}n\Big)\geq \frac{2\rho(v)}n.\]

Applying also Theorem \ref{thmprincbis}  and Theorem \ref{Thmstructure} together with Lemma \ref{lemcubeK_v} easily provides the next result. 

\begin{proposition}\label{prophypercube}
On the discrete hypercube  $\X=\{0,1\}^n$, let $m_v$ denotes the measure with density $e^{-v}$ with respect to the counting measure $m_0$, with $v:\{0,1\}^n\to \R$. 
If $\rho(v)>0$ then,  denoting $\kappa_1^v$ (respectively $\tilde\kappa^v$, $\overline{\kappa}^v$,$\tilde\kappa_2^v$) the $W_1$-entropic curvature of the space $(\X, d,L_1,m_v)$  (respectively the $\widetilde{T}$, $T_{\overline c}$ and $\widetilde T_2$-entropic curvature of the space), one has 
\[\kappa_1^v\geq  \min_z r_1^v(z)\geq \frac{4\rho(v)}n,\qquad \tilde\kappa^v \geq \frac{2\rho(v)}{n}\Big(1-\frac{\rho(v)}{2n}\Big),\]
\[\overline{\kappa}^v \geq \min_z \overline{r}^v(z)\geq \frac{4\rho(v)}n,\qquad \widetilde\kappa_2^v \geq\min_z \widetilde r_2^v(z) \geq\rho(v).\]
\end{proposition}
The lower bounds on $r_1^v(z)$ and $\overline{r}^v(z)$ follow from the first part of Lemma \ref{lemcubeK_v}, by adapting the  arguments we have used at the beginning of this section in order to estimate $r_1(z)$ and $\overline{r}(z)$. It suffices to observe that according to Lemma \ref{lemcubeK_v}, if $\rho(v)>0$ then for any $z\in\{0,1\}^n$, $W\subset S_2(z)$,
\[  K^v(z, W)\leq 1-\rho(v)\inf_{\alpha}\Big\{\sum_{i\in A^1}\alpha_i^2\Big\}=1-\frac{\rho(v)}{|A^1|}=1-\frac{\rho(v)}{|]z,W[|}.\]
The other details of the proofs are left to the reader.

{\bf Comments:}
\begin{enumerate}[label=(\roman*)]
\item As an example, if  the potential $v$ is given by \eqref{vacances} then 
\[\rho(v)=1+ \frac{\lambda_{\min}(V)}2 - \frac{\lambda_{\max}|V|}2 \,k\Big(\frac{|V|_{\max}}{2}\Big) .\]

If $V=0$, then $\rho(v)=1$ and $\mu_v=m_v/m_v(\{0,1\}^n)$ is the product of Bernoulli measures with parameter $p_i=\frac{e^{u_i}}{1+e^{u_i}}$. Therefore, all the entropic curvature lower bounds we get are the same as for the uniform probability measure $\mu_0$ on $\{0,1\}^n$. 

Observe that contrarily to the result \eqref{mauvais} following from Theorem \ref{entropicperturb}, we don't need $|\lambda_{\min}(V)|$ to be of order constant over $n$ for positive entropic curvature, we just need $|\lambda_{\min}(V)|$ and $\lambda_{\max}|V|k\Big(\frac{|V|_{\max}}{2}\Big)$ to be  bounded.

\item Corollary \ref{Transport} and \eqref{transportentC_0^D} provide transport-entropy inequalities associated to different types of concentration inequalities for the measure $\mu_v$. As example, since $\widetilde\kappa^v_2\geq \rho(v)$, according to \eqref{transportentC_0^D}, for any $\nu\in\Pc(\{0,1\}^n)$
\[\frac{\rho(v)}2 \inf_{\pi\in\Pi(\mu_v,\nu)}\Big\{\int \sum_{i\in [n]}    h_0\left(\Pi^{\sigma_i}_\rightarrow(x)\right) d\mu_v(x)+\int \sum_{i\in [n]}  h_1\left(\Pi^{\sigma_i}_\leftarrow(y)\right) d\nu(y)\Big\}\leq \HR(\nu|\mu_v).\]
If $\rho(v)>0$ then by using usual duality arguments as in \cite{GRST14}, one  gets exponential inequalities  for the class of real function   $g$  satisfying : for all $x,y\in \{0,1\}^n$,
\begin{equation}\label{hypof}
g(y)-g(x) \leq \sum_{i\in[n]} a_i(y)\1_{x_i\neq y_i}+\sum_{i\in[n]} b_i(x)\1_{x_i\neq y_i},
\end{equation}
where the $a_i$'s and $b_i$'s are non-negative functions. Integrating this inequality with respect to $\pi\in \Pi(\mu_v,\nu)$ provides,
\begin{align*}
    \int g\,d\nu-\int g\, d\mu_v &\leq \int \sum_{i\in[n]} a_i(y) \Pi^{\sigma_i}_\leftarrow(y)  d\nu(y)+\int \sum_{i\in[n]} b_i(x) \Pi^{\sigma_i}_\rightarrow(x)  d\mu_v(y)\\
    &\leq \frac{\rho(v)}2 \int \sum_{i\in [n]}  h_1\left(\Pi^{\sigma_i}_\leftarrow\right) d\nu + \rho(v)\int \sum_{i\in [n]} \Big(\frac{a_i}{\rho(v)}-\log\big(1+\frac{a_i}{\rho(v)}\big)\Big) d\nu\\
    &\quad + \frac{\rho(v)}2 \int \sum_{i\in [n]}  h_0\left(\Pi^{\sigma_i}_\rightarrow\right) d\mu_v + \rho(v)\int \sum_{i\in [n]} \Big(e^{-\frac{b_i}{\rho(v)}}+\frac{b_i}{\rho(v)}-1\Big) d\mu_v
\end{align*}
since  $\sup_{u\geq 0} \{a u- h_1(u)/2\}=a-\log(1+a):=\ell_1(a)\leq a^2/2$ and $\sup_{u\geq 0} \{b u- h_0(u)/2\}=e^{-b}+b-1:=\ell_0(b)\leq b^2/2$ for any $a,b\geq 0$.
Optimizing the last inequality over all $\pi\in \Pi(\mu_v,\nu)$, using the above transport entropy inequality, and then optimizing over all probability measures $\nu$, the well known duality formulae, $\log\int e^h d\mu=\sup_{\nu\in{\mathcal P}(\X)} \big\{\int h\,d\nu-\HR(\nu|\mu_v)\big\}$, implies 
\begin{align*} 
\log\int \exp\Big[g-\rho(v) \sum_{i\in [n]} \ell_1\Big(\frac{a_i}{\rho(v)}\Big)\Big] d\mu_v&=\sup_{\nu\in \Pc(\{0,1\}^n)} \Big\{\int \Big[g-\rho(v) \sum_{i\in [n]} \ell_1\Big(\frac{a_i}{\rho(v)}\Big)\Big] d\nu -\HR(\nu|\mu_v)\Big\}\\
&\leq \int g\, d\mu_v + \int \rho(v) \sum_{i\in [n]} \ell_0\Big(\frac{b_i}{\rho(v)}\Big) \,d\mu_v\\
\end{align*}
Then Chebychev inequality implies the following new general deviation inequality,
\[\mu_v\Big(g\geq \mu_v[g]+u +\rho(v) \sum_{i\in[n]}\ell_1\Big(\frac{a_i}{\rho(v)}\Big)+\rho(v) \mu_v\Big[\sum_{i\in[n]}\ell_0\Big(\frac{b_i}{\rho(v)}\Big)\Big]\Big)\leq e^{-u},\qquad u\geq 0.\]  
As a byproduct, this exponential inequality gives convex concentration properties for the measure $\mu_v$. Namely, if $f$ is a smooth 1-Lipschitz convex function on $[0,1]^n\supset\{0,1\}^n$, then  the hypothesis \eqref{hypof} holds for $g=\lambda f$, $\lambda> 0$, with $a_i=\lambda |\partial_i f|$ and $b_i=0$, and it also  holds  for $g=-\lambda f$ with  $b_i=\lambda |\partial_i f|$ and $a_i=0$. Since  either  $\rho(v) \sum_{i\in[n]}\ell_1\big(\frac{a_i}{\rho(v)}\big)\leq \frac{\lambda^2|\nabla f|^2}2\leq \frac{ \lambda^2}{2}$ or $\rho(v)\sum_{i\in[n]}\ell_0\big(\frac{b_i}{\rho(v)}\big)\leq \frac{\lambda^2|\nabla f|^2}2$, setting $s=\frac{u}\lambda+\frac{\lambda}{2\rho(v)}$ or $s=\frac{u}\lambda+\frac{\lambda\mu_v[|\nabla f|^2]}{2\rho(v)}$ and  optimizing  over $\lambda$, provides the following exponential concentration bounds for the deviations of $f$ above or under its mean, for any $s\geq 0$. 
\[\mu_v(f\geq \mu_v[f]+s)\leq e^{-\frac{\rho(v)s^2}2} \quad\mbox{and}\quad \mu_v(f\leq \mu_v[f]-s)\leq e^{-\frac{\rho(v)s^2}{2\mu_v[|\nabla f|^2]}}\leq e^{-\frac{\rho(v)s^2}{2}}.\]

Note that another application of the above exponential inequality is concentration inequalities for suprema of empirical processes under $\mu_v$, by following the lines of proof \cite[Corollary 3.3]{Sam07} reached in the independent case (for  product measures $\mu$ on $\{0,1\}^n$).   

\item 
 Modified logarithmic Sobolev  and Poincaré type of inequalities for the measure $\mu_v$ are consequences of Theorem \ref{Logsob} and Theorem \ref{Logsobbis}. As example, since   $\widetilde\kappa^v_2\geq \rho(v)$,  Theorem \ref{Logsobbis} ensures that if $\rho(v)>0$ then for any positive function $f$ on $\{0,1\}^n$,
\begin{equation*}
{\rm Ent}_{\mu_v}(f)\leq  \int \sum_{i\in[n]} \frac{ \rho(v) }2 \, h^*\left(\frac{2}{\rho(v)} [\partial_{\sigma_i} \log f]_-\right)  f\,d\mu_v\leq \frac 1{\rho(v)} \int \sum_{i\in[n]} {[\partial_{\sigma_i} \log f]_-[\partial_{\sigma_i} f]_-} \,d\mu_v,
\end{equation*}
and also   for any  $g:\{0,1\}^n\to \R$, 
\[{\rm Var}_{\mu_v}(g)\leq \frac{1}{2 \rho(v)} 
\int \sum_{i\in[n]}(\partial_{\sigma_i} g)^2 d\mu_v
.\]
\end{enumerate}
Let us adapt these results for Ising models. Let $\Lambda$ be a finite set of vertices, $|\Lambda|=n$, of a finite graph {$G_{\Lambda}=(\Lambda,E_\Lambda)$} without  multiple edges and without loops. If two vertices $i$ and $j$ are neighbours, one denotes $i\sim_{\!\Lambda} j$. Let $m_0$ denote the counting measure on $\{-1,1\}^\Lambda$ and  $m_w=e^{-w}m_0$ be the measure on $\{-1,1\}^\Lambda$ with potential of interaction $w$ defined by \eqref{defwising}.
 
Since $z_i\in \{-1,1\} $ if and only if $\frac{z_i+1}2\in\{0,1\}$, by a simple change of variable, the results  of Proposition \ref{prophypercube} for the measure $m_v$ transpose to the measure $m_w$ by replacing the quantity $\rho(v)$ with the quantity $\rho_\beta(W)$, defined  by \eqref{defrhoising},
where $W$ is the $n$ by $n$ symmetric matrix of interaction coefficients with off diagonal.

Assume first that all coefficients $W_{ij}$ are non-negative, $W=|W|$, and 
\[\rho_\beta(W)=1-\lambda_{\max}(W) \big(e^{2\beta|W|_{\max}}-1\big) .\] 
Due to the Perron-Frobenius theorem,  one has 
\begin{align*}\lambda_{\max}(W)&=\max\big\{|\lambda|\,\big|\,\lambda\mbox{ is an eigenvalue of }W\big\}=\sup_{x\in \R^\Lambda_+, |x|_2=1} |Wx|_2 \\
&\leq S_\infty(W):=\max_{i\in \Lambda}\sum_{j, j\sim_{\! \Lambda}i} W_{ij},
\end{align*}
where $|y|_2^2=\sum_{i\in\Lambda} y_i^2$ for $y\in \R^\Lambda$.
It follows that 
$\widetilde\rho(W)>0$ as soon as 
\[\beta\leq \frac{1}{2|W|_{\max}(1+S_\infty(W))}<\frac{1}{2|W|_{\max}}\log\Big(1+\frac 1{S_\infty(W)}\Big).\]

As example, in the simplest ferromagnetic model, $W$ is the adjacency matrix $A$ of the graph $G_\Lambda$, 
\[A_{ij} := \hspace{0.2cm }\begin{cases}
1\hspace{0.2cm}  \text{if} \hspace{0.2cm} i{\sim}_{\!\Lambda} j,\\
 \text{0} \hspace{0.2cm}   \text{otherwise} \hspace{0.06cm}.
\end{cases}\]
In that case, $|W|_{\max}=1$ and  $S_\infty(W)=\Delta(G_\Lambda)$ is the maximal degree of $G_\Lambda$. Therefore,  $\widetilde\rho(W)>0$ as soon as 
\[\beta< \frac1{2(1 + \Delta(G_\Lambda))}.\]
Observe that for  $G_\Lambda$ being  the complete graph,  the so-called Curie-Weiss model, $W_{ij}=A_{ij}=1 $ for $i\neq j$, $\lambda_{\max}(W)=S_\infty(W)=n-1$ and $\beta< \frac1{2n}$ is a sufficient condition for $\rho_\beta(W)>0$. 
The critical $\overline\beta_n$ for the Curie-Weiss model for which some Poincaré 
inequality  is known to fail beyond this point is $\overline\beta_n=\frac1{n}$ (see \cite{DLP09}). It remains a challenging problem to get positive entropic curvature up to this  critical value. For $G_\Lambda$ being  a subgraph of a  graph $G_\infty$, by selecting a set of vertices $\Lambda$, if $G_\infty$ has uniform bounded degree,  then $\Delta(G_\Lambda)\leq \Delta(G_\infty)<\infty$ provides a the uniform condition over all boxes $\Lambda$ for positive entropic curvature, namely $\beta< \frac1{2(1 + \Delta(G_\infty))}$. As a classical example, one may choose  $G_\infty=\Z^d$  with $\Delta(G_\infty)=d$ and $\Lambda=[-N,N]^d$, $N\in \N$. 
 
Without any assumption on the sign of the $W_{ij}$'s, the Perron-Frobenius theorem implies $\lambda_{\max}(M)\leq \lambda_{\max}|M|$ and therefore 
\[
\rho_\beta(W)\geq 1-\lambda_{\max}|W| \big(e^{2\beta|W|_{\max}}-1\big) >0,
\]
as soon as 
\begin{eqnarray}\label{condtilderho}
\beta\leq\frac{1}{2|W|_{\max}(1+\lambda_{\max}|W|)},
\end{eqnarray} and therefore under the following stronger Dobrushin-type  condition
\begin{eqnarray}\label{condtilderhobis}
\beta\leq \frac{1}{2|W|_{\max}(1+S_\infty|W|)}, \quad S_\infty|W|:=\max_{i\in \Lambda}\sum_{j, j\sim_{\! \Lambda}i} |W_{ij}|.
\end{eqnarray}
Note that  Erbar-Maas entropic curvature for general Ising models have been also studied in \cite[Theorem 4.1]{EHMT17}. Their condition for positive Erbar-Maas entropic curvature is rather comparable to a Dobrushin-type  conditions such  as in \eqref{condtilderhobis}. In particular for the ferromagnetic Cury-Weiss model, Erbar-Maas entropic positive curvature is reached for $\beta\leq 0.284/n$ (see \cite[Corollary 4.5]{EHMT17}). 

Actually, if $\lambda_{\max}(W)$ is very small as regard to $\lambda_{\max}|W|$, then  \eqref{condtilderho} is strongly improved by keeping the expression of $\rho_\beta(W)$. As  example, let us consider the case where $W=-A$ where $A$ is the adjacency matrix of a complete graph $G_{\Lambda}$, that we call  antiferromagnetic Curie-Weiss model.  Since $\lambda_{\max}(W)=-\lambda_{\min}(A)=1$ and  $\lambda_{\max}|W|=\lambda_{\max}(A)=n-1$, one has 
\[\rho_\beta(W)=1-2\beta-2(n-1)\beta\,k(2\beta)\geq 1-2n\beta\,k(2\beta)\geq 1-2n\big(e^\beta-1)^2.\]
One  checks that $\rho_\beta(W)>0$ as soon as $\beta\leq \frac1{1+\sqrt{2n}}$, which is a weaker condition than \eqref{condtilderho}, $\beta\leq \frac{1}{2n}$.

Finally, let us consider the case where the  matrix of interaction coefficient  $W$ is random, with ``good'' concentration's properties as its size $n$  increases. The parameter $\beta=\beta_n$ will therefore depend on $n$. Our result reads as follows, if with high probability as $n$ goes to $\infty$, $\lambda_{\max}(W)$ is of order $\gamma_n$,   $\frac{|W|_{\max}}{\gamma_n}$ is bounded and $\lambda_{\max}|W|$ is of order $o(\gamma_n^2)$, then   with high probability
$\rho_\beta(W)>0$  as soon as $\beta_n\leq \frac{1-\varepsilon}{2\gamma_n}$, $\varepsilon >0$.
More precisely, assume that with high probability as $n$ goes to $\infty$,
\begin{equation}\label{condSP}
\lambda_{\max}(W)\sim \gamma_n  \quad\mbox{and} \quad \lambda_{\max}|W|\,{|W|_{\max}}e^{{|W|_{\max}}/{\gamma_n}}=o(\gamma_n^2),
\end{equation}
then  with high probability for $n$ sufficiently large, for any $\beta_n< \frac{1-\varepsilon}{2\gamma_n}$,
\begin{equation*}
\rho_\beta(W)\geq \varepsilon-o(1)-\frac{\lambda_{\max}|W|\,{|W|_{\max}}}{\gamma_n^2}\,e^{{|W|_{\max}}/{\gamma_n}}\geq \varepsilon-o(1)>0,
\end{equation*}
 since $k(\alpha s)\leq \alpha k(s)\leq \alpha s e^s/2$, for  $s>0$ and $0\leq \alpha\leq 1$.
  As example, assume that all the $W_{ij}$'s are independent standard gaussian random variables as in the celebrated Sherrington-Kirkpatrick (SK) model from spin glass theory \cite{Tal10}. Setting $\widetilde W_{ij}=\frac{|W_{i,j}|-\mathbb{E}(|W_{i,j}|)}{\sqrt{\mathbb{V}{\rm ar}(|W_{i,j}|)}}$, $i,j\in[n], $ one easily checks that  
  \[\lambda_{\max}|W|\leq \sqrt{\mathbb{V}{\rm ar}(|G|)} \lambda_{\max}(\widetilde W)+\mathbb{E}(|G|),\]
  where $G$ is a standard gaussian random variable.
  According to concentration results of the largest  eigenvalue for symmetric random matrices with independent entries (see \cite{AGZ10}), it holds  
  \[\lim_{n\to\infty}\frac1{\sqrt{n}}\lambda_{\max}(\widetilde W)=2=\lim_{n\to\infty}\frac1{\sqrt{n}}\lambda_{\max}( W),\]
  almost surely, and  therefore, setting $\gamma_n=\sqrt{n}$,
  \[\lambda_{\max}(W)\sim 2 \gamma_n  \quad\mbox{and} \quad \limsup_{n\to \infty} \frac{\lambda_{\max}|W|}{\gamma_n}=O(1),\]
  almost surely.
  Moreover, applying the union bound directly, the subgaussianity of the random variables $W_{ij}$ implies 
  \[\mathbb P\big(|W|_\infty \geq \sqrt{2\log(n(n-1))}+u\big)
  \leq n(n-1)\exp\Big\{-\frac12 \big(\sqrt{2\log(n(n-1))}+u\big)^2\Big\}\leq e^{-u^2/2}.\]
  From this concentration bound, applied with $u=u_n=o(\gamma_n)=o(\sqrt n)$ with $\sum e^{-u_n^2/2}<\infty$, Borell-Cantelli's Lemma ensures that
  $\lim_{n\to \infty }\frac{|W|_\infty}{\gamma_n}=0$ almost surely.
  As a consequence, putting all together the last concentration results, the conditions \eqref{condSP} holds almost surely for $\beta_n\leq\frac{1-\varepsilon}{4\sqrt n}$, which implies   $\rho_\beta(W)>0$ almost surely for $n$ sufficiently large. These comments   extend if the $W_{ij}$'s have subgaussian tails
  from concentration results for the operator norm  $\|W\|_{op}$  since $|\lambda_{\max}(W)|\leq\|W\|_{op}$ (see e.g. \cite{RJSS23}).  
  Recent  results for operator norm in  random matrix concentration theory \cite{LHY18,BBH23} also allow to consider gaussian interaction matrices $W$ (with  non necessary independent identically distributed entries). 
 
As commented above, first   consequences of the condition $\rho_\beta(W)>0$ 
are  concentration properties for the measure $\mu_w$. It is well known  that the transport-entropy inequality of Corollary \ref{Transport} with cost $W_1$ implies concentration inequalities for real Lipschitz functions on $\{-1,1\}^\Lambda$ with respect to the graph distance on $\{-1,1\}^\Lambda$, $d(x,y):=  \sum_{i\in\Lambda} \1_{x_i\neq y_i} $(see e.g. \cite{BG99}). Namely, we get for any  1-Lipschitz function with respect to $d$, for any $s\geq 0$,
\[\mu_w(f\geq \mu_w(f)+s)\leq e^{-2\rho_\beta(W)s^2/|\Lambda|}.\]

Considering now $\{-1,1\}^\Lambda$ as a subset of $[-1,1]^n$, according to the above comment $(ii)$, the transport-entropy inequality \eqref{transportentC_0^D} with cost $\widetilde C_t^1$ (or $\widetilde T_2)$ gives the following stronger convex concentration property, for any 1-Lipschitz smooth convex function $f$  on $[-1,1]^n$ (with respect to the euclidean distance on $[-1,1]^n$), for any $s\geq 0$, 
\[\mu_w(f\geq \mu_w(f)+s)\leq e^{-\rho_\beta(W)s^2/8} \quad\mbox{and}\quad \mu_w(f\leq \mu_w(f)-s)\leq e^{-\rho_\beta(W)s^2/8}.\]
This convex concentration property is identical to the one reached by Adamczak \& al. in \cite[Proposition 5.4]{AKPS18}. Their result is a consequence of approximate tensorization property of entropy under  the Dobrushin condition $|W|_{\max} S_\infty|W|< 1$. As already mentioned,  this condition is  stronger than $\rho_\beta(W)>0$ in some cases. Our  result also improved the one of \cite[Proposition 5.4]{AKPS18} since we give an explicit expression of the constant $\rho_\beta(W)(W)$ in the deviation bound. Moreover this constant   does not depend on the maximal value of the coefficients  $|T_i|$, $i\in \Lambda$, as opposed to their deviation result in \cite{AKPS18}. 

Poincaré and modified logarithmic Sobolev inequalities are other consequences of the condition $\rho_\beta(W)>0$. Our result is comparable to the one of Bauerschmidt-Bodineau \cite{BB19}: namely, setting $\langle W\rangle :=\lambda_{\max}(W)-\lambda_{\min}(W)$, if $\beta_n \langle W\rangle<1$ then for all positive functions $f$ 
\begin{equation*}
{\rm Ent}_{\mu_w}(f)\leq  \frac12\Big(1+\frac{2\beta_n\langle W\rangle}{1-\beta_n \langle W\rangle}\Big)  \int \sum_{i\in[n]} (\partial_{\sigma_i} \sqrt f)^2 \,d\mu_w = \Big(1+\frac{2\beta_n\langle W\rangle}{1-\beta_n \langle W\rangle}\Big) {\mathcal{E}}_{L_2}\big(\sqrt f, \sqrt f\big).
\end{equation*}
according to \eqref{grrr}. Recall that from the easy bound $4(\sqrt{a}-\sqrt b)^2\leq (a-b)(\log(a)-\log(b))$, one has  $4{\mathcal{E}}_{L_2}\big(\sqrt f, \sqrt f\big)\leq {\mathcal{E}}_{L_2}\big( f, \log f\big)$. For small $\beta_n$, our condition $\widetilde \rho(W)$ 
is close to $\beta_n \lambda_{\max}(W)\lesssim 1/2$, which is of same order as $\beta_n \langle W\rangle<1$ if the spectrum of $W$ is symmetric. Our condition is weaker when $|\lambda_{\min}|$ is much bigger than $\lambda_{\max}$, like in the case of the antiferromagnetic Curie-Weiss model as discuss before (in that case  $\rho_\beta(W)>0$ for $\beta_n\leq \frac1{1+\sqrt{2n}}$, whereas  
$\beta_n \langle W\rangle<1$ means $\beta_n <1/n$).   
  Let us mention also that Eldan-Koehler-Zeitouni \cite{EKZ22} proved a Poincaré inequality with improved  Dirichlet form (associated to the Glauber dynamic) under the same condition   $\beta_n \langle W\rangle<1$ 
by using localization techniques (see also \cite{YR22}).

 \subsection{The lattice $\Z^n$}
 
 In this part,  $m_0$ denotes the counting measure on $\X:=\Z^n$. Recall that the lattice $\Z^n$ is a structured graph with set of moves $\Sc:=\{\sigma_{i+},\sigma_{i-}\,|\, i\in[n]\}$ with 
 $\sigma_{i+}(z):=z+e_i$ and $\sigma_{i-}(z):=z-e_i$ for any $z\in\Z^n$. By Theorem \ref{thmstructure}, it has non-zero entropic curvature and since it is not a finite graph by Bonnet-Myers Theorem \ref{BM}, $\kappa=0$. 
 As an illustrative example, it is easy to see that 
 \[K_{0}\big(z,S_2(z)\big)= \sup_\alpha\Bigg[\sum_{i=1}^n \big(\alpha_{i+}^2+\alpha_{i-}^2\big) + \!\!\!\!\!\sum_{1\leq i<j\leq n} \!\!\! \!\!2(\alpha_{i-}\alpha_{j-}+\alpha_{i-}\alpha_{j+}+ \alpha_{i+}\alpha_{j-}+\alpha_{i+}\alpha_{j+}\big) \Bigg]=1,\] 
 where the supremum runs over all $\alpha=(\alpha_{i+},\alpha_{i-})_{i\in [n]}$ with non negatives coordinates satisfying $\sum_{i\in [n]}(\alpha_{i+}+\alpha_{i-})=1$.
 
 For any integers $d, k_1,\ldots ,k_n$ such that $d=k_1+\cdots +k_n$, let $\binom{d }{k_1,\ldots,k_n}=\frac{n!}{k_1!\cdots k_n!}$ denote the multinomial coefficient.  Since 
 \begin{equation}\label{expL0Zn}
L_0^{d(x,y)}(x,y)= \# G(x,y)= \binom{d(x,y)}{|y_1-x_1|,\ldots,|y_n-x_n|},
\end{equation}
the Schr\"odinger bridge on the space $(\Z^n,d,L_0,m_0)$ between Dirac measures at $x$ and $y$ in $\Z^n$ is given by 
 \begin{align}\label{bridgeZn}
  \nu_t^{x,y}(z)&=\frac{\binom{d(x,z)}{|z_1-x_1|,\ldots,|z_n-x_n|}\binom{d(z,y)}{|y_1-z_1|,\ldots,|y_n-z_n|}}{ \binom{d(x,y)}{|y_1-x_1|,\ldots,|y_n-x_n|} } \binom{d(x,y)}{d(x,z)} \;t^{d(x,z)} (1-t)^{d(z,y)} \1_{[x,y]}(z)\nonumber\\
  &=\binom{|y_1-x_1|}{|z_1-x_1|} \cdots \binom{|y_n-x_n|}{|z_n-x_n|}\;t^{d(x,z)} (1-t)^{d(z,y)} \1_{[x,y]}(z),\qquad z\in \Z^n.
  \end{align}

Let  $m_v$ denote a measure  with  potential   $v:\Z^n\to\R$  with respect to $m_0$, $m_v=e^{-v} m_0$. 
As in the case of the discrete hypercube,  one  defines coefficients  that can be interpreted as local second partial derivatives. For any $i\in [n]$, let 
\[\partial^2_{ii}v(z):=v(z+e_i)+v(z-e_i)-2 v(z),\]
and for $\{i,j\}\subset [n]$,
\[\partial^2_{ij}v(z):=v(z+e_i+e_j)+v(z)-v(z+e_i)-v(z+e_j).\]
One may  check that for any $i,j\in[n]$  and any $\varepsilon_i,\varepsilon_j\in\{-1,1\}$, 
\begin{equation}\label{DvZn}
Dv(z,z+\varepsilon_i e_i+\varepsilon_j e_j)=\varepsilon_i\varepsilon_j \partial_{ij}v(z\wedge(z+\varepsilon_i e_i+\varepsilon_j e_j)),
\end{equation}
where for $z,w\in \Z^n$, $z\wedge w$ denotes the vector with coordinates $\min(z_i,w_i),i\in[n]$.

As an example, if  the potential $v$ is given by the sum of a quadratic and a linear form,
\begin{equation}\label{vquad}
    v(z):=\sum_{i\in [n]} T_i z_i +\frac{1}{2} \sum_{(i,j)\in [n]^2} V_{ij}\, z_iz_j,\quad z\in \Z^n,
\end{equation}
with $(T_1, \ldots, T_n)\in \R^n$ and  $V=(V_{ij})_{i,j\in [n]}$  a symmetric matrix of real coefficients, then easy computations give  
\[\partial_{ij}v(z\wedge(z+\varepsilon_i e_i+\varepsilon_j e_j))=V_{ij}\quad \mbox{and}\quad Dv(z,z+\varepsilon_i e_i+\varepsilon_j e_j)=\varepsilon_i\varepsilon_j V_{ij}\]
for  any $\varepsilon_i,\varepsilon_j\in\{-1,1\}$.

Let us first apply Theorem \ref{entropicperturb} in case the potential $v$ is given by \eqref{vquad}. From the expression \eqref{expL0Zn} and after few computations one exactly gets for any $x,y\in \Z^n$ and for any $t\in(0,1)$, 
\[d(x,y)(d(x,y)-1)\,D_t v(x,y):=2 \sum_{\{i,j\}\subset [n]} V_{ij} (y_j-x_j)(y_i-x_i) +\sum_{i\in[n]} V_{ii}|x_i-y_i|\big(|x_i-y_i|-1\big) .\]
and therefore the relative entropy with respect to $m_v$, $\nu\in \Pc(\X)\mapsto H(\nu|m_v)$, satisfies the $C^v$-displacement convexity property \eqref{deplacebis} along the Schr\"odinger bridge at zero temperature of the space $(\Z^n,d,m_0,L_0)$, with for any $t\in (0,1)$,
\begin{align*}
  C^v_t(\widehat \pi)&:=\iint \Big(2\sum_{\{i,j\}\subset [n]} V_{ij} (y_j-x_j)(y_i-x_i) +\sum_{i\in[n]} V_{ii}|x_i-y_i|\big(|x_i-y_i|-1\big)\Big)\, d\widehat \pi(x,y) \\
  &=\iint \langle y-x,V(y-x)\rangle \, d\widehat \pi(x,y) -\sum_{i\in[n]} V_{ii} |y_i-x_i| \, d\widehat \pi(x,y),
\end{align*}
where $\langle\cdot,\cdot\rangle$ denotes the usual scalar product on $\R^n$.
If the smallest eigenvalue of $V$ is positive then cost function is clearly positive for big values of the euclidean norm of $(x,y)$. Therefore the positivity of $\lambda_{\min}(V)$  plays a central role for the positivity of entropic curvature.

In order to give mild conditions for positive entropic curvature,
let us apply Theorem \ref{thmprinc} on the space $(\Z^n, d,L_1,m_v)$ (the generator $L_1$ is given by \eqref{defL1L2}).  The entropic curvature $\kappa^v$ of this space is bounded from below by $r^v=-2\log\big(\max_{z\in\{0,1\}^n} K^v(z,S_2(z))\big)$ with, according to \eqref{defKLv}, 
\begin{multline*}
K^v(z, S_2(z)):=\sup_\alpha\Big\{2\sum_{\{i,j\}\subset[n]}
\sum_{\varepsilon_i,\varepsilon_j\in\{-1,+1\}}  e^{-\varepsilon_i\varepsilon_j\,\partial_{ij}v(z\wedge(z+\varepsilon_i e_i+\varepsilon_j e_j))/2} \alpha_{\varepsilon_i}\alpha_{\varepsilon_j}\\
+ \sum_{i\in [n]} \sum_{\varepsilon_i\in\{-1,+1\}}
e^{-\partial_{ii}v(z\wedge(z+2\varepsilon_i ))/2} \alpha_{\varepsilon_i}^2 \Big\},
\end{multline*}
where the supremum runs over all vectors $\alpha$ with non-negative coordinates $\alpha_{i+},\alpha_{i-},i\in [n]$ such that $\sum_{i\in[n]} (\alpha_{i+}+\alpha_{i-})=1$ (with the notation $\alpha_{\varepsilon_i}=\alpha_{i+}$ for $\varepsilon_i=1$ and $\alpha_{\varepsilon_i}=\alpha_{i-}$ for $\varepsilon_i=-1$).

For a better understanding, the next lemma present one way to upper bound $K^v(z, S_2(z))$.
\begin{lemma}\label{KLvZ} Given a potential $v:\Z^n\to \R$, for any $z\in \Z^n$, let $Av(z)$  denotes the $n$ by $n$ symmetric matrix  defined by 
\[(Av(z))_{ii}:=e^{-a_{ii}(z)/2}-1,\quad i\in[n],\]  
with $a_{ii}(z):=\min\big(\partial_{ii}v(z),\partial_{ii}v(z-2e_i)\big), i\in[n]$ and 
\[(Av(z))_{ij}:=e^{a_{ij}(z)/2}-1,\quad \{i,j\}\subset[n],\]
with $\displaystyle a_{ij}(z):=\max\big\{-\partial_{ij}v(z), -\partial_{ij}v(z-e_i-e_j),\partial_{ij}v(z-e_i),\partial_{ij}v(z-e_j)\big\}.$\\
If $\lambda_{\max}(Av(z))\leq 0$, then one  has 
\[K^v(z, S_2(z))\leq 1+\frac{\lambda_{\max}(Av(z))}n \quad\mbox{and}\quad 
\widetilde K^v(z)=\sup_{W\subset S_2(z)}\widetilde K^v(z,W)\leq 1+  \lambda_{\max}(Av(z)) .\]
\end{lemma}
The hypothesis $\lambda_{\max}^\infty(Av):=\sup_{z\in \Z^n}\lambda_{\max}(Av(z)) \leq 0$ is not empty. For example,  if  $v$ is given by \eqref{vquad} then for any $z\in \Z^n$, one has $(Av(z))_{ii}=e^{-V_{ii}/2}-1$ and $(Av(z))_{ij}=e^{|V_{ij}|/2}-1$.  Clearly if for any $i\in [n]$, $V_{ii}>0$ and for any $\{i,j\}\subset[n]$, $V_{ij}=0$, then 
\[\lambda_{\max}^\infty(Av)\leq -\big(1-e^{-\,\min_{i\in[n]}(V_{ii})/2}\big)<0.\]
By a continuity argument if for any $i\in [n]$, $V_{ii}>0$ and the values of $|V_{ij}|$, $\{i,j\}\subset[n]$, are sufficiently small then  
$\lambda_{\max}^\infty(Av)<0$. 

For instance, assume that for any $z\in \Z^n$, the matrix $B:=-Av(z)$ has  non negative diagonal entries and  is \textit{diagonally dominant}, that is 
\begin{equation}\label{condGershgorin}
    B_{ii}\geq R_i :=\sum_{j, j\neq i} |B_{ij}|,\quad \forall i\in [n],
\end{equation}
then the Gershgorin's circle theorem ensures that for any eigenvalue $\lambda$ of $B$, there  exists $i\in[n]$ such that $|\lambda-B_{ii}|\leq R_i$. 
As a consequence, one has 
\[\lambda^\infty_{\min}(-Av)\geq \min_{i\in[n],z\in\Z^n}\Big\{(-Av(z))_{ii}-\sum_{j, j\neq i} |Av(z)|_{ij}\Big\}\geq 0,\]
and therefore  $\lambda^\infty_{\max}(Av)\leq 0$.
In particular, if the potential $v$ is given by \eqref{vquad} with  $V_{ii}\geq 0$, then \eqref{condGershgorin} reads as
\begin{eqnarray*}
1- e^{{-V_{ii}}/{2}}\geq \sum_{j\in [n],j\neq i}\big|e^{{V_{ij}}/{2}}-1\big|,\quad \forall i\in [n].
\end{eqnarray*}
 Observe that in case all $V_{ij}$ are non-negative, this inequality implies the  diagonal dominance property for the matrix $V$. And conversely, for sufficiently small coefficients $V_{ij}$,  it  is also close to the  diagonal dominance property of the matrix $V$. As a comment,   the diagonal dominance of the matrix $B$ is known to be equivalent to a \textit{discrete midpoint convexity} property of the quadratic form $v(z)=\sum_{i,j\in [n], i\neq j} B_{ij}\, z_iz_j$ on $\Z^n$ (see details in \cite[Theorem 9]{TT21}). 

Note that the bounds on $K^v\big(z, S_2(z)\big)$ and $\widetilde K^v(z)$ given by Lemma \ref{KLvZ} can be  improved if the potential $v$ is simply given by \eqref{vquad} by introducing other matrices, and even more if the matrix $V$ is specified. 

From the estimates given by Lemma \ref{KLvZ},  Theorem \ref{thmprinc}, \ref{thmprincbis} and \ref{Thmstructure} provide the following  lower bounds for entropic curvatures.

\begin{proposition}\label{propZn} With the above notations, the entropic curvature $\kappa^v$ of the graph space $(\Z^n,d,L_1,m_v)$ satisfies
\[\kappa^v\geq 2\log\big(1+{\lambda_{\max}^\infty(Av)}/n\big)\geq -\frac{2\lambda_{\max}^\infty(Av)}n.\]
Moreover, if  $\lambda_{\max}^\infty(Av)\leq 0$ then the  $W_1$-entropic curvature $\kappa^v_1$, the $\widetilde{T}$-entropic curvature $\widetilde\kappa^v$ and the $\widetilde T_2$-entropic curvature $\widetilde\kappa^v_2$ of this space are non-negative, one has 
\[\kappa^v_1\geq -\frac{2\lambda_{\max}^\infty(Av)}n,\qquad 
\widetilde\kappa^v\geq -\frac{\lambda_{\max}^\infty(Av)}n,\qquad \widetilde\kappa_2^v\geq -\lambda_{\max}^\infty(Av).\]
\end{proposition}
 If $\lambda_{\max}^\infty(Av)< 0$ then 
Theorem \ref{PL} applies with cost $c(d)=d(d-1)$ and $\kappa_c=\kappa^v$  and provides a Prékopa-Leindler type of inequality for the measure $m_v$ with the Schr\"odinger bridges between Dirac measures given by \eqref{bridgeZn}.

Assuming moreover that $m_v(\Z^n)<\infty$, Corollary \ref{Transport} and the first part of Theorem \ref{Logsobbis} provide transport-entropy inequalities for the probability measure $\mu_v=m_v/m_v(\Z^n)$ involving the costs $W_1,T_2,\widetilde{T}$ and $\widetilde T_2$ given by 
\begin{multline*}\widetilde T_2(\pi):=\sum_{i\in [n]} \Big(\int [y_i-x_i]_+ \,d\pi_{_\rightarrow}(y|x)\Big)^2 d\nu_0(x)+\sum_{i\in [n]}\Big(\int [y_i-x_i]_- \,d\pi_{_\rightarrow}(y|x)\Big)^2 d\nu_0(x)\\
 + \sum_{i\in [n]} \Big(\int [x_i-y_i]_+ \,d\pi_{_\leftarrow}(x|y)\Big)^2 d\nu_1(y)+\sum_{i\in [n]}\Big(\int [x_i-y_i]_- \,d\pi_{_\leftarrow}(x|y)\Big)^2 d\nu_1(y),
 \end{multline*}
 for any $\pi\in \Pi(\nu_0,\nu_1)$.
 This expression is a consequence of the following identities, for any $x,y\in \Z^n$,
\[\1_{\sigma_{i+}(x)\in]x,y]}r(x,\sigma_{i+}(x),\sigma_{i+}(x), y)=[y_i-x_i]_+ \quad\mbox{and}\quad\1_{\sigma_{i+}(x)\in]x,y]}r(x,\sigma_{i-}(x),\sigma_{i-}(x), y)=[y_i-x_i]_-.\]

Theorem \ref{Logsobbis} also gives a modified logarithmic Sobolev inequality and a Poincaré inequality. Namely if $\lambda_{\max}^\infty(Av)< 0$, then $\widetilde \kappa_2>-\lambda_{\max}^\infty(Av)$ and for any bounded function $f:\Z^n\to [0,+\infty)$, 
 \begin{equation*}
{\rm Ent}_{\mu_v}(f)\leq \frac{1}{2\widetilde \kappa_2} \int \sum_{i\in [n]}   \big([\partial_{\sigma_{i+}} \log f]_-^2 + [\partial_{\sigma_{i-}} \log f]_-^2\big) f\,d\mu_v,
\end{equation*}
and for any real bounded function $g:\Z^n\to \R$, 
\[{\rm Var}_{\mu_v}(g)\leq \frac{1}{2\widetilde{\kappa}_2} 
\int \sum_{i\in[n]}\big((\partial_{\sigma_i+} g)^2+(\partial_{\sigma_i-} g)^2 \big)\,d\mu_v
.\]

Note that according  Section \ref{sectionrestconv}, all these results extend to restrictions of the measure $m_v$  to a convex subset ${\mathcal C}$ of $\Z^n$, as example any product of segment of $\Z$, ${\mathcal C}:=[c_1, d_1]\times \cdots\times  [c_n, d_n]$ with $c_i,d_i\in \Z$ for all $i\in [n]$.

\section{Non positively curved graphs}\label{sectionneg}

This section and the next ones only concern graph spaces equipped with the uniform measure $G=(\X,d,L_0,m_0)$ denoted also by $(\X,E)$ where $E$ is the set of edges of the graph. This part more specifically deals with graphs with negative entropic curvature, in particular the so called \textit{geodetic graphs} introduced by Ore (see \cite{Oys87}).

\begin{definition}
A graph $G=(\mathcal{X},E)$ is called geodetic if for every two vertices $z$ and $w$ in $G$ there exists a unique 
geodesic connecting $z$ and $w$.
\end{definition}
For example, every tree, every complete graph, every odd-length cycle and the Petersen graph are geodetic graphs. 
The following proposition provides a non-positive lower bound on entropic curvature.  \begin{proposition}\label{propgeodetic}
Let $G=(\mathcal{X},E)$ be a geodetic graph with diameter greater or equal to 2, then the entropic curvature is lower bounded by $\displaystyle \big(2-\max_{z^{\prime},\hspace{0.1cm} z^{\prime}\sim z}\textnormal{deg}(z^{\prime})\big)$.
\end{proposition}
The proof of this proposition is an easy consequence of Theorem \ref{thmprinc}. Since for any $z\in \X$ and any $\ttz\in S_2(z)$, the set $]z,\ttz[$ is reduced to a single vertex, one has
\begin{equation*}
K_{0}\big(z,S_{2}(z)\big)=\sup_{\alpha} \!\!\!\sum_{\tz\sim \ttz, \hspace{0.1cm} d(z,\ttz)=2} \!\!\!\!\!\!\!\!\!\!\!\!\alpha(\tz)^{2}= \sup_{\alpha} \sum_{\tz,\hspace{0.1cm} \tz\sim z}\!\!\!(\textnormal{deg}(z^{\prime})-1)\alpha(z^{\prime})^{2}=\max_{z^{\prime},\hspace{0.1cm} z^{\prime}\sim z}\textnormal{deg}(z^{\prime})-1 \hspace{0.1cm} ,
\end{equation*}
and given that Theorem \ref{thmprinc} ensures that the entropic curvature is lower bounded by $\inf_{z\in \X} r(z)$ with 
$r(z)=-2\log K_{0}(z,S_{2}(z))\geq  2\big(1-K_{0}(z,S_{2}(z))\big)$, the conclusion follows.

{\bf Comments:}
\begin{enumerate}[(i)]

 \item Let us observe that the fact that $\displaystyle \big(2-\max_{z^{\prime},\hspace{0.1cm} z^{\prime}\sim z}\textnormal{deg}(z^{\prime})\big)\leq 0 $ is consistent to the geometry of the underlying generic geodetic graph. Similar considerations have been taken into account for the Ollivier curvature obtaining that trees reach some negative lower bound with respect to this curvature (see \cite[Theorem 2, Proposition 2]{JL14}).

 \item There are other graphs, which are non geodetic  whose lower entropic curvature bound given by $\inf_{z\in \X} r(z)$ is non-positive.  For instance, the hexagonal tiling of the plane is not a geodetic graph, however locally it looks like a 3-regular tree and therefore  $r(z)\geq 2\big(1-\sup_{z\in \X}K_{0}(z,S_{2}(z))\big)=-2$.
 
\begin{figure}[!h]
\centering
\includegraphics[scale=0.2]{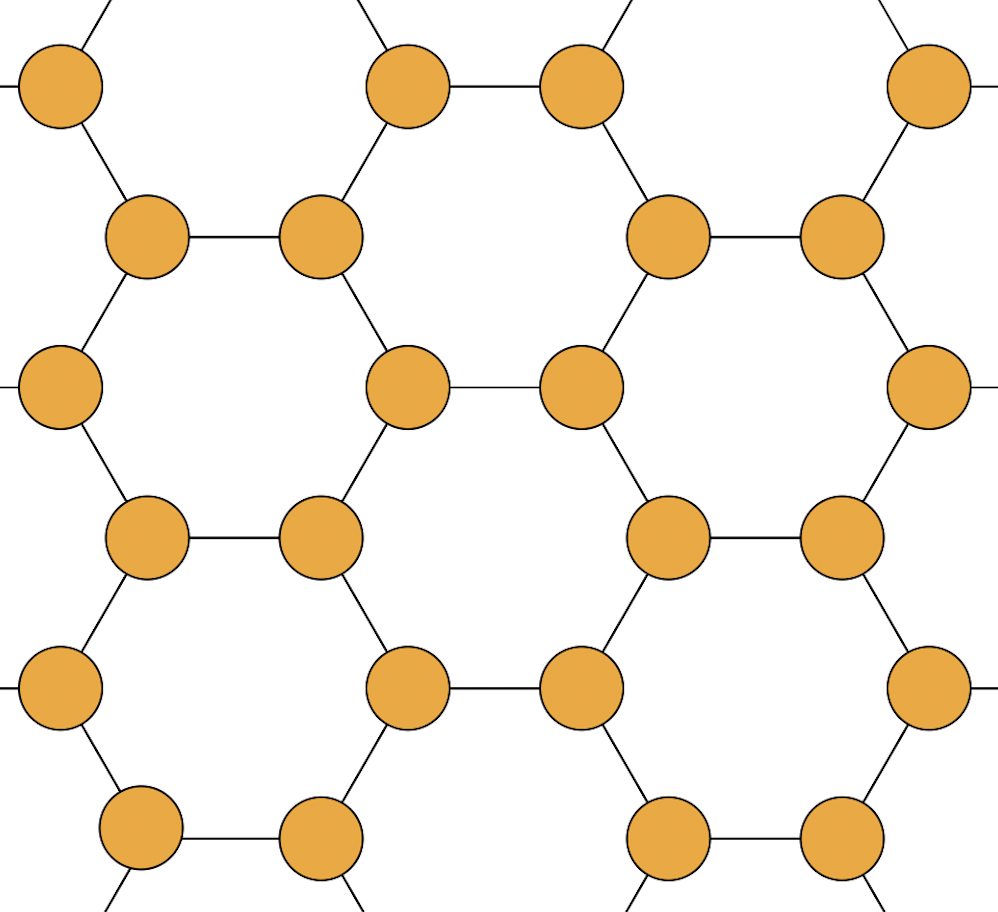}
\caption{Figure of the hexagonal tiling of the plane. } 
\end{figure}
\FloatBarrier

\end{enumerate}

\section{Study of the main criteria}\label{studycrit}
\subsection{Tensorization properties}
In this part, we will study the \textit{tensorization properties} of the constants $r(z)$ and $\widetilde r_2(z)$, $z\in \X$, given by \eqref{defr} and \eqref{defr_3}  with respect to the usual \textit{Cartesian product of graphs}. Recall that from the comments of Theorem \ref{thmprinc} and from Theorem \ref{Thmstructure}, the quantity  $\inf_{z\in\X} r(z)$ is known to be  a lower bound on the entropic curvarture $\kappa$, and $\inf_{z\in\X} \widetilde r_2(z)$ a lower bound on the  $\widetilde T_2$-entropic curvature $\widetilde\kappa_2$. 

For the sake of completeness, let us recall the standard definition of Cartesian product of graphs. The \textit{Cartesian product} of two graphs $(\mathcal{X}_1, E_{1})$ and $ (\mathcal{X}_2, E_{2})$ endowed with distances $d_1$ and $d_2$ respectively is a graph 
\begin{equation*}
(\X,E)=(\mathcal{X}_1, E_{1})\square (\mathcal{X}_2,E_{2}) ,
\end{equation*}
where $\X=\mathcal{X}_1 \times \mathcal{X}_2 $ and the set of edges $E$ is defined by 
\[ (x_1,x_2)\sim (x^\prime_1, x^\prime_2) \hspace{0.2cm } \text{if} \hspace{0.2cm}
\begin{cases}
    \text{either} \hspace{0.2cm} x_1\sim x^\prime_1 \hspace{0.2cm}  \text{and} \hspace{0.2cm} x_2=x^\prime_2 \hspace{0.06cm} , \\
    \text{or}  \hspace{0.2cm}  x_2\sim x^\prime_2  \hspace{0.2cm} \text{and}  \hspace{0.2cm} x_1=x^\prime_1 \hspace{0.1cm}.
\end{cases}
\]
As a consequence if $d_1$, respectively $d_2$, denotes the graph distance on the graph $(\mathcal{X}_1, E_{1})$, respectively $(\mathcal{X}_2, E_{2})$, then the graph distance $d:=d_1\square d_2$ on $(\mathcal{X}, E)$ is given by 
\[d\big((x_1,x_2),(y_1,y_2)\big):=d_1(x_1,y_1)+d_2(x_2,y_2), \qquad  (x_1,x_2),(y_1,y_2)\in \X.\]

Similarly, one defines the  \textit{Cartesian product} of  two graph spaces $(\mathcal{X}_1,d_1,m_{1},L_{1})$, $(\mathcal{X}_2,d_2,m_{2},L_{2})$ the graph space 
\[(\mathcal{X}_1,d_1,m_{1},L_{1})\square(\mathcal{X}_2,d_2,m_{2},L_{2}):=(\X_1\times \X_2, d_1\square d_2, m_1\times m_2, L_1\oplus L_2),\]
where the generator $L=L_1\oplus L_2$ on $\X=\X_1\times \X_2$ is given by: for $z=(z_1,z_2),z'=(z_1',z_2')\in \X$
\[L(z,z^{\prime}) := \hspace{0.2cm }\begin{cases}
0\hspace{0.2cm}  \text{if} \hspace{0.2cm} d(z,z')\geq 2,\\
   L_{1}(z_1,z^{\prime}_1) \hspace{0.2cm}  \text{if} \hspace{0.2cm} z_1\sim z^{\prime}_1 \hspace{0.2cm} \text{and} \hspace{0.2cm}   z'=(z_1',z_2) \hspace{0.06cm} ,\\
   L_{2}(z_2,z_2^{\prime}) \hspace{0.2cm}   \text{if} \hspace{0.2cm} z_2\sim z_2^{\prime} \hspace{0.2cm} \text{and} \hspace{0.2cm}   z'=(z_1,z_2') \hspace{0.06cm} ,\\
 -\Big( \sum_{z_1',z_1'\sim z_1} L_{1}(z_1,z^{\prime}_1) +\sum_{z_2', z_2'\sim z_2} L_{2}(z_2,z_2')\Big)\hspace{0.2cm} \text{if} \hspace{0.2cm}   z=z^{\prime} \hspace{0.06cm}.
\end{cases}\]
Since for $i=1,2$ the measures $m_i$ is reversible with respect to the generators $L_i$, the product measure $m$ is also reversible with respect to $L$. This definition can be iterated to define the Cartesian product of  a finite sequences of graph spaces. 

\begin{theorem}\label{ttens}
Let $(\X, d,m,L)$ be the cartesian product of $n$ graph spaces $(\X_{i},d_i,m_{i},L_{i}), i\in [n]$,
\[(\X, d,m,L):=\big( \mathcal{X}_{1}\times\cdots\times \X_n,d_1\square\cdots  \square d_n ,m_1\otimes\cdots\otimes m_n, L_1\oplus\cdots \oplus  L_n\big).\]
For $z=(z_1.\ldots ,z_n)\in \X$, let $r(z)$ denotes the quantity \eqref{defr} defined on the space $(\X, d,m,L)$ and for $i\in [n]$, let $r(z_i)$ be the same quantity defined on the space $(\X_i, d_i,m_i,L_i)$. Identically one denotes $\widetilde r_2(z)$, $\widetilde r_2(z_i), i\in [n]$ the quantities whose definition is given by \eqref{defr_3}. If $\min \big({r}(z_{1}),\ldots ,{r}(z_{n})  \big)\leq 0$, then 
\[r(z)\geq \min \big({r}(z_{1}),\ldots ,{r}(z_{n})  \big),\]
and if $\min \big({r}(z_{1}),\ldots ,{r}(z_{n})  \big)\geq 0$
then 
\[r(z)\geq -2\log \Big(1-\frac1n\big(1-e^{-\min({r}(z_{1}),\ldots ,{r}(z_{n}))/2}\big)\Big)\geq \frac1n \min \big({r}(z_{1}),\ldots ,{r}(z_{n})  \big).\]
We also have 
\[\widetilde r_2(z)\geq \min \big(\widetilde{r}_{2}(z_{1}),\ldots ,\widetilde{r}_{2}(z_{n})  \big).\]
\end{theorem}
The proof of this theorem is postponed in Appendix B.

{\bf Comments:} If $\X_1=\cdots=\X_n=\{0,1\}$ is the two points space equipped with the counting measure $m_1=\cdots=m_n=m_0$ then the graph space $(\X,d,m,L)$ is the discrete hypercube studied in Section \ref{sectionhypecube} equipped with the counting measure also, for which we prove that for any $z\in \X$
\[r(z)=-2\log(1-1/n)\quad   \mbox{for} \;n\geq 2, \quad \widetilde r_2(z)=1 \quad\mbox{for} \;n\geq 1,\] 
and $r(z)=+\infty$ for $n=1$. 
It follows that the results of Theorem \ref{ttens} can not be improved in full generality.     

\subsection{Geometric conditions for positive entropic curvature - the Motzkin-Strauss Theorem. }\label{sectionpos}

In positive curvature midpoints spread out. These considerations have been taken into account for the Ollivier's coarse curvature (see \cite{Oll11,OllVill12}). The following proposition shows that this property is also a necessary condition so that the lower bound $\inf_{z\in \X} r(z)$ of the entropic curvature of the graph space $(\X,d,L_0,m_0)$, due to Theorem \ref{thmprinc}, is positive. 

\begin{proposition}\label{distanciados}
Let $G=(\X,E)$ be a graph and let $z\in \X$ and let us suppose that $r(z)>0$ where $r(z):=-2 \log K_0(z,S_2(z))$. Then the following properties hold :
\begin{enumerate}
    \item for all $W\subset S_{2}(z)$ with $|W|\leq 2$ one has 
$|]z,W[|>|W|$. 
\item for all $W\subset S_{2}(z)$ satisfying $]z,\ttz[=]z,w''[$ for all $\ttz,w''\in W$, one has  $|]z,W[|>|W|$.
\end{enumerate}
\end{proposition}

We already have proved in comments \ref{distance1} of Theorem \ref{thmprinc} that Proposition \ref{distanciados} holds for $|W|=1$. For  $W:=\{z_{1}^{\prime\prime},z_{2}^{\prime\prime}\}$, if $\big|]z,W[\big|=2$ with $r(z)>0$ then necessarily
$]z,W[=]z,z_{1}^{\prime\prime}[=]z,z_{2}^{\prime\prime}[=\{z_{1}^{\prime},z_{2}^{\prime} \}$ and therefore 
\[K_0\big(z,S_{2}(z)\big)\geq  \sup_{\alpha(z_{1}^{\prime})+\alpha(z_{2}^{\prime})=1}
\big\{2\alpha(z_{1}^{\prime})\alpha(z_{2}^{\prime})+2\alpha(z_{1}^{\prime})\alpha(z_{2}^{\prime})\big\}= 1.\]
This is a contradiction to the assumption that $r(z)>0$ and thus  $\big|]z,W[\big|>2$. The proof of the last point of Proposition \ref{distanciados} is similar. Let $W\subset S_{2}(z)$ satisfying $]z,\ttz[=]z,w''[$ for all $\ttz,w''\in W$, then $]z,\ttz[=]z,W[$ for all $\ttz\in W$ and therefore  
\[K_0\big(z,S_{2}(z)\big)\geq  \sup_{\alpha, \sum_{\tz\in S_1(z)} \alpha(\tz)=1}
\Big\{\sum_{\ttz\in W} \big|]z,W[\big| \prod_{\tz\in ]z,\ttz[} \alpha(\tz)^{\frac2{\big|]z,W[\big|}}\Big\}\geq \frac{|W|}{\big|]z,W[\big|},\]
where for the last inequality we choose $\alpha(\tz)=\frac1{\big|]z,W[\big|}$ for all $\tz \in ]z,W[$. As a consequence $r(z)>0$ implies $K_0\big(z,S_{2}(z)\big)<0$ and therefore $|W|<\big|]z,W[\big|$.

Actually, a natural guess for positive entropic curvature is the following one. 

\begin{conjecture}\label{conj1} Let $G=(\X,E)$ be a graph endowed with the counting measure. If for all  $z\in \X$ and for all $W\subset S_{2}(z)$, 
$
\big|]z,W[\big| > \big|W \big| , 
$
then the graph space $G$ has positive entropic curvature.
\end{conjecture}
The following two remarks show that the relationships between the cardinality of the midpoints and the lower bound on the entropic curvature given by $\inf_{z\in \X} r(z)$ are  subtle.
\begin{itemize}

\item  Given $z\in \X$, $r(z)>0$ does not imply that  for all $W\subset S_{2}(z), \big|]z,W[\big|>\big|W\big|$. 
Indeed, assume that  the graph $G=(V,E)$ restricted to the ball $B_2(z)$ is given by
$S_1(z):=\{z_{1}^{\prime},z_{2}^{\prime},z_{3}^{\prime}\}$ and $S_2(z):=\{z_{1}^{\prime\prime},z_{2}^{\prime\prime},z_{3}^{\prime\prime}\}$ with $]z,z''_1[:=\{z_1',z_2'\}$, $]z,z''_2[:=\{z_2',z_3'\}$ and $]z,z''_3[:=\{z_1',z_3'\}$. 
Then it holds
    $K_0(z,S_{2}(z))=\sup_{\alpha}\big\{2\alpha(z_{1}^{\prime})\alpha(z_{2}^{\prime})+2\alpha(z_{2}^{\prime})\alpha(z_{3}^{\prime})+2\alpha(z_{1}^{\prime})\alpha(z_{3}^{\prime})\big\}=\frac{2}{3}<1$,
and however for $W:=S_2(z)$ one has $\big|W\big|=\big|]z,W[\big|$.

\item  The following example shows that it is possible that the assumption of the conjecture \ref{conj1} holds for a fixed vertex $z_0$ and that $r(z_0)<0$. 
In fact, one will construct a family of balls $B_2(z_0)$  indexed by $n\in \mathbb{N}^*$ for which for all $W\subset S_{2}(z_0)$,
$\big|]z_0,W[\big| > \big|W \big|$ and show that for sufficiently large 
$n$, $r(z_0)<0$. 
 Let 
$
 S_{1}(z_0):=\{z_{1}^{\prime},z_{2}^{\prime},\ldots, z_{n+2}^{\prime}\}$ and $S_{2}(z_0):=\{z_{1}^{\prime\prime},z_{2}^{\prime\prime},\ldots ,z_{n}^{\prime\prime}\}$ with
 for all $i\in [n]$, 
$]z_0,z_{i}^{\prime\prime}[=\{z_{1}^{\prime},z_{2}^{\prime},z_{i+2}^{\prime}\}$. 

\vspace{- 0,5 cm}
\begin{figure}[!h]
\centering
\includegraphics[scale=0.4]{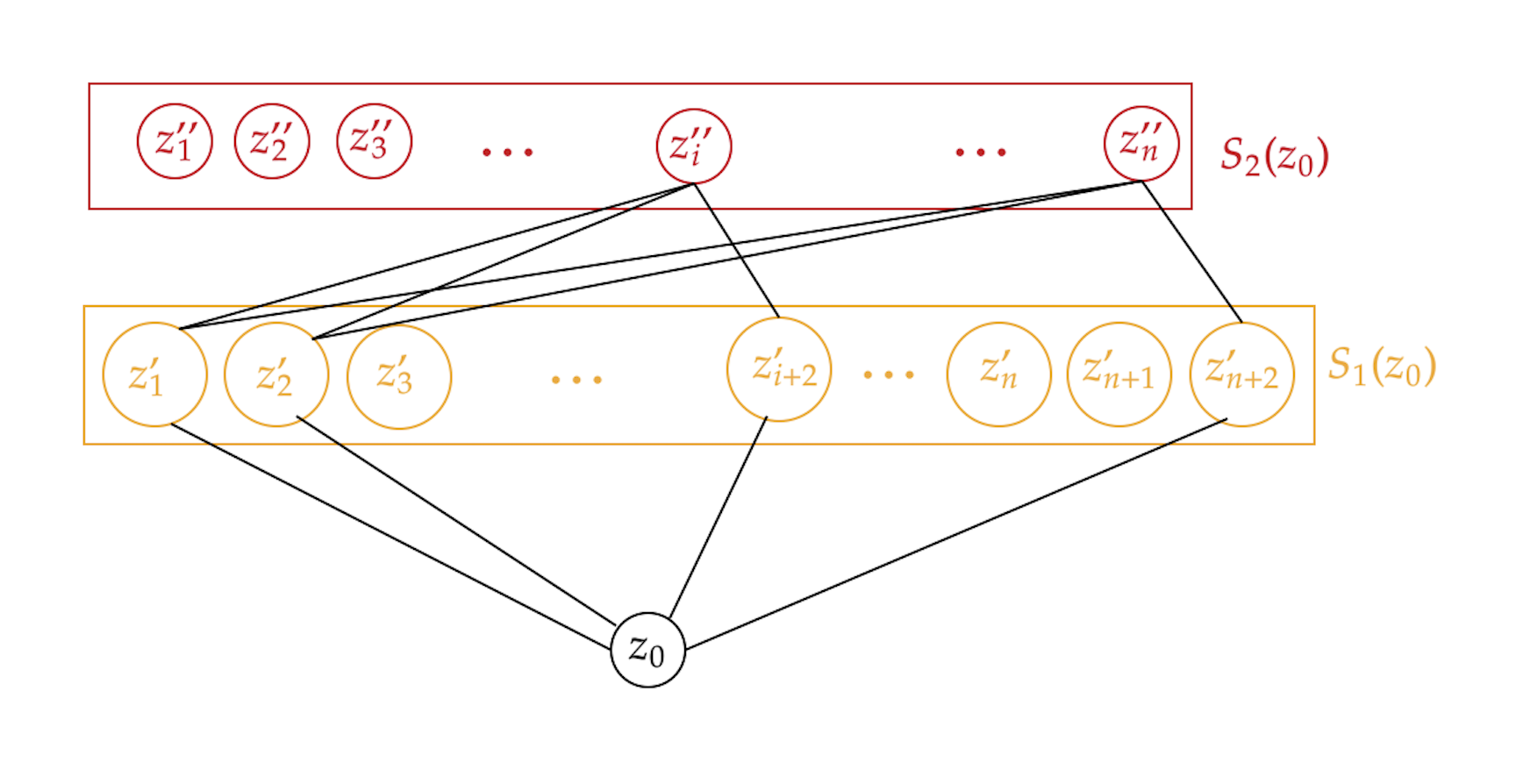}
\end{figure}
\FloatBarrier
\vspace{- 0,5 cm}
One easily check that for all $W\subset S_2(z_0)$, $|]z_0,W[|=|W|+1>|W|$.
Moreover, it holds
\[
K_0\big(z_0,S_{2}(z_0)\big)=\sup_{\alpha, \sum_{i=1}^{n+2}\alpha(z_{i}^{\prime})=1}\sum_{i=1}^{n} 3\big(\alpha(z_{1}^{\prime})\alpha(z_{2}^{\prime})\alpha(z_{i+2}^{\prime})\big)^{\frac{2}{3}} 
\geq 3n\big(\frac{1}{16}\cdot \frac{1}{2n}\big)^{\frac{2}{3}}=\frac{3}{32^{\frac{2}{3}}}n^{\frac{1}{3}} 
\]
where for the last inequality  we  choose 
$\alpha(z_{1}^{\prime})=\alpha(z_{2}^{\prime})=\frac{1}{4}$ and $\alpha(z_{i}^{\prime})=\frac{1}{2n}$ for $i\in \{3,4,\ldots, n+2\}$. 
It follows that for $n\geq 38$, $K_0\big(z_0,S_{2}(z_0)\big)>1$ and therefore 
$r(z_0)<0$.
\end{itemize} 

\begin{remark}
Note that the above construction is not a counterexample of conjecture \ref{conj1} since we only assume that the hypotheses of the conjecture holds for a single $z_0$, and also since  our criteria only gives  a lower bound on the entropic curvature. 
\end{remark}

Let us now propose   sufficient geometric conditions for positive entropic curvature related to the criteria $\inf_{z\in \X} r(z)<0$. One of the key ingredient of the results given below is the so-called \textit{Motzkin Strauss Theorem}. In a seminal 1965 paper \cite{MS65}, Motzkin and Straus found an elegant relationship between the maximum clique of a graph and the global maxima of a quadratic optimization problem defined on the standard simplex. This connection produced another new proof of Turán's theorem \cite{MS65}. The Motzkin Strauss Theorem allows us to interpret the optimization problem which defines $K_{0}\big(z,S_{2}(z)\big)$ for some class of graphs 
as a problem of finding the maximum clique in a related graph. In order to be precise, one needs to introduce certain preliminary definitions.
\begin{definition}[Clique, maximum clique and clique number] Given a simple undirected graph $G^*=(V^*,E^*)$,
\textit{a clique}  is a subset of vertices  such that every two distinct vertices of this subset are adjacent. In other words, a clique  is an induced subgraph of the graph $G^{*}$ that is complete. A \textit{maximum clique} of  $G^{*}$ is a clique with maximum cardinality. This maximum cardinality, denoted by $\omega(G^{*})$, is called the \textit{clique number} of  $G^{*}$.
\end{definition}

\begin{theorem}\label{cliquethm}\cite{MS65} Let $G^{*}=(V^{*},E^{*})$ be a simple undirected graph with clique number $\omega(G^{*})$.
Then the following relation holds,
\begin{equation*}
    2\max_\alpha\big \{\sum_{\{i,j\}\in E^{*}}\alpha_{i}\alpha_{j}\big\} =1-\frac{1}{\omega(G^{*})}\hspace{0.2cm}, 
\end{equation*}
where the supremum runs overs all vectors $\alpha=(\alpha_i)_{i\in V^*}$ with non negative coordinates and such that $\sum_{i\in V^*}\alpha_{i}=1$. Moreover the maximum is achieved by a characteristic vector $\alpha$ of a maximum clique $S$ of the graph $G^{*}$, that is : $\alpha_{i}=\frac{1}{\mid S\mid}$ for $i\in S$ and $\alpha_i=0$ otherwise.
\end{theorem}

Let  $\mathcal{C}$ be  the class of graphs $G=(\X,E)$ satisfying that any pair of vertices in $\X$ at distance two share two midpoints and two midpoints can not be shared by more than two vertices. Let $G$ be a graph of $\mathcal C$, and $z$ be a fixed vertex $z\in \X$. Let $G_z^*=(V_z^*,E_z^*)$ be the graph with set of vertices $V_z^*=S_1(z)$ and such that $\{z',w'\} $ is an edge of $E_z^*$ if $\{z', w'\}=]z,\ttz[$ for some $\ttz\in S_2(z)$. According to this construction, one exactly has
\[ K_{0}\big(z,S_{2}(z)\big)=\sup_{\alpha}\Big\{2\sum_{\{z',w'\}\in E_z^*} \alpha(z')\alpha(w')\Big\}.\]
As a consequence,  the next result is an easy application of Motzkin Strauss Theorem \ref{cliquethm} together with Theorem \ref{thmprinc}.
\begin{proposition}\label{graphMS}
Let $G=(\X,E)$ be a graph belonging to $\mathcal C$ and let $z$ be an
arbitrary vertex of $\X$. Then, one has  $K_{0}\big(z,S_{2}(z)\big)=1-\frac{1}{\omega(G^{*}_{z})}<1$, 
where $\omega(G^{*}_{z})$ is the clique number of the graph $G_z^*$ as defined above. As a consequence the entropic curvature of the graph space $G$ endowed with the counting measure is positive bounded from below by 
\[\inf_{z\in \X} r(z)=\inf_{z\in \X} - 2\log \Big(1-\frac{1}{\omega(G^{*}_{z})}\Big)\geq \frac 2{\sup_{z\in \X}\omega(G^{*}_{z})}.\]
\end{proposition}

 {\bf Comments:}
\begin{itemize}
\item The following figure illustrates the construction of $G_z^*$ for a generic $G\in \mathcal C$:  
\begin{figure}[!ht]
\centering
\includegraphics[scale=0.25]{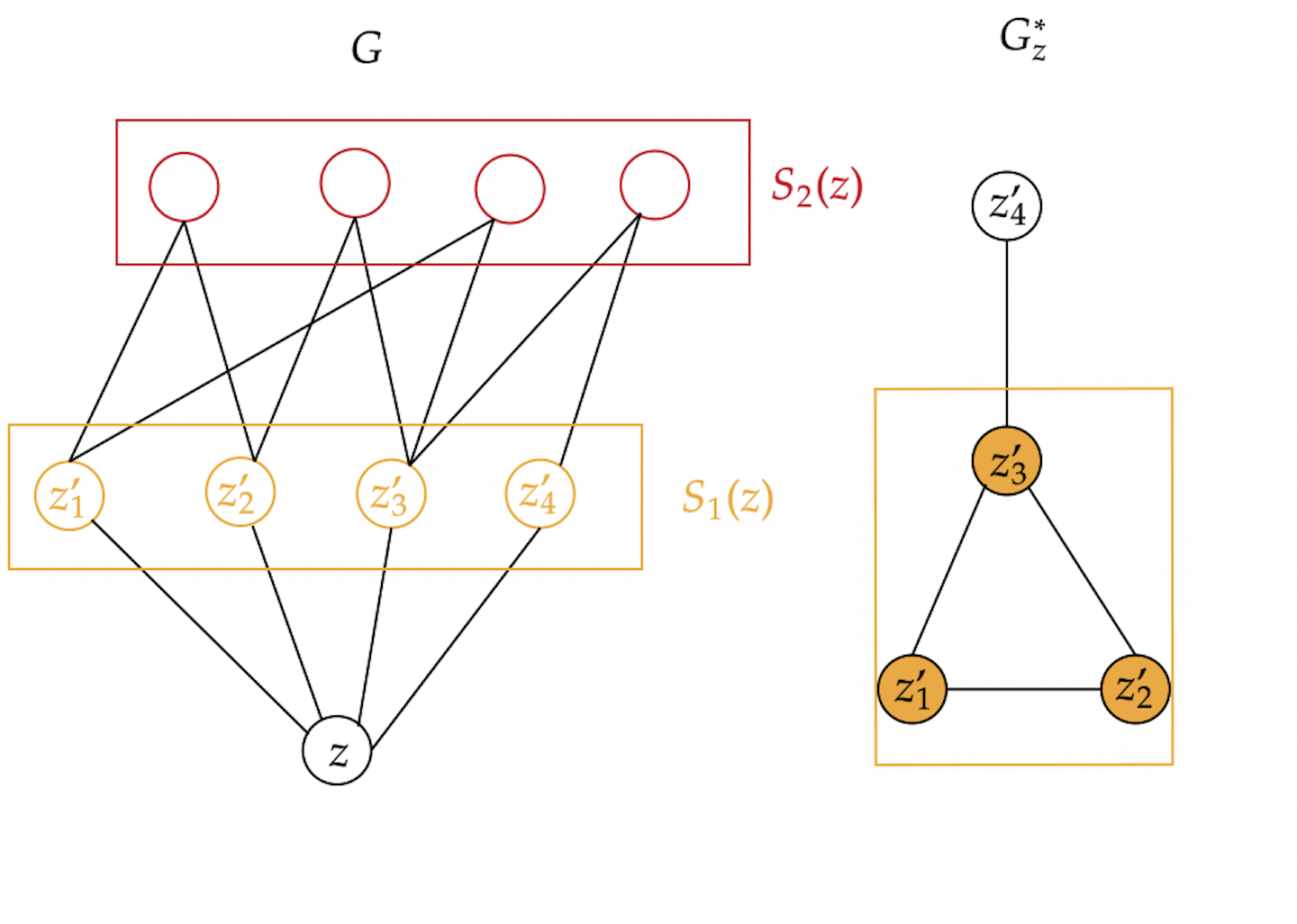}
\end{figure}
\FloatBarrier
\vspace{-0,5 cm}
For this drawing example, one has
$K_0(z,S_2(z))= 1-\frac1{\omega(G_z^*)}=\frac23$.

\item  As an example, the hypercube $\X=\{0,1\}^n$  belongs to the class of graphs  $\mathcal{C}$ and and for any $z\in \X$, $G_z^*$ is the complete graph on $n$ vertices $K_n$ with clique number $n$. One recovers the lower bound $-2\log(1-1/n)$ on the entropic curvature given in Section \ref{sectionhypecube}.

\item From a complexity point of view, the problem of computing the clique number is one of Karp’s 21 NP-hard problems \cite{Karp}. Thus, it is immediate that the problem of calculating the entropic curvature of a graph is an NP-hard problem.

\end{itemize}

In the next proposition we consider another class of graphs  
satisfying the assumptions of  Conjecture \ref{conj1}, together with a covering condition. For this class of graphs we also derive positive entropic curvature by applying the Motzkin Strauss Theorem. 

\begin{proposition}\label{Dieppe}  Let $G=(\X,E)$ be a graph. Assume that for any  arbitrary vertex $z$ of $\X$ and for any three distinct vertices $z_1^{\prime\prime},z_2^{\prime\prime},z_3^{\prime\prime}$ of $S_2(z)$,
\begin{equation}\label{Hypo1}
]z,z_1^{\prime\prime}[\cap]z,z_2^{\prime\prime}[\cap ]z,z_3^{\prime\prime}[=\emptyset.
\end{equation}
Assume also that for all subsets $W\subset S_2(z)$ of cardinality one or two, one has
\begin{equation}\label{Hypo2}
\big|]z,W[\big| >|W|.
\end{equation}
Then one gets $K\big(z,S_2(z)\big)\leq 7/8$ and Theorem \ref{thmprinc} ensures that the entropic curvature of the graph space $G$ endowed with the counting measure is bounded from below by $1/4$.
\end{proposition}

\begin{remark}
Note that the hypothesis \eqref{Hypo1} implies that 
there is no overlapping of more than three midpoint sets. 
It is still an open problem to generalize this type of results to larger overlapping.
\end{remark}

\section{Few comparisons with other notions of curvature.}\label{sectcomp} 
Recall that according to the comments of Theorem \ref{thmprinc}, the entropic curvature of the graph space $G=(\X,d,m_{0},L_{0})$ is lower bounded by $\inf_{z\in \X} r(z)$ where $r(z):=-2\log K_0(z,S_2(z)$ is interpreted as a local lower bound on the entropic curvature. In this section we give  some comparative remarks between this lower bound and the notions of curvature by  \textit{Lin-Lu-Yau} and  \textit{Bakry-{\'E}mery}.\\

\subsection{Entropic curvature and the Lin-Lu-Yau curvature.}\
The Lin-Lu-Yau curvature is a modified notion of the coarse Ollivier's Ricci curvature introduced by Lin, Lu and Yau in \cite{LLY11}.  In \cite{MW19}, Florentin M{\"u}nch and Rados{\l}aw K.Wojciechowski, generalized the notion of Lin-Lu-Yau curvature  for any graph Laplacian.
\begin{definition}[Lin-Lu-Yau Ricci curvature]
Given $G=(\mathcal{X},E)$ a graph endowed with its graph distance $d$ and with a Markov chain defined by $m:=\{m_{z}(\cdot)\}_{z\in \mathcal{X}}$. 
For $0\leq \alpha<1$, the \textit{$\alpha$-lazy random walk} $m_{x}^{\alpha}$ associated to the graph Laplacian $L_{0}$ is defined as
\begin{equation*}
m_{x}^{\alpha}(y)=
\begin{cases}
\frac{\alpha}{\Delta(G)} \hspace{0.2cm} \text{if} \hspace{0.2cm} x\sim y ,\\
 1- \alpha \frac{\text{deg
 }(x)}{\Delta(G)}
\hspace{0.2cm} \text{if} \hspace{0.2cm} y= x ,\\
0 \hspace{0.2cm} \text{otherwise}.
\end{cases}
\end{equation*}

For every $x,y\in \mathcal{X}$, one defines
\begin{equation*}\kappa_{\alpha}(x,y):=1-\frac{W_{1}(m_{x}^{\alpha},m_{y}^{\alpha})}{d(x,y)}.
\end{equation*}

As shown in \cite{MW19}, the limit as $\alpha\rightarrow 0$ exists for any graph Laplacian and therefore one can define the Ricci Lin-Lu-Yau curvature along the edge $\{x,y\}$ denoted as $\kappa_{LLY}(x,y)$ by
\begin{equation*}
\kappa_{LLY}(x,y):=\lim_{\alpha \rightarrow 0} \frac{\kappa_{\alpha}(x,y)}{\alpha} .
\end{equation*}
\end{definition}

The next proposition establishes a link between the Lin-Lu-Yau curvature and  the graph-theoretical notion of \textit{girth}. 

\begin{definition}
The girth of a graph $G=(\mathcal{X},E)$, denoted $g(G)$ is the length of the shortest cycle contained in $G$. Acyclic graphs are considered to have infinite girth.
\end{definition}
Adapting the proof \cite[Theorem 2.b(ii)]{HS13}, provides the following proposition with the measures $m_{x}^{\alpha}$ associated to the generator $L_{0}$. Its proof is postponed in Appendix B.

\begin{proposition}\label{propcycles}
Let $G=(\mathcal{X},E)$ be a graph. If for all $x,y\in \mathcal{X}$ with $d(x,y)=1$ 
\begin{equation*}
\kappa_{LLY}(x,y)<\frac{6-\text{deg}(x)-\text{deg}(y)}{\Delta(G)}
\end{equation*}
then $g(G)\geq 5$. 
\end{proposition}

\begin{remark}\label{girth2}
Let  $G=(\mathcal{X},E)$ be a graph . If $g(G)\geq 5$, then there cannot be two midpoints between two vertices at distance two and thus as already noted in the introduction for all $z\in \mathcal{X}$, $r(z)\leq 0$.
\end{remark}
Thanks to the above remark, by contraposition, we immediately obtain the following corollary.
\begin{corollary}
Let $G=(\mathcal{X},E)$ be a graph. If for all $z\in \X$,
    $r(z)> 0$,
then for all $x,y\in \mathcal{X}$ with $d(x,y)=1$,
\begin{equation*}
\kappa_{LLY}(x,y)\geq \frac{6-\text{deg(x)}-\text{deg(y)}}{\Delta(G)} \hspace{0.1cm} .
\end{equation*}
\end{corollary}

Observe that  $\kappa_{LLY}(x,y)>0$  does not imply $r(x)>0$ nor that $r(y)>0$. Indeed, let us consider the so-called  \textit{windmill graph} $W_{d}(4,2)$, consisting of 2 copies of the complete graph $K_{4}$ at a shared universal vertex:
\vspace{-0,7 cm}
\begin{figure}[!ht]
\centering
\includegraphics[scale=0.2]{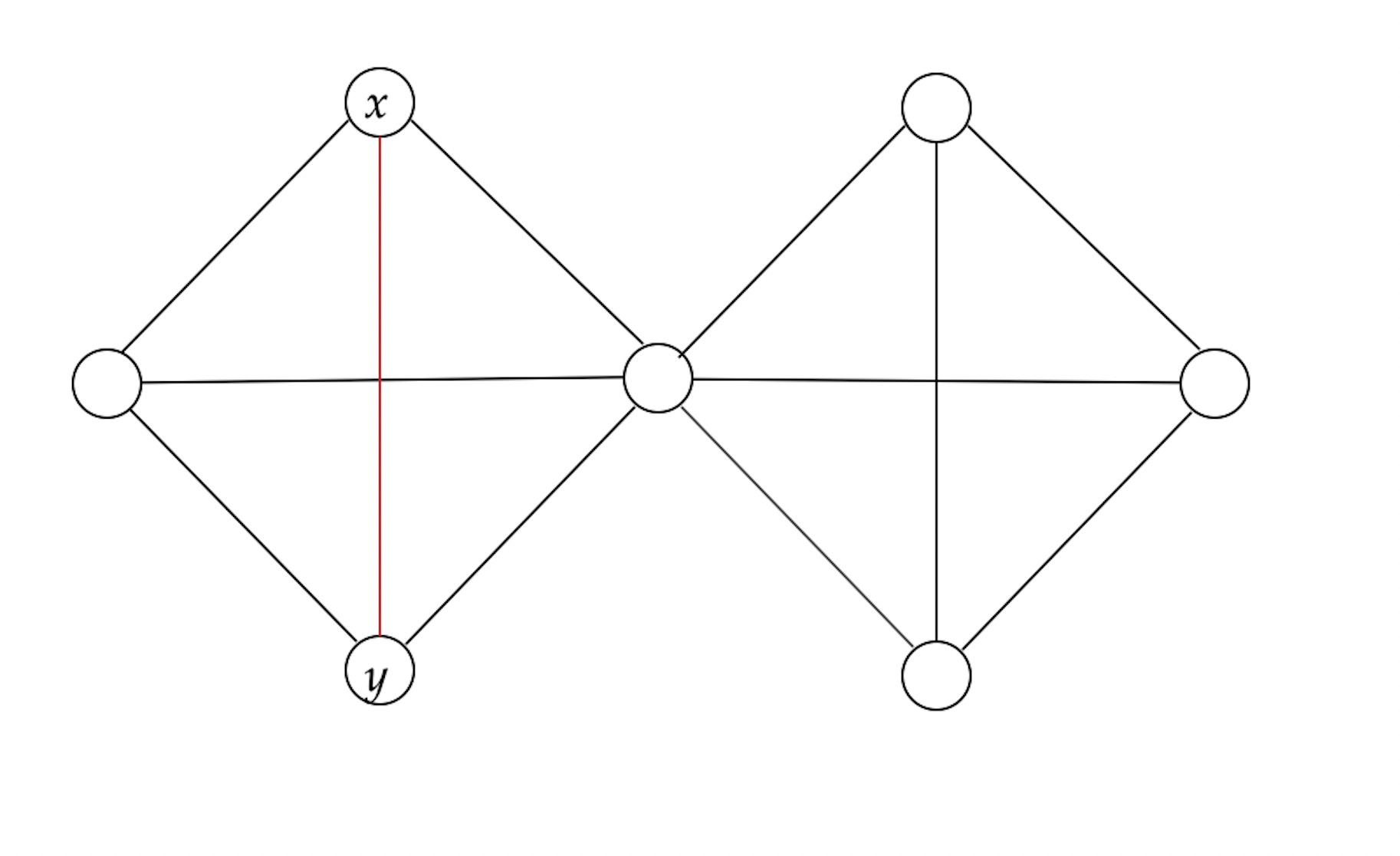}
\label{Figure1.1}
\end{figure}
\FloatBarrier
\vspace{-1,5 cm}
 For the edge $\{x,y\}$ of the graph $W_{d}(4,2)$ on the  figure, $\kappa_{LLY}(x,y)=2/3$ (one may easily check that 
$W_{1}(m_{x}^{\alpha},m_{y}^{\alpha})=1-\frac{\alpha}{2}-\frac{\alpha}{6}$). 
However, one has $K_0\big(x,S_2(x)\big)=K_0\big(y,S_2(y)\big)=3$ and thus $r(x)=r(y)<0$.\\

\subsection{Entropic curvature and the Bakry-Émery curvature condition.}\

There are  relationships between the entropic curvature and the \textit{Bakry-Émery curvature}. The notion of Bakry-Émery curvature was first introduced by Bakry and \'Emery in \cite{BE85}. The Bakry-Émery curvature is motivated by the Bochner's identity in Riemannian Geometry and has been extensively studied in discrete spaces recently  \cite{CKL21,CLP20}. Let us make some qualitative  remarks on the similarities with respect to the local structure and the negativity of curvature for both notions.\

In the case where $L=L_{0}$, the works  \cite{CKL21,CLP20} show that the Bakry-Émery conditions are also related to the local structure of balls of radius 2. More precisely, \textit{the curvature matrix} for a vertex $z\in \mathcal{X}$ is completely determined by \textit{the incomplete ball of radius 2 around $z$}, that is, the graph induced by $B_{2}(z)$ removing all edges connecting vertices within $S_{2}(z)$ {(see \cite[Remark 2.2]{CKL21})}. Similarly, the lower bound $r(z)$ interpreted as the local entropic curvature at vertex $z$ only depends on {this incomplete ball of radius 2} removing also all edges connecting vertices within $S_{1}(z)$. 

Moreover, according to \cite[Theorem 6.4]{CLP20},   if \textit{the punctured 2-ball  around  $z$},  
defined as 
the incomplete ball of radius 2  from which we remove all edges connected to $z$,
has more than one connected component then the Bakry-{\'Emery}  curvature criterion 
at the vertex $z$ is negative with five exceptions (see \cite[Theorem 6.4]{CLP20}). In this configuration, choosing $z'\in S_1(z)$  in one component of the punctured 2-ball, the quantity $r(z')$ is less or equal to zero. Indeed, since the punctured 2-ball around $z$ has more than one connected component there exists $w'\in S_2(z')\cap S_1(z)$ such that $z\in S_1(z')$ is the unique midpoint between $z'$ and $w'$.
As an immediate consequence (see \cite[Corollary 6.8]{CLP20}), if a graph has girth greater than or equal to five, then the Bakry-{\'Emery}  curvature criterion at  each vertex is less than zero. Recall that in our setting, girth greater or equal to  five implies $r(z)\leq 0$ for all vertices $z$ (see Remark \ref{girth2}).

\section{Appendix A}

\begin{lemma}\label{lemL} Let $\X$ be  a structured graph associated to a set of moves $\Sc$. 
Then the following properties hold :
\begin{enumerate}[label=(\roman*)]
\item \label{1} Given $d\in \N$, $\tau \in \Sc$ and $\alpha_0,\ldots,\alpha_d\in \X$, if $\big(\alpha_0,\ldots,\alpha_d, \tau(\alpha_d)\big)\in G\big(\alpha_0,\tau(\alpha_d)\big)$ then for any $k\in\{0,\ldots,d\}$ one has $\big(\alpha_0,\ldots , \alpha_k,\tau(\alpha_k),\ldots,\tau(\alpha_d)\big)\in G\big(\alpha_0,\tau(\alpha_d)\big)$.
If moreover the generator $L$  on $\X$ satisfies condition \eqref{condL}, then  one has 
 \[L\big(\alpha_0,\ldots,\alpha_d, \tau(\alpha_d)\big)=L\big(\alpha_0,\ldots , \alpha_k,\tau(\alpha_k),\ldots,\tau(\alpha_d)\big).\]

\item \label{3} For any $x,y,z\in \X$ and $\tau\in \Sc$ with $\tau(z)\sim z$, if  $(z,\tau(z))\in [x,y]$ then  $\tau(x)\in ]x,y]$ and $\tau(z)\in [\tau(x),y]$.

\item \label{6} Let $z\in \X$ and $\tau,\sigma_1,\sigma_2\in \Sc $ such that 
$d\big(z,\tau\sigma_1(z)\big)=2 $ and $\tau\sigma_1(z)=\tau\sigma_2(z)$. Then one has $\sigma_1=\sigma_2$.

\item \label{2}
Let $d\in \N^*$, $\tau \in \Sc$ and $\alpha_0,\beta_2,\ldots \beta_{d+1}\in \X$. If $(\alpha_0, \tau(\alpha_0), \beta_2, \ldots ,\beta_{d+1})\in G(\alpha_0, \beta_{d+1})$ then for any $k\in \{0,\ldots, d\}$, there exists a single $(\alpha_0,\alpha_1,\ldots, \alpha_k)\in \X^{k+1}$ such that \[(\alpha_0,\alpha_1,\ldots, \alpha_k,\tau(\alpha_k), \beta_{k+2},\ldots, \beta_{d+1})\in G(\alpha_0, \beta_{d+1}),\] and for any $\ell\in [k]$, $\beta_{\ell+1}=\tau(\alpha_\ell)$.

\item \label{4} Given $d\in \N^*$, $\tau \in \Sc$ and $\alpha_0\in \X$, let 
\[{\mathcal Y}(\alpha_0,\tau ,d):=\big\{(\alpha_1,\ldots,\alpha_{d+1})\in \X^{d+1}\,\big|\, (\alpha_0,\ldots,\alpha_d,\alpha_{d+1})\in G(\alpha_0,\alpha_{d+1}), \alpha_{d+1}=\tau (\alpha_d)\big\},\]
and 
\[{\mathcal W}(\alpha_0,\tau ,d):=\big\{(\beta_2,\ldots,\beta_{d+1})\in \X^{d}\,\big|\, (\alpha_0,\tau(\alpha_0),\beta_2,\ldots,\beta_{d+1})\in G(\alpha_0,\beta_{d+1})\big\}.\]
The map $\Psi:(\alpha_1,\ldots, \alpha_{d+1})\mapsto (\tau(\alpha_1),\ldots, \tau(\alpha_d))$ is one to one from the set ${\mathcal Y}(\alpha_0,\tau ,d)$  to the set ${\mathcal W}(\alpha_0,\tau ,d)$.
\item \label{5} Assume that $L$ is a generator on $\X$ satisfying condition \eqref{condL}. 
 Let $x,y\in \X$ and $\tau\in \Sc$ such that $\tau(x)\in ]x,y]$ and let  $k\in \{0,\ldots,d-1\}$ where $d=d(x,y)-1$. Then  one has 
 \[\sum_{z\in \X,d(x,z)=k,(z,\tau(z))\in[x,y] } L^k(x,z)L(z,\tau(z))L^{d-k-1}(\tau(z),y)=L(x,\tau(x)) L^{d-1}(\tau(x), y).\]
\end{enumerate}
\end{lemma}
\begin{proof} 
The proof of \ref{1} is by induction over $k\in \{0,\ldots,d\}$. The property holds for $k=d$ by assumption. Assume that for some fixed $k\in [d]$,   $\big(\alpha_0,\ldots , \alpha_k,\tau(\alpha_k),\ldots,\tau(\alpha_d)\big)$ is a discrete geodesic. Then since $d(\alpha_{k-1},\alpha_k)=1$, there exists a single  $\sigma_k\in \Sc$ such that 
$\sigma_k(\alpha_{k-1})=\alpha_k$. Since $d\big(\alpha_{k-1},\tau(\alpha_k)\big)=d\big(\alpha_{k-1},\tau\sigma_k(\alpha_{k-1})\big)=2$, 
according to the definition of structured graphs, there exists $\psi(\sigma_k)\in \Sc_{\tau(\alpha_{k-1})}$ such that 
\[\tau\sigma_k(\alpha_{k-1})=\psi(\sigma_k)\tau(\alpha_{k-1}).\]
Therefore one has $\tau(\alpha_{k-1})\in]\alpha_{k-1},\tau(\alpha_k)[$ and  $(\alpha_0,\ldots , \alpha_{k-1},\tau(\alpha_{k-1}),\ldots,\tau(\alpha_d))$ is a discrete geodesic. If moreover condition \eqref{condL} holds, then by induction hypothesis
\begin{align*}
&L\big(\alpha_0,\ldots , \alpha_{k-1},\tau(\alpha_{k-1}),\ldots,\tau(\alpha_d)\big)\\
&=L\big(\alpha_0,\ldots , \alpha_{k-1}\big)L\big(\alpha_{k-1},\tau(\alpha_{k-1})\big)L\big(\tau(\alpha_{k-1}),\tau\sigma_k(\alpha_{k-1}) \big)L\big(\tau(\alpha_k),\ldots, \tau(\alpha_d)\big)\\
&=L\big(\alpha_0,\ldots , \alpha_{k-1}\big)L\big(\alpha_{k-1},\sigma_k(\alpha_{k-1})\big)L\big(\sigma_k(\alpha_{k-1}),\tau\sigma_k(\alpha_{k-1}) \big)L\big(\tau(\alpha_k),\ldots, \tau(\alpha_d)\big)\\
&=L\big(\alpha_0,\ldots , \alpha_k,\tau(\alpha_k),\ldots,\tau(\alpha_d)\big)\\
&=L\big(\alpha_0,\ldots,\alpha_d, \tau(\alpha_d)\big).
\end{align*}

Item \ref{3} is an easy consequence of \ref{1}. Indeed, if $(z,\tau(z))\in [x,y]$ then, setting $d=d(x,y)$ and $k=d(x,z)$, there exists $(\alpha_0,\alpha_1,\ldots,\alpha_k, \tau(\alpha_k), \beta_{k+2}, \ldots, \beta_d)\in G(x,y)$ with $\alpha_k=z$. Item \ref{1} implies that \[(x,\tau(x),\tau(\alpha_1),\ldots, \tau(\alpha_k), \beta_{k+2}, \ldots, \beta_d)=(\alpha_0,\tau(\alpha_0),\ldots, \tau(\alpha_k), \beta_{k+2}, \ldots, \beta_d)\in G(x,y),\]
and therefore $\tau(x)\in ]x,y]$ and $\tau(z)\in[\tau(x),y]$.

For the proof of \ref{6}, let $z\in\X$ and $\sigma_1,\sigma_2,\tau \in \Sc$ such that $d\big(z,\tau\sigma_1(z)\big)=2$. If $\tau\sigma_1(z)=\tau\sigma_2(z)$ then according to the definition of structured graphs $\psi(\sigma_1)\tau(z)=\psi(\sigma_2)\tau(z)$. It follows that $\psi(\sigma_1)=\psi(\sigma_2)$ and therefore $\sigma_1=\sigma_2$ since the map $\psi$ is one to one. 

The proof of \ref{2} is  by induction over $k\in \{0,\ldots,d\}$. The property holds for $k=0$ by assumption. Assume that for a fixed $k\in\{0,\ldots,d-1\}$, there exists a single $(\alpha_0,\alpha_1,\ldots, \alpha_k)\in \X^{k+1}$ such that 
\[(\alpha_0,\alpha_1,\ldots, \alpha_k,\tau(\alpha_k), \beta_{k+2},\ldots, \beta_{d+1})\in G(\alpha_0, \beta_{d+1}),\] 
and for any $\ell\in [k]$, $\beta_{\ell+1}=\tau(\alpha_\ell)$. According to the definition of structured graphs, since $d(\tau(\alpha_k), \beta_{k+2})=1$, there exists a single $\sigma_{k+1}'\in \Sc_{\tau(\alpha_k)}$ such that $\beta_{k+2}=\sigma_{k+1}'\tau(\alpha_k)$, and since $d\big(\alpha_k,\sigma_{k+1}'\tau(\alpha_k)\big)=2$, one has
\[\beta_{k+2}=\sigma_{k+1}'\tau(\alpha_k)=\tau\psi^{-1}(\sigma_{k+1}')(\alpha_k).\]
Setting $\alpha_{k+1}=\psi^{-1}(\sigma_{k+1}')(\alpha_k)$, one has $\beta_{k+2}=\tau(\alpha_{k+1})$ and since $\alpha_{k+1}\in ]\alpha_k,\beta_{k+2}[$ it follows that  
\[(\alpha_0,\alpha_1,\ldots, \alpha_{k+1},\tau(\alpha_{k+1}), \beta_{k+3},\ldots, \beta_{d+1})\in G(\alpha_0, \beta_{d+1}).\] 
Moreover if $\alpha'_{k+1}$ is such that $\beta_{k+2}=\tau(\alpha_{k+1}')$ and 
\[(\alpha_0,\alpha_1,\ldots,\alpha_k, \alpha_{k+1}',\tau(\alpha_{k+1}'), \beta_{k+3},\ldots, \beta_{d+1})\in G(\alpha_0, \beta_{d+1}),\]
then there exists $\sigma_{k+1}$ such that $\alpha_{k+1}'=\sigma_{k+1}(\alpha_k)$. Applying  \ref{6},
it follows that $\psi(\sigma_{k+1})=\sigma_{k+1}'$ and therefore 
\[\alpha_{k+1}'=\sigma_{k+1}(\alpha_k)=\psi^{-1}(\sigma_{k+1}')(\alpha_k)=\alpha_{k+1}.\]
This ends the proof of \ref{2}.

We now turn to the proof of \ref{4}. Let $d\in \N^*$, $\tau \in \Sc$ and $\alpha_0\in \X$. If $(\alpha_1,\ldots,\alpha_{d+1})\in {\mathcal Y}(\alpha_0,\tau, d)$ then according to \ref{1} with $k=0$, one has 
\[(\alpha_0,\tau(\alpha_0),\ldots,\tau(\alpha_d))\in G(\alpha_0,\tau(\alpha_d)),\]
and therefore $\Psi(\alpha_1,\ldots,\alpha_{d+1})\in {\mathcal W}(\alpha_0,\tau, d)$.
Conversely if $(\beta_2,\ldots,\beta_{d+1})\in {\mathcal W}(\alpha_0,\tau, d)$, then 
\[(\alpha_0,\tau(\alpha_0),\beta_2,\ldots,\beta_{d+1})\in G(\alpha_0,\beta_{d+1}),\]
and according to \ref{2} for $k=d$, there exists a single $(\alpha_0,\ldots, \alpha_d)$ such that 
\[(\alpha_0,\alpha_1,\ldots,\alpha_d, \tau(\alpha_d))\in  G(\alpha_0,\beta_{d+1}),\]
and for all $\ell\in[d]$, $\beta_{\ell+1}=\tau(\alpha_\ell)$.
Therefore, there exists a single  $(\alpha_1,\ldots, \alpha_d)\in {\mathcal Y}(\alpha_0,\tau, d)$ such that $\psi(\alpha_1,\ldots, \alpha_d)=(\beta_2,\ldots ,\beta_{d+1})$.

For the proof of \ref{5}, let $L$ be a generator on $\X$ satisfying condition \eqref{condL}, let $x,y\in \X$ and $\tau\in \Sc$ such that $\tau(x)\in ]x,y]$ and let  $k\in \{0,\ldots,d-1\}$ where $d=d(x,y)-1$. By definition,  one has 
\begin{multline*}
 \sum_{z\in \X,d(x,z)=k,(z,\tau(z))\in[x,y] } L^k(x,z)L(z,\tau(z))L^{d-k-1}(\tau(z),y)\\
 =\sum_{(\alpha_1,\ldots, \alpha_{k+1})\in {\mathcal Y}(x,\tau,k)}\sum_{\gamma\in G(\tau(\alpha_k),y)} L(x,\alpha_1,\ldots, \alpha_{k+1})\, L(\gamma).
\end{multline*}
Applying \ref{1} and then \ref{4} it follows that 
\begin{align*}
 \sum_{z\in \X,d(x,z)=k,(z,\tau(z))\in[x,y] }& L^k(x,z)L(z,\tau(z))L^{d-k-1}(\tau(z),y)\\
 &=\sum_{(\alpha_1,\ldots, \alpha_k)\in {\mathcal Y}(x,\tau,k)}\sum_{\gamma\in G(\tau(\alpha_k),y)} L(x,\tau(x),\tau(\alpha_1),\ldots, \tau(\alpha_k))\, L(\gamma)\\
 &=\sum_{(\beta_2,\ldots,\beta_{k+1})\in {\mathcal W}(x,\tau,k)} \sum_{\gamma\in G(\beta_{k+1},y)} L(x,\tau(x),\beta_2,\ldots, \beta_{k+1})\, L(\gamma)\\
 &=L(x,\tau(x))\sum_{\gamma'\in G(\tau(x),y)} L(\gamma')=L(x,\tau(x))L^{d-1}(\tau(x),y).
\end{align*}
The proof of Lemma \ref{lemL} is completed.
\end{proof}

\begin{lemma}\label{deriveseconde} Let   $(\X,d,m,L)$ be a graph space. Let $v:\X\to \R$ be a bounded function and given  $x,y\in\X$  let 
\[R(t)=\int v \, d\nu_t^{x,y},\qquad t\in (0,1),\]
where $(\nu_t^{x,y})_{t\in [0,1]}$ is the Schr\"odinger path between Dirac measures at $x$ and $y$ defined  by \eqref{pathdirac}. 
One has for any $t\in[0,1]$,
\[R''(t):=
d(x,y)\big(d(x,y)-1\big)\, D_t v(x,y),
\]
with 
\begin{multline*}
  D_t v(x,y):=\sum_{(z,\ttz)\in [x,y], d(z,\ttz)=2} \left[\sum_{\tz\in  ]z,\ttz[} \left(v(\ttz)+v(z)-2v(\tz)\right) L(z,\tz)L(\tz,\ttz) \right] \\r(x,z,\ttz,y)\, \rho_t^{d(x,y)-2}(d(x,z))  \end{multline*}
\end{lemma}

\begin{proof}[Proof of Lemma \ref{deriveseconde}]
Let $d:=d(x,y)$. For $t\in [0,1]$, one has 
\[R(t)=\sum_{k=0}^d\rho_t^d(k) R_{k},\quad\mbox{
with}\quad 
R_{k}:=\sum_{z\in [x,y], d(x,z)=k} v(z) \,\frac{L^{d(x,z)} (x,z)L^{d(z,y)}(z,y)}{L^d(x,y)}\]
Simple computations give for any $t\in [0,1]$,
\[R'(t)=d\,\sum_{k=0}^{d-1}\rho_t^{d-1}(k) \big(R_{k+1}-R_{k}\big),\]
and therefore
\[R''(t)=d(d-1)\,\sum_{k=0}^{d-2}\rho_t^{d-2}(k) \big(R_{k+2}+R_{k}-2R_{k+1}\big).\]
Then the result follows  observing that for any  $k\in\{0,\ldots ,d-2\}$, 
\begin{align*}
&R_{k+2}+R_{k}-2R_{k+1}\\
&=\!\!\!\!\!\!\!\!\! \sum_{(z,\tz,\ttz)\in [x,y], d(x,z)=k} \left(v(\ttz)+v(z)-2v(\tz)\right) \frac{L^{d(x,z)}(x,z)L(z,\tz)L(\tz,\ttz)L^{d(\ttz,y)}(\ttz,y)}{L^{d}(x,y)}\\
&=\!\!\!\!\!\!\!\!\! \sum_{(z,\ttz)\in [x,y] d(x,z)=k, d(z,\ttz)=2} \left[\sum_{\tz\in  ]z,\ttz[} \left(v(\ttz)+v(z)-2v(\tz)\right) L(z,\tz)L(\tz,\ttz) \right] r(x,z,\ttz,y) .
\end{align*}
\end{proof}

\begin{proposition}\label{BE}
 Let $G=(\X,E)$ be a structured graph with  finite set of moves $\Sc$. Let us suppose that moves in $\Sc$ commute, that is, $\sigma\tau=\tau\sigma$ for all
 $\sigma,\tau \in \Sc$. Then the Bakry-{\'E}mery curvature criterion $CD(0,\infty)$ is satisfied for every $z\in \X$.
 \end{proposition}

 Note that the commutativity condition is not necessary for structured graphs  to satisfy the $CD(0,\infty)$  criterion. Indeed, the Bernoulli-Laplace model corresponds to a non-commutative structured graph   with positive Bakry-\'Emery curvature as shown in \cite[Theorem 2.7]{KKRT16}.
\begin{proof}[Proof of Proposition \ref{BE}]
Let us recall the definition of the Bakry-Émery curvature condition in a graph space equipped with the generator $L_0$. Let  $\Gamma$ and $\Gamma_{2}$ be symmetric operators defined respectively as 
\begin{align*}
    2\Gamma(f,g)&:=L_{0}(fg)-fL_{0} g-gL_{0} f,\\
    2\Gamma_{2}(f,g)&:=L_{0}(\Gamma(f,g))-\Gamma(f,L_{0} g)-\Gamma(g,L_{0} f)
\end{align*}
for all real function $f$ and $g$ on $\X$,
where $L_0(f)$ is the discrete Laplace operator. As a convention $\Gamma(f):=\Gamma(f,f)$ and $\Gamma_{2}(f):=\Gamma_{2}(f,f)$.
\begin{definition}\cite[Bakry-Émery curvature condition]{YS10,Sch99}
A vertex $z\in \mathcal{X}$ satisfies the Bakry-Émery curvature-dimension inequality $CD(\kappa_{BE},N)$ with $\kappa_{BE}\in \mathbb{R}$ and $N\in (0,\infty]$, if for any $f:\mathcal{X}\rightarrow \mathbb{R}$
\begin{equation*}
    \Gamma_{2}(f)(z)\geq \frac{1}{N}(L_{0} f(z))^{2}+\kappa_{BE}\Gamma(f)(z),
\end{equation*}
where $N$ is  a dimension parameter and $\kappa_{BE}$ is regarded as a lower Ricci bound at $z\in \X$.  
\end{definition}

 We want to prove that for any vertex $z$ of a structured graph whose moves in $\Sc$ commute,  one has $\Gamma_{2}(f)(z)\geq 0$. Following the same computations as in \cite{CY96,YS10,CKKLP21}, one has 
\begin{align*}
L_0 \Gamma(f)(z)&=\sum_{\sigma \in \Sc,\sigma (z)\sim z}\big( \Gamma(f)(\sigma (z))-\Gamma(f)(z) \big)\\
&=\frac{1}{2} \sum_{\sigma,\sigma(z)\sim z}\;\sum_{\tau,\tau(z)\sim z}\big(f(\tau\sigma(z))-f(\sigma(z))-f(\tau(z))+f(z)\big)^{2}\\
&+\sum_{\sigma,\sigma(z)\sim z}\;\sum_{\tau,\tau(z)\sim z}(f(\tau(z))-f(z))(f(\tau\sigma(z))-f(\sigma(z))-f(\tau(z))+f(z))
\end{align*}
where we have used the identity $A^{2}-B^{2}=(A-B)^{2}+2B(A-B)$, and 
\begin{align*}
   -2\Gamma(f,L_{0} f)(z)&=-\sum_{\tau,\tau(z)\sim z}(f(\tau(z))-f(z))(L_{0} f(\tau(z))-L_{0} f(z))\\
   &=-\sum_{\sigma,\sigma(z)\sim z}\;\sum_{\tau,\tau(z)\sim z} (f(\tau(z))-f(z))(f(\sigma\tau(z))-f(\sigma(z))-f(\tau(z))+f(z)).
\end{align*}
From the commutativity  assumption  $\sigma\tau=\tau\sigma$, and  summing up one gets
\begin{equation*}
    2\Gamma_{2}(f)(z)=\frac{1}{2}\sum_{\sigma,\sigma(z)\sim z}\;\sum_{\tau,\tau(z)\sim z} (f(\tau\sigma(z))-f(\sigma(z))-f(\tau(z))+f(z))^{2}\geq 0.
\end{equation*}
\end{proof}

\section{Appendix B}
\subsection{Proofs of Theorem \ref{thmprinc}, Theorem \ref{thmprincbis} and  Theorem \ref{Thmstructure}
}
\begin{proof}[Proof of Theorem \ref{thmprinc}]
Theorem \ref{thmprinc} is a consequence of Lemma 3.1 and Theorem  3.5 of \cite{Sam21}. 
The results of this paper \cite{Sam21} are given for graph spaces and the two  following additional assumptions : 
the measure $m$ is uniformly upper bounded and lower bounded away from 0,
\begin{equation*}
\sup_{x\in \X} m(x)<\infty,\qquad \inf_{x\in \X} m(x)>0,
\end{equation*}
and the generator $L$ is uniformly upper bounded, and uniformly lower bounded away from zero on the set of neighbours,
\[\sup_{x\in \X} |L(x,x)|<\infty,\qquad \inf_{x,y\in \X, d(x,y)=1} L(x,y)>0.\]
These conditions are not  in the setting of Theorem \ref{thmprinc}. To overcome this difficulty, one will consider a well chosen space $(\Cc,d,L_\Cc, m_\Cc)$ defined as the restriction of the   
space $(\X,d,m,L)$ to a well chosen finite  convex subset $\Cc$ as defined in Section \ref{sectionrestconv}. 

Let $\nu_0$ and $\nu_1$ be two probability measures on $\X$ with   bounded support. Since each vertex has bounded degree, there exists a finite convex subset $\Cc$ of $\X$ that contains all the balls of radius 2 with center in the finite subset $[\supp(\nu_0),\supp(\nu_1)]$. Choose for example the convex subset $\Cc$ with minimal elements. Let $(\widehat \nu_t)_{t\in[0,1]}$ denotes the Schr\"odinger bridge at zero temperature selected from the slowing down procedure on the space $(\Cc,d,L_\Cc, m_\Cc)$. As explained in \cite{Sam21} there exists  a $W_1$-optimal coupling $\widehat \pi$ with marginals $\nu_0$ and $\nu_1$ such that the expression of $(\widehat \nu_t)_{t\in[0,1]}$ is given by \eqref{defhatnut} on the space $(\Cc,d,L_\Cc, m_\Cc)$. Due to the
assumption on the subset $\Cc$, for any $(x,y)\in \supp \widehat\pi$, the set $[x,y]$ is the same on the space $(\X,d,m,L)$ and on the space $(\Cc,d,L_\Cc, m_\Cc)$. Moreover, since $L_\Cc(x,y)=L(x,y)$ for $x\neq y$, the expression of $r(x,z,z,y)$ for $z\in[x,y]$ and $(x,y)\in \supp (\widehat\pi)$ is also the same on $(\X,d,m,L)$ and on $(\Cc,d,L_\Cc, m_\Cc)$.  Therefore the expression of the Schr\"odinger bridges between Dirac measure 
$\delta_x$ and $\delta_y$ is for us given by \eqref{pathdirac} does not depend on the chosen convex subset $\Cc$. Up to now we are working on $(\Cc,d,L_\Cc, m_\Cc)$ but for most of all expressions we write, there is no dependence in $\Cc$ or the subset $\Cc$ may be replaced by $\X$. 

As already mentioned, the subset $\supp \,\widehat\pi$ is $d$-cyclically monotone. According to the definitions introduced in Section \ref{sectmainresults}, for any $t\in (0,1)$ the support of  $\widehat  \nu_t$ is given by   
\[\supp(\widehat  \nu_t)=\widehat Z:=Z(\supp \,\widehat\pi)=\bigcup_{(x,y)\in \supp\,\widehat \pi}[x,y].\]
 and one denotes $C_{_\rightarrow}:=C_{_\rightarrow}(\supp \, \widehat \pi)$, $C_{_\leftarrow}:=C_{_\leftarrow}(\supp\, \widehat \pi)$, and for any $z\in \widehat Z$
\[V_{_\rightarrow}(z):=V_{_\rightarrow}^{\supp\,\widehat \pi}(z),\quad
V_{_\leftarrow}(z):=V_{_\leftarrow}^{\supp\,\widehat \pi}(z),\quad
\V_{_\rightarrow}(z):=\V_{_\rightarrow}^{\supp\,\widehat \pi}(z),\quad
\V_{_\leftarrow}(z):=\V_{_\leftarrow}^{\supp\,\widehat \pi}(z).\]

 For any $z\in \widehat Z$
let also
 \[ \widehat Y_z:=\Big\{y\in\supp(\nu_1)\,\Big|\, \exists x\in \X,  (x,y)\in \widehat{\pi}, z\in[x,y]\Big\},\]
 and identically let 
 \[ \widehat X_z:=\Big\{x\in\supp(\nu_0)\,\Big|\, \exists y\in \X,  (x,y)\in \widehat{\pi}, z\in[x,y]\Big\}.\]
 For $y\in \supp(\nu_1)$, $z\in \X$ and $t\in [0,1]$, the quantity \begin{equation*}
 a_t(z,y) :=\int\nu_t^{w,y}(z) \,d\widehat{\pi}_{_\leftarrow}(w|y),
 \end{equation*}
 is positive if and only if $z\in \widehat Z$ and 
 $y\in \widehat Y_z$. Identically, for  $x\in \supp(\nu_0)$,  $z\in \X$ and $t\in [0,1]$, the quantity 
 \[b_t(z,x):=\int \nu_t^{x,w}(z) \,d\widehat{\pi}_{_\rightarrow}(w|x),\]
  is positive if and only if $z\in \widehat Z$ and $x\in \widehat X_z$. Actually $a_t$ and $b_t$ represent conditional laws, $\sum_{z\in \X}a_t(z,y)=\sum_{z\in \X}b_t(z,x) =1$.

For $t\in[0,1]$, $z\in \widehat Z$,   $\tz\in S_1(z)$ and $y\in\supp(\nu_1)$, let 
\begin{equation}\label{a_t'}
{\mathrm a}_t(z,\tz,y):= \sum_{w\in \X, (z,\tz)\in[y,w]}  
r(y,z,\tz,w) \,d(y,w)\, \B_t^{d(y,w)-1}(d(z,w)-1)\,\widehat{\pi}_{_\leftarrow}(w|y),
\end{equation}
and for any $x\in \supp(\nu_0)$, let 
\begin{equation*}
{\mathrm b}_t(z,\tz,x):=\sum_{w\in \X, (z,\tz)\in[x,w]}  
\,r(x,z,\tz,w) \,d(x,w)\, \B_t^{d(x,w)-1}(d(x,z))\,\widehat{\pi}_{_\rightarrow}(w|x).
\end{equation*}
where the function $r$ is given by \eqref{defrpont}.
For $t\in(0,1)$, the quantity ${\mathrm a}_t(z,\tz,y)$ is positive if and only if $\tz\in V_{_\leftarrow}(z)$ and 
$y\in \widehat Y_{(z,\tz)}$ with
 \[ \widehat Y_{(z,\tz)}=\Big\{y\in\supp(\nu_1)\,\Big|\, \exists x\in \X,  (x,y)\in \widehat{\pi}, (z,\tz)\in[y,x]\Big\}\subset \widehat Y_z\cap  \widehat Y_{\tz} ,\] 
 According to \cite[Lemma 3.4]{Sam21}, given $z\in \widehat Z$ and $\tz\in V_{_\leftarrow}(z)$  the ratio ${\mathrm a}_t(z,\tz,y)/a_t(z,y)$ does not depend on $y\in \widehat Y_{(z,\tz)}$. Therefore, for any $z\in \X$
 and $\tz\in S_1(z)$, one may define 
 \[A_t(z,\tz):=\left\{\begin{array}{ll}
  \frac{{\mathrm a}_t(z,\tz,y)}{a_t(z,y)}& \mbox{ for } (z,\tz) \in \widehat Z\times V_{_\leftarrow}(z) \mbox{ and } y\in \widehat Y_{(z,\tz)}\neq \emptyset, \\
  0 & \mbox{ otherwise.} 
 \end{array}\right.\] 
 Identically, for $t\in(0,1)$, the quantity ${\mathbbm b}_t(z,\ttz,x)$ is positive if and only if $\tz\in V_{_\rightarrow}(z)$ and $x\in \widehat X_{(z,\tz)}$ with
 \[ \widehat X_{(z,\tz)}=\Big\{y\in\supp(\nu_1)\,\Big|\, \exists x\in \X,  (x,y)\in \widehat{\pi}^0, (z,\tz)\in[x,y]\Big\}\subset \widehat X_z\cap  \widehat X_{\tz},\]
 and according to \cite[Lemma 3.4]{Sam21},  the ratio ${\mathrm b}_t(z,\tz,x)/b_t(z,x)$ does not depend on
$x\in X_{(z,\tz)}$. Therefore, for any $z\in \X$
 and $\tz\in S_1(z)$, one  defines 
 \[B_t(z,\tz):=\left\{\begin{array}{ll}
  \frac{{\mathrm b}_t(z,\tz,x)}{b_t(z,x)}& \mbox{ for } (z,\tz) \in \widehat Z\times V_{_\rightarrow}(z) \mbox{ and } x\in \widehat X_{(z,\tz)}\neq \emptyset, \\
  0 & \mbox{ otherwise.} 
 \end{array}\right.\] 
 Observe that by reversibility, for any $(z,\tz)\in C_{_{\rightarrow}}$ with $d(z,\tz)=1$,
 \begin{align*}
   B_t(z,\tz)&L(z,\tz)\widehat{\nu}_t(z) =\sum_{x\in \widehat X_{(z,\tz)}} B_t(z,\tz)L(z,\tz) b_t(z,x)\nu_0(x)
   =\sum_{x\in \widehat X_{(z,\tz)}} b_t(z,\tz,x)L(z,\tz) \nu_0(x)\\
   &=\sum_{(x,y)\in \supp(\widehat \pi), (z,\tz)\in[x,y]} r(x,z,\tz, y)L(z,\tz) d(x,y) \B_t^{d(x,y)-1}(d(x,z))\,\widehat{\pi}(x,y)\\
   &=\sum_{(x,y)\in \supp(\widehat \pi), (\tz,z)\in[y,x]} r(y,\tz,z, x)L(\tz,z) d(y,x) \B_t^{d(y,x)-1}(d(x,\tz)-1)\,\widehat{\pi}(x,y)\\
   &=\sum_{y\in \widehat Y_{(\tz,z)}} a_t(\tz,z,y)L(\tz,z) \nu_1(y)
   =A_t(\tz,z)L(\tz,z)\widehat{\nu}_t(\tz).
 \end{align*}

For $t\in[0,1]$,  $z\in \widehat Z$, $\ttz\in S_2(z)$  and  $y\in\supp(\nu_1)$, define also 
 \begin{equation}\label{a_t''}
{\mathbbm a}_t(z,\ttz,y):= \!\!\!\!\!\sum_{w\in \X, (z,\ttz)\in[y,w]}  
\!\!\!\!\!r(y,z,\ttz,w) \,d(y,w)(d(y,w)-1)\, \B_t^{d(y,w)-2}(d(z,w)-2)\,\widehat{\pi}_{_\leftarrow}(w|y),
\end{equation}
and for $x\in\supp(\nu_0)$
\begin{eqnarray*}
 {\mathbbm b}_t(z,\ttz,x):=\!\!\!\!\!\sum_{w\in \X, (z,\ttz)\in[x,w]} \!\!\!\!\! 
r(x,z,\ttz,w) \,d(x,w)(d(x,w)-1)\, \B_t^{d(x,w)-2}(d(x,z))
\,\widehat{\pi}_{_\rightarrow}(w|x).
\end{eqnarray*}
For $t\in(0,1)$, we also have ${\mathbbm a}_t(z,\ttz,y)>0$ if and only if $\ttz\in \V_{_\leftarrow}(z)$ and 
$y\in \widehat Y_{(z,\ttz)}$, 
and ${\mathbbm b}_t(z,\ttz,x)>0$ if and only if $\ttz\in \V_{_\rightarrow}(z)$ and 
$x\in  \widehat X_{(z,\ttz)}$. Since according to \cite[Lemma 3.4]{Sam21}, the ratio ${\mathbbm a}_t(z,\ttz,y)/a_t(z,y)$ does not depend on $y\in \widehat Y_{(z,\ttz)}$, and the ratio ${\mathbbm b}_t(z,\ttz,x)/b_t(z,x)$ does not depend on
$x\in X_{(z,\tz)}$. Therefore one may define for any $z\in \X$ and $\ttz\in S_2(z)$, 
\[{\mathbb A}_t(z,\ttz):=\left\{\begin{array}{ll}
  \frac{{\mathbbm a}_t(z,\ttz,y)}{a_t(z,y)}& \mbox{ for } (z,\ttz) \in \widehat Z\times \V_{_\leftarrow}(z) \mbox{ and } y\in \widehat Y_{(z,\ttz)}\neq \emptyset, \\
  0 & \mbox{ otherwise,} 
 \end{array}\right.\] 
and
\[{\mathbb B}_t(z,\ttz):=\left\{\begin{array}{ll}
  \frac{{\mathbbm b}_t(z,\ttz,x)}{b_t(z,x)}& \mbox{ for } (z,\ttz) \in \widehat Z\times \V_{_\rightarrow}(z) \mbox{ and } x\in \widehat X_{(z,\ttz)}\neq \emptyset, \\
  0 & \mbox{ otherwise.} 
 \end{array}\right.\] 
One also observe that  for $t\in(0,1)$ and $z\in\widehat Z$, if $\ttz\in \V_{_\leftarrow}(z)$ (or equivalently  ${\mathbb A}_t(z,\ttz)>0$), then   $A_t(z,\tz)>0$ for any $\tz\in S_1(z)$ with  $\tz\sim \ttz$ (since $\tz\in V_{_\leftarrow}(z))$. Therefore for any $t\in(0,1)$, $z\in \widehat Z$, $\tz\in S_1(z),\ttz \in S_2(z)$, one has $(A_t(z,\tz),{\mathbb A}_t(z,\ttz))\in (0,+\infty)\times [0,+\infty)\cup\{(0,0)\}$.   Identically, one has $(B_t(z,\tz),{\mathbb B}_t(z,\ttz))\in (0,+\infty)\times [0,+\infty)\cup\{(0,0)\}$. 
 
 As above, one simply check that by reversibility, for any $(z,\ttz)\in C_{_{\rightarrow}}$ with $d(z,\ttz)=2$,
 \begin{align}\label{RelAB}
   {\mathbb B}_t(z,\ttz)L^2(z,\ttz)\widehat{\nu}_t(z)={\mathbb A}_t(\ttz,z)L^2(\ttz,z)\widehat{\nu}_t(\ttz). 
 \end{align}
 
One will apply the following theorem which is a direct result of Lemma 3.1 and the main Theorem 3.5 of \cite{Sam21}. For $z\in \widehat Z$ and $t\in (0,1)$, let 
\begin{multline*}
H_t(z):=\Big(\sum_{\tz\in V_{_\leftarrow}(z)}
A_t(z,\tz) \,L(z,\tz)\Big)^2\\  
+ \sum_{\tz\in V_{_\leftarrow}(z),\, \ttz\in \V_{_\leftarrow}(z), \,\tz\sim \ttz} \rho \Big(A_t^2(z,\tz),{\mathbb A_t}(z,\ttz)\Big)  \, L(\tz, \ttz)L(z,\tz),
\end{multline*}
and let
\begin{multline*}
K_t(z):=\Big(\sum_{\tz\in V_{_\rightarrow}(z)}
B_t(z,\tz) \,L(z,\tz)\Big)^2 \\+  \sum_{\tz\in V_{_\rightarrow}(z),\, \ttz\in \V_{_\rightarrow}(z),\, \tz\sim \ttz}  \rho\Big(B_t^2(z,\tz),\mathbbm{B}_t(z,\ttz)\Big)  \,
 L(\tz,\ttz)L(z,\tz),
 \end{multline*}
where the function $\rho:(0,+\infty)\times[0,+\infty)\cup\{(0,0)\}\to \R$ is defined by 
\begin{equation*}
\rho(a,b):=\left(\log b-\log a-1\right) b,  \qquad a>0, b>0,
\end{equation*}
and $\rho(a,0)=0$ for $a\geq 0$.
According to \cite[Lemma 3.1]{Sam21} and \cite[Theorem 3.5]{Sam21} the following result holds.
\begin{theorem}\label{thmsam21}  We assume that the discrete space $(\X,d,m,L)$ is a graph space. 
Let $(\widehat \nu_t)_{t\in [0,1]}$ be the Schr\"odinger bridge at zero temperature
between two  probability measures $\nu_0,\nu_1\in \Pc(\X)$ with bounded support defined above and given by \eqref{defhatnut}. For any $t\in (0,1)$, let $q_t$ be the kernel on $[0,1]$ defined by \[q_t(s)=\frac{2s}t \1_{[0,t]}(s)+ \frac{2(1-s)}{1-t} \1_{[t,1]}(s),\qquad s\in[0,1].\] 
Then, one has 
\[(1-t)H(\nu_0|m)+tH(\nu_1|m)-H(\widehat \nu_t|m)\geq \int_0^1 \left(\int \big(H_s+K_s\big)  \,d \widehat\nu_s\right) q_t(s)\,ds.\]
As a consequence if there exists a real function $\zeta:[0,1]\to \R$ 
such that for any $s\in(0,1)$,
\begin{equation*}
\int \big(H_s+K_s\big)  \,d \widehat\nu_s \geq \zeta(s),
\end{equation*}
and if $\zeta q_t$ is integrable with respect to the Lebesgue measure on $[0,1]$, then the convexity property of entropy \eqref{deplacebis} holds with, for any $t\in (0,1)$,  
 \[
 C_t(\widehat \pi)= \int_0^1 \zeta(s) q_t(s)\, ds.\]
 \end{theorem}
 
 The proof of Theorem \ref{thmprinc} will therefore follows from an appropriate lower bound $\zeta(t)$ of $\int \big(H_t+K_t\big)  \,d \widehat\nu_t$ for any $t\in(0,1)$. Observe that  if  $\zeta$ is a constant function then the cost $C_t(\widehat \pi)$ is equal to this constant since $\int_0^1 q_t(s)\, ds=1$. And if $\zeta=\xi''$ where $\xi$ is a real continuous functions   on $[0,1]$, twice differentiable on $(0,1)$, then one has
 \[
 C_t(\widehat \pi)= \frac2{t(1-t)}\Big[
(1-t)\xi(0)+t\xi(1)-\xi(t)\Big].\]
 
 Let us  first rewrite the quantity $\int \big(H_t+K_t\big)  \,d \widehat\nu_t$.
 Using  the following  identity, for any integer $N$, for any $b\geq 0$,  and any positive $L_1,\ldots, L_N, a_1,\ldots, a_N$,
 \[\sum_{i=1}^N  \rho(a_i^2,b) L_i = L \,\rho\Big(\prod_{i=1}^N a_i^{2L_i/L},b\Big), \qquad \mbox{with} \quad L=\sum_{i=1}^N L_i,\]
 one gets for any $\ttz\in \widehat Z$,
  \begin{align}\label{epuise1}
H_t(\ttz )&=\Big(\sum_{\tz\in V_{_\leftarrow}(\ttz)}
A_t(\ttz,\tz) \,L(\ttz,\tz)\Big)^2+ \sum_{ z\in \V_{_\leftarrow}(\ttz)} \quad \sum_{\tz\in ]z,\ttz[ } \rho \Big(A_t^2(\ttz,\tz),{\mathbb A_t}(\ttz,z)\Big)  \, L(\tz, z)L(\ttz,\tz)\nonumber\\
&=\overline{ A}_t^2(\ttz)+ \sum_{ z\in \V_{_\leftarrow}(\ttz)} L^2(\ttz,z) \,\rho \Bigg(\prod_{\tz\in ]z,\ttz[ } A_t(\ttz,\tz)^{2\ell(\ttz,\tz,z)},{\mathbb A_t}(\ttz,z)\Bigg),
\end{align}
where $\ell(\ttz,\tz,z)= \frac{L(\ttz,\tz)L(\tz, z)}{ L^2(\ttz,z)}$, and 
\[\overline{ A}_t(\ttz):=\sum_{\tz\in V_{_\leftarrow}(\ttz)}
A_t(\ttz,\tz) \,L(\ttz,\tz).\]
Identically one gets for any $z\in\widehat Z$,
\begin{equation*}
K_t(z)= \overline{B}^2_t(z)+  \sum_{\ttz\in \V_{_\rightarrow}(z)} L^2(z,\ttz) \,\rho \Bigg(\prod_{\tz\in ]z,\ttz[ } B_t(z,\tz)^{2\ell(z,\tz,\ttz)},{\mathbb B_t}(z,\ttz)\Bigg),
\end{equation*}
with 
\[\overline{ B}_t(z):=\sum_{\tz\in V_{_\rightarrow}(z)}
B_t(z,\tz) \,L(z,\tz).\]
The reversibility property ensures that $\ell(z,\tz,\ttz)=\ell(\ttz,\tz,z)$. Setting \[\mathbb C_t(z,\ttz):=\sqrt{\mathbb A_t(\ttz,z)\mathbb B_t(z,\ttz)}\] and since $\mathcal L^2(z,\ttz)=\sqrt{L^2(z,\ttz)L^2(\ttz,z)}$, the above property of the function $\rho$ and the symmetric property \eqref{RelAB} imply that  
\begin{align}\label{epuise0}
    \int &\big(H_t+K_t\big)  \,d \widehat\nu_t = \sum_{\ttz \in\widehat Z}\overline{ A}_t^2(\ttz) \widehat\nu_t(\ttz) + \sum_{z \in\widehat Z}\overline{ B}_t^2(z) \widehat\nu_t(z)\nonumber\\
   &+\!\!\!\!\!\!\!\!\!\sum_{(z,\ttz)\in C_{_\rightarrow}, d(z,\ttz)=2}\!\!\!\!\!\!\!\!\!  2\, \rho \Big(\prod_{\tz\in ]z,\ttz[ } \big[A_t(\ttz,\tz)B_t(z,\tz)\big]^{\ell(z,\tz,\ttz)},{\mathbb C_t}(z,\ttz)\Big) \mathcal L^2(z,\ttz) \sqrt{\widehat\nu_t(z)\widehat\nu_t(\ttz)}.
\end{align}
Let 
\[\overline{{\mathbbm A}}_t(\ttz):=\sum_{z\in \V_{_\leftarrow}(\ttz)}
A_t(\ttz,z) \,L^2(\ttz,z)\quad \mbox{and}\quad \overline{{\mathbbm B}}_t(z):=\sum_{\ttz\in \V_{_\rightarrow}(z)}
B_t(z,\ttz) \,L^2(z,\ttz).\]
According to \eqref{a_t''} and \eqref{RelAB},  one has
\begin{align*}
&\sum_{(z,\ttz)\in C_{_\rightarrow}, d(z,\ttz)=2} {\mathbb C_t}(z,\ttz) \mathcal L^2(z,\ttz) \sqrt{\widehat\nu_t(z)\widehat\nu_t(\ttz)} = \int  \overline{{\mathbbm B}}_t  \,d \widehat\nu_t=\int  \overline{{\mathbbm A}}_t  \,d \widehat\nu_t \\
&=\int \sum_{z\in\widehat Z} \sum_{\ttz\in \V_{_\leftarrow}(z)}
\frac{{\mathbbm a}_t(z,\ttz,y)}{a_t(z,y)} \,L^2(z,\ttz) a_t(z,y)\,d \nu_1(y)\\
&= \iint\!\!\!  \sum_{(z,\ttz), (z,\ttz)\in [y,w]}  \!\!\!\!\!\!\!\!\!\!\!\!
r(y,z,\ttz,w)L^2(z,\ttz) \,d(y,w)(d(y,w)-1)\, \B_t^{d(y,w)-2}(d(z,w)-2)\,d\widehat{\pi}(w,y)\\
&= \iint \sum_{k=2}^{d(y,w)} \Big(\!\!\!\sum_{ (z,\ttz)\in[y,w], \ttz\in \V_{_\leftarrow}(z), d(z,w)=k} \!\!\!\!\!\!\!\!\!\!\!\!\!\!\!\!\!\!r(y,z,\ttz,w)  L^2(z,\ttz) \Big) \\
&\qquad\qquad\qquad\qquad\qquad\qquad\qquad\qquad\qquad\B_t^{d(y,w)-2}(k-2) \,d(y,w)(d(y,w)-1)\,d\widehat{\pi}(w,y)\\
&= \iint \sum_{k=2}^{d(y,w)} \B_t^{d(y,w)-2}(k-2) \,d(y,w)(d(y,w)-1)\,d\widehat{\pi}(w,y)=T_2(\widehat\pi).\\
\end{align*}
As a consequence, according to  equality \eqref{epuise0},  the convexity property of the function $\rho$, and the identity $\rho(\lambda a,\lambda b)=\lambda \rho(a,b)$, $a>0,b,\lambda\geq 0$, imply
\begin{align}\label{epuise}
 &   \int \big(H_t+K_t\big)  \,d \widehat\nu_t \geq  \sum_{\ttz \in\widehat Z}\overline{ A}_t^2(\ttz) \widehat\nu_t(\ttz) + \sum_{z \in\widehat Z}\overline{ B}_t^2(z) \widehat\nu_t(z)\\
   &+2\, \rho \Big(\!\!\!\!\!\!\sum_{(z,\ttz)\in C_{_\rightarrow}, d(z,\ttz)=2}\!\!\!\!\!\!\mathcal L^2(z,\ttz) \sqrt{\widehat\nu_t(z)\widehat\nu_t(\ttz)} \prod_{\tz\in ]z,\ttz[ } \big[A_t(\ttz,\tz)B_t(z,\tz)\big]^{\ell(z,\tz,\ttz)} ,T_2(\widehat\pi)\Big) .\nonumber
\end{align}
From the definition of constant $K(\supp\,\widehat\pi)$ and since the function $a\mapsto \rho(a,b)$ is decreasing on $(0,+\infty)$ for any $b\geq 0$, it follows that 
\begin{equation}\label{epuisebis}
    \int \big(H_t+K_t\big)  \,d \widehat\nu_t \geq A^2+B^2+2\rho\big(K(\supp\,\widehat\pi)AB ,T_2(\widehat\pi)\big),
\end{equation}
with 
$A^2:=\int \overline{ A}_t^2\,d\widehat\nu_t$ and $B^2=\int\overline{ B}_t^2\,d\widehat\nu_t$.
Applying then the inequality 
\begin{equation}\label{ineprat}
\rho(K a,b)= \rho(a,b)-b\log K\geq -a -b\log K,\quad K,a>0 ,b\geq 0,
\end{equation}
one gets
 \[\int \big(H_t+K_t\big)  \,d \widehat\nu_t \geq A^2+B^2-2AB-2\log\big( K(\supp(\widehat\pi))\big) \, T_2(\widehat\pi)\geq -2\log (K) \, T_2(\widehat\pi).\]
 The proof of Theorem \ref{thmprinc} then ends by applying Theorem \ref{thmsam21}.
 Observe that in the last inequalities the definition the constant $K$ should be first given on the space  $(\Cc,d,L_\Cc, m_\Cc)$. But due to the construction of  $\Cc$  the definition of the constant $K$ on $(\Cc,d,L_\Cc, m_\Cc)$ does not depend on $\Cc$, but  only depends on the geometric structure of the space $(\X,d,L, m)$ and the values of the jump rates $L(x,y)$ for $x\neq y$. 
 \end{proof}

\begin{proof}[Proof of Theorem \ref{thmprincbis}] Let us assume that $K\leq 1$. Applying 
the inequality $\rho(c,b)\geq -c,$
for $c>0$, $b\geq 0$, \eqref{epuisebis} provides  
\[\int \big(H_t+K_t\big)  \,d \widehat\nu_t  \geq A^2 +B^2-2KAB\geq \Big(1-\frac K\lambda\Big) A^2 + (1-K\lambda) B^2, \]
for all $\lambda >0$.
 By Cauchy-Schwarz inequality, one has 
\begin{align*}
A^2 &=\int \sum_{z\in\widehat Z} \Big(\sum_{\tz\in V_{_\leftarrow}(z)}
\frac{{\mathrm a}_t(z,\tz,y)}{a_t(z,y)} \,L(z,\tz)\Big)^2 a_t(z,y)\,d\nu_1(y)\\&\geq \int  \Big(\sum_{z\in\widehat Z} \sum_{\tz\in V_{_\leftarrow}(z)}
{\mathrm a}_t(z,\tz,y)\,L(z,\tz)\Big)^2\,d \nu_1(y).
 \end{align*}
 Moreover,  according to \eqref{a_t'}, easy computations give
 \begin{align*}
 &\sum_{z\in\widehat Z}\;\; \sum_{\tz\in V_{_\leftarrow}(z)}
{\mathrm a}_t(z,\tz,y)\,L(z,\tz)\\
&= \sum_{w\in \X}\;\;\sum_{ (z,\tz)\in[y,w], \tz\in V_{_\leftarrow}(z)}  
r(y,z,\tz,w) \,d(y,w)\, \B_t^{d(y,w)-1}(d(z,w)-1)\,\widehat{\pi}_{_\leftarrow}(w|y) L(z,\tz)\\
&= \sum_{w\in \X}\sum_{k=1}^{d(y,w)}\Big(\!\!\!\sum_{ (z,\tz)\in[y,w], \tz\in V_{_\leftarrow}(z), d(z,w)=k}\!\!\!\!\!\!\!\!\! \!\!\!\!\!\!\!\!\!r(y,z,\tz,w)  L(z,\tz) \Big) \B_t^{d(y,w)-1}(k-1) \,d(y,w)\,\widehat{\pi}_{_\leftarrow}(w|y)\\
&= \sum_{w\in \X}\sum_{k=1}^{d(y,w)} \B_t^{d(y,w)-1}(k-1) \,d(y,w)\,\widehat{\pi}_{_\leftarrow}(w|y)
= \sum_{w\in \X}  \,d(y,w)\,\widehat{\pi}_{_\leftarrow}(w|y),\end{align*}
and therefore $A^2\geq  \widetilde T_{_\leftarrow}(\widehat \pi)$. Identically, one gets $B^2\geq \widetilde T_{_\rightarrow}(\widehat \pi)$.
It follows that for any $\lambda \in [K,1/K]$
\[
\int \big(H_t+K_t\big)  \,d \widehat\nu_t 
\geq \Big(1-\frac K\lambda\Big) \widetilde T_{_\leftarrow}(\widehat \pi) + (1-K\lambda) \widetilde T_{_\rightarrow}(\widehat \pi).
\]
Choosing then either $\lambda=K$,$\lambda=1/K$ or $\lambda=1$ gives 
\[
\int \big(H_t+K_t\big)  \,d \widehat\nu_t 
\geq \max\Big((1-K^2)\widetilde T(\widehat \pi), (1-K) \big(\widetilde T_{_\leftarrow}(\widehat \pi)+\widetilde T_{_\leftarrow}(\widehat \pi)\big)\Big).\]
Then the result of the first item of Theorem \ref{thmprincbis} follows by applying Theorem \ref{thmsam21}.

Due to the above computations, we know that  $\int \overline{A}_t \,d \widehat\nu_t=\int \overline{B}_t \,d \widehat\nu_t= W_1(\nu_0,\nu_1)$. Let $\alpha:=\frac{\overline{A}_t \widehat\nu_t}{ W_1(\nu_0,\nu_1)}$ and $\beta:=\frac{\overline{B}_t \widehat\nu_t}{ W_1(\nu_0,\nu_1)}$. Applying the inequality 
$\rho(c,b)\geq -c,$
for $c>0 ,b\geq 0$, the inequality \eqref{epuise} provides  
\begin{align*}
&\int \big(H_t+K_t\big)  \,d \widehat\nu_t\\
&\geq  W_1^2(\nu_0,\nu_1) \bigg[ \sum_{\ttz\in \widehat Z} \frac{\alpha^2(\ttz)}{\widehat\nu_t(\ttz)} + \sum_{z\in \widehat Z} \frac{\beta^2(z)}{\widehat\nu_t(z)} \\
&-2 \sum_{(z,\ttz)\in C_{_\rightarrow}, d(z,\ttz)=2}\mathcal L^2(z,\ttz) \prod_{\tz\in ]z,\ttz[ } \bigg(\frac{\beta(z,\tz)}{L(z,\tz)\sqrt{\widehat\nu_t(z)}}\frac{\alpha(\ttz,\tz)}{L(\ttz,\tz)\sqrt{\widehat\nu_t(\ttz)}}\bigg)^{\ell(z,\tz,\ttz)}\bigg],
\end{align*}
where $\alpha(\ttz,\tz):=\frac{A_t(\ttz,\tz)L(\ttz,\tz)\widehat\nu_t(\ttz)}{W_1(\nu_0,\nu_1)}$ and 
$\beta(\ttz,\tz):=\frac{B_t(z,\tz)L(z,\tz)\widehat\nu_t(z)}{W_1(\nu_0,\nu_1)}$.
Since \[\sum_{z\in \widehat Z}\sum_{\tz\in V_{_\rightarrow}(z)} \beta(z,\tz)=1\quad\mbox{and}\quad\sum_{\ttz\in \widehat Z} \sum_{\tz\in V_{_\leftarrow}(\ttz)} \alpha(\ttz,\tz) =1,\] according to the definition of $R_1(\supp(\widehat \pi))$, it follows that 
\[\int \big(H_t+K_t\big)  \,d \widehat\nu_t\geq R_1(\supp(\widehat \pi)) \, W_1^2(\nu_0,\nu_1)\geq r_1\,W_1^2(\nu_0,\nu_1).\]
Applying Theorem \ref{thmsam21} ends  the proof of the main part of \ref{(ii)thmprincbis} in Theorem \ref{thmprincbis}. 

We now turn to the proof of inequality \eqref{R_1R_1K}. Let $S$ be a $d$-cyclically monotone subset, and let $\alpha:\X\times\X\to \R_+$, $\beta:\X\times\X\to \R_+$, $\nu:\X\to \R_+$ be functions satisfying  conditions \eqref{condalphabetanu}.
From the inequality $2uv\leq u^2 +v^2$, $u,v\in \R$, and setting $\alpha(\ttz):=\sum_{\tz\in V^S_{_\leftarrow}(\ttz)} \alpha(\ttz,\tz)$, $\beta(z):=\sum_{\tz\in V^S_{_\rightarrow}(z)} \beta(z,\tz)$, one has 
\begin{align*}
2 \sum_{(z,\ttz)\in C_{_\rightarrow}(S), d(z,\ttz)=2}&\mathcal L^2(z,\ttz) \prod_{\tz\in ]z,\ttz[ } \left(\frac{\beta(z,\tz)}{L(z,\tz)\sqrt{\nu(z)}}\frac{\alpha(\ttz,\tz)}{L(\ttz,\tz)\sqrt{\nu(\ttz)}}\right)^{\ell(z,\tz,\ttz)}\\
\leq &\sum_{z\in Z(S)} \frac{\beta^2(z)}{\nu(z)}\sum_{\ttz\in \V^S_{_\rightarrow}(z)} L^2(z,\ttz)  \prod_{\tz\in ]z,\ttz[ }\left(\frac{\beta(z,\tz)}{\beta(z)}\right)^{2\ell(z,\tz,\ttz)}\\
& +\sum_{\ttz\in Z(S)} \frac{\alpha^2(\ttz)}{\nu(\ttz)}\sum_{z\in \V^S_{_\leftarrow}(\ttz)} L^2(\ttz,z)  \prod_{\tz\in ]z,\ttz[ }\left(\frac{\alpha(\ttz,\tz)}{\alpha(\ttz)}\right)^{2\ell(z,\tz,\ttz)}\\
&\leq \sum_{z\in Z(S)} \frac{\beta^2(z)}{\nu(z)} K(z,\V^S_{_\rightarrow}(z)) +\sum_{\ttz\in Z(S)} \frac{\alpha^2(\ttz)}{\nu(\ttz)} K(\ttz,\V^S_{_\leftarrow}(\ttz)),
\end{align*}
where the last inequality follows from the definition of $K(z,\V^S_{_\rightarrow}(z))$ and $K(\ttz,\V^S_{_\leftarrow}(\ttz))$. From this upper bound and using then using Cauchy-Schwarz inequality, it follows that
\begin{align*}
    &K(S)\geq \inf_{\alpha, \beta,\nu} \Big\{ \sum_{z\in Z(S)} \frac{\beta^2(z)}{\nu(z)} \big(1-K(z,\V^S_{_\rightarrow}(z))\big) +\sum_{z\in Z(S)} \frac{\alpha^2(z)}{\nu(z)} \big(1-K(z,\V^S_{_\leftarrow}(z))\big)\Big\}\\
    &\geq \inf_{\nu} \Big\{\frac{1}{ \sum_{z\in Z(S)} \big[1-K(z,\V^S_{_\leftarrow}(z))\big]^{-1} \1_{V^S_{_\leftarrow}(z)\neq \emptyset}\, \nu(z)}\\
    &\qquad\qquad\qquad\qquad\qquad+\frac{1} {\sum_{z\in Z(S)}  \big[1-K(z, \V^S_{_\rightarrow}(z))\big]^{-1}  \1_{V^S_{_\rightarrow}(z)\neq \emptyset}\, \nu(z)}\Big\},
\end{align*}
since the infimum now runs over all $\alpha:\X\to \R_+$, $\beta:\X\to \R_+$ such that $\sum_{z\in Z(S)} \alpha(z)=1$ and $\sum_{z\in Z(S)} \beta(z)=1$. Finally the inequality \eqref{R_1R_1K} follows by using  the identity 
\[\inf_{u,v>0, u+v\leq w^{-1}}\left\{ \frac{1}u+\frac{1}v\right\}= 4w,\quad w> 0.\]

In order to prove the last part of Theorem \ref{thmprincbis}, one extends  to any graphs ideas from the proof of \cite[Theorem 2.5]{Sam21} on the discrete hypercube. 
Coming back to \eqref{epuise1}, the convexity property of the function $\rho$ gives
\[H_t(\ttz )\geq \overline{ A}_t^2(\ttz)+  \,\rho \Bigg(\sum_{ z\in \V_{_\leftarrow}(\ttz)} L^2(\ttz,z)\prod_{\tz\in ]z,\ttz[ } A_t(\ttz,\tz)^{2\ell(\ttz,\tz,z)},\overline{{\mathbbm A}}_t(\ttz)\Bigg).
\]
The inequality \eqref{ineprat} together with the definition of the quantity $K\big(\ttz,\V_{_\leftarrow}(\ttz)\big)$ then provides  
\[H_t(\ttz )\geq -\log K\big(\ttz,\V_{_\leftarrow}(\ttz)\big)\,\overline{\mathbb A}_t(\ttz),
\]
and one may identically shows that 
$K_t(z )\geq -\log K\big(z,\V_{_\rightarrow}(z)\big)\,\overline{\mathbb B}_t(z)$. As a consequence, one gets
\begin{align*}
   \int \big(H_t+K_t\big)  &\,d \widehat\nu_t  \geq \int -\log K\big(z,\V_{_\leftarrow}(z)\big)\,\overline{\mathbb A}_t(z)  -\log K\big(z,\V_{_\rightarrow}(z)\big)\, \overline{\mathbb B}_t(z) \,d \widehat\nu_t(z)\\
   &=\int \sum_{z\in\widehat Z} \sum_{\ttz\in \V_{_\leftarrow}(z)}
-\log K\big(z,\V_{_\leftarrow}(z)\big)\, \1_{\V_{_\leftarrow}(z)\neq \emptyset}\, {\mathbbm a}_t(z,\ttz,y)\,L^2(z,\ttz) \,d \nu_1(y)
\\
& \quad+\int \sum_{z\in\widehat Z} \sum_{\ttz\in \V_{_\rightarrow}(z)}
-\log K\big(z,\V_{_\rightarrow}(z)\big) \,\1_{\V_{_\rightarrow}(z)\neq \emptyset}\, {\mathbbm b}_t(z,\ttz,x)\,L^2(z,\ttz) \,d \nu_0(x)\\
&= \iint C_t(x,y) \,d\widehat \pi(x,y),
\end{align*}
with, setting $d(x,y)=d$, \begin{align*}
&C_t(x,y):=\sum_{z\in[x,y]} -\log K\big(z,\V_{_\leftarrow}(z)\big)\,\1_{\V_{_\leftarrow}(z)\neq \emptyset} \,r(x,z,z,y) \,d(d-1)\, \B_t^{d-2}(d(x,z)-2) \\
&\qquad\qquad+\sum_{z\in[x,y]}  -\log K\big(z,\V_{_\rightarrow}(z)\big)\,\1_{\V_{_\rightarrow}(z)\neq \emptyset} \,r(x,z,z,y) \,d(d-1)\, \B_t^{d-2}(d(x,z))\\
&=\sum_{k=0}^{d}
\sum_{z\in [x,y], d(x,z)=k} -\log K\big(z,\V_{_\leftarrow}(z)\big) \1_{\V_{_\leftarrow}(z)\neq \emptyset} \,r(x,z,z,y)  \,\frac{k(k-1)}{t^2}\,{\rho_t^d}(k)\\
&\quad+\sum_{k=0}^{d}
\sum_{z\in [x,y], d(x,z)=k}-\log K\big(z,\V_{_\rightarrow}(z)\big)\1_{\V_{_\rightarrow}(z)\neq \emptyset} \,r(x,z,z,y)  \,  \frac{(d-k)(d-k-1)}{ (1-t)^2}\,  {\rho_t^d}(k).
\end{align*} 
A   lower bound on $C_t(x,y)$ as a function of $d=d(x,y)$ can be obtained as follows.
\begin{align*}
&C_t(x,y)\\
&\geq -\log \big(\sup_{z\in \X}K(z,S_2(z))\big)\Big(\frac{d(d-1)}{t^2}\,\rho_t^d(d)+\frac{(d-1)(d-2)}{t^2}\,\rho_t^d(d-1)\\&\qquad\qquad\qquad\qquad\qquad\qquad\qquad\qquad\qquad+\frac{d(d-1)}{(1-t)^2}\,\rho_t^d(0)+\frac{(d-1)(d-2)}{(1-t)^2}\,\rho_t^d(1)\Big)\\
&+\1_{d\geq 4} \sum_{k=0}^{d} \sum_{z\in [x,y], d(x,z)=k} 
\!\!\!\!\!\!\Big(\!\!-\log K\big(z,\V_{_\leftarrow}(z)\big) \frac{k(k-1)}{t^2}  -\log K\big(z,\V_{_\rightarrow}(z)\big) \frac{(d-k)(d-k-1)}{ (1-t)^2}\Big)\\
&\qquad\qquad\qquad\qquad\qquad\qquad\qquad\qquad\qquad\qquad\qquad\qquad\qquad\qquad\qquad\, r(x,z,z,y) \,{\rho_t^d}(k)\\
&\geq \frac{\overline{r}}4\, d(d-1)\Big(t^{d-2}+(1-t)^{d-2}+(d-2) t^{d-3}(1-t)+(d-2)(1-t)^{d-3}t\Big)\\
&\qquad\qquad\qquad+ \frac{\overline{r}}4 \,\1_{d\geq 4} \sum_{k=2}^{d-2} 
\left( \frac{\sqrt{k(k-1)}}{t} + \frac{\sqrt{(d-k)(d-k-1)}}{ 1-t}\right)^2{\rho_t^d}(k)
 \end{align*}
 where  for the last inequalities we use the fact that $\sum_{z\in[x,y],d(x,z)=k} r(x,z,z,y)=1$ and the inequality
 \[aA+bA=(a^{-1}+b^{-1})^{-1}\Big(\frac A\alpha+\frac B\beta\Big)\geq (a^{-1}+b^{-1})^{-1}\Big(\sqrt A+\sqrt B\Big)^2,\]
 with $a=-\log K\big(z,\V_{_\leftarrow}(z)\big)>0$, $b=-\log K\big(z,\V_{_\rightarrow}(z)\big)>0$, $A=\frac{k(k-1)}{t^2}$, $B=\frac{(d-k)(d-k-1)}{ (1-t)^2}$, $\alpha=\frac{a^{-1}}{a^{-1}+b^{-1}}$, $\beta=\frac{b^{-1}}{a^{-1}+b^{-1}}$, so that $\alpha+\beta=1$ and according to the definition of the constant $\overline{r}$, $4(a^{-1}+b^{-1})^{-1}\geq \overline{r}$
 (since for $z\in [x,y],2\leq d(x,z)\leq d-2$, $\V_{_\leftarrow}(z)\neq \emptyset$ and $\V_{_\rightarrow}(z)\neq \emptyset$ and $\V_{_\leftarrow}(z)\times \V_{_\rightarrow}(z)$ is a $d$-cyclically monotone subset of $\X\times \X$).
 
 Observing that  
 \begin{align*}\sum_{k=2}^{d-2}k(k-1){\rho_t^d}(k)&=\sum_{k=0}^{d}k(k-1){\rho_t^d}(k)-d(d-1)\big[t^d+(d-2)t^{d-1}(1-t)\big]\\&=d(d-1)\big[t^2-t^d-(d-2)t^{d-1}(1-t)\big] ,
\end{align*}
and according to the definition \eqref{defvt} of $u_t(d)$, one gets 
\begin{align*}
&d(d-1)\Big(t^{d-2}+(1-t)^{d-2}+(d-2) t^{d-3}(1-t)+(d-2)(1-t)^{d-3}t\Big)\\
&\qquad\qquad+  \1_{d\geq 4}
\sum_{k=2}^{d-2} 
\left( \frac{\sqrt{k(k-1)}}{t} + \frac{\sqrt{(d-k)(d-k-1)}}{ 1-t}\right)^2{\rho_t^d}(k)=4u_t(d),
\end{align*}
and therefore 
$C_t(x,y)\geq \overline{r} u_t(d)$. Then Theorem \ref{thmsam21} ensures that the displacement convexity property holds with 
\begin{equation*}
C_t(\widehat \pi)=\overline{r}\iint \overline{c}_t\big(d(x,y)\big)\,d\widehat \pi(x,y),
\end{equation*}
with $\overline{c}_t(d):=\int_0^1 u_s(d)\, q_t(s)\,ds $.
\end{proof}

\begin{proof}[Proof of Theorem \ref{Thmstructure}] We start with the proof of inequality \eqref{rr_3}. Let $z\in \X$. According to \eqref{defKtilde} and \eqref{defr_3}
\begin{align}\label{vent}
  & \widetilde r_2(z)=\inf_{W\subset S_2(z)}\left\{1-\widetilde K_L(z,W)\right\}\nonumber\\
   &=\inf_{W\subset S_2(z)}
   \inf_\beta  \Biggl\{ \left(\sum_{\sigma\in \Sc_{]z,W[}} \sqrt{\beta(\sigma)}\right)^2 -\sum_{\ttz\in W} L^2(z,\ttz) \prod_{\sigma\in \Sc_{]z,\ttz[}}\left(\frac{\beta(\sigma)}{\big(L(z,\sigma(z))\big)^2}\right)^{\ell(z,\sigma(z),\ttz)}\nonumber\Biggr\}\\
   &=\inf_{W\subset S_2(z)}
   \inf  \Biggl\{ \left(\sum_{\sigma\in \Sc_{]z,W[}} \sqrt{\beta(\sigma)}\right)^2 \left[1-\sum_{\ttz\in W} L^2(z,\ttz) \prod_{\sigma\in \Sc_{]z,\ttz[}}\left(\frac{\alpha(\sigma)}{L(z,\sigma(z))}\right)^{2\ell(z,\sigma(z),\ttz)}\right]\nonumber\\
   &\quad\qquad\qquad\qquad\qquad\qquad\, \Bigg|\,{\beta}:\Sc_{]z,W[} \to  \mathbb{R}_{+},\sum_{\sigma\in \Sc_{]z,W[}} \beta(\sigma)=1,\alpha(\sigma):=\frac{\sqrt{\beta(\sigma)}}{\sum_{\sigma\in \Sc_{]z,W[}} \sqrt{\beta(\sigma)}}  \Biggr\}
\end{align}
According to the definition of $K(z, W)$, one gets 
\begin{align*}
    \widetilde r_2(z)&\geq \inf_{W\subset S_2(z)} \left[\inf  \Biggl\{ \left(\sum_{\sigma\in \Sc_{]z,W[}} \sqrt{\beta(\sigma)}\right)^2\,\Bigg|\,{\beta}:\Sc_{]z,W[} \to  \mathbb{R}_{+},\sum_{\sigma\in \Sc_{]z,W[}} \beta(\sigma)=1\Biggr\}\Big(1-K(z,W)\Big)\right]\\
    &=\inf_{W\subset S_2(z)}\Big(1-K(z,W)\Big) =1-K(z,S_2(z))
\end{align*}
If $r(z)\geq 0$, by the Cauchy Schwarz inequality, \eqref{vent} provides
\begin{align*}
   \widetilde r_2(z)&\leq \inf_{W\subset S_2(z)}
   \inf  \Biggl\{ \big| \Sc_{]z,W[}\big| \, \left[1-\sum_{\ttz\in W} L^2(z,\ttz) \prod_{\sigma\in \Sc_{]z,\ttz[}}\left(\frac{\alpha(\sigma)}{L(z,\sigma(z))}\right)^{\frac{2L(z,\sigma(z))L(\sigma(z),\ttz)}{L^2(z,\ttz)}}\right]\nonumber\\
   &\qquad\qquad\qquad\qquad\, \Bigg|\,{\beta}:\Sc_{]z,W[} \to  \mathbb{R}_{+},\sum_{\sigma\in S_{]z,W[}} \beta(\sigma)=1,\alpha(\sigma):=\frac{\sqrt{\beta(\sigma)}}{\sum_{\sigma)\in \Sc_{]z,W[}} \sqrt{\beta(\sigma)}}  \Biggr\}\\
   &\leq \big|S_1(z)\big|\inf_{W\subset S_2(z)}
   \inf  \Biggl\{ 1-\sum_{\ttz\in W} L^2(z,\ttz) \prod_{\sigma\in S_{]z,\ttz[}}\left(\frac{\alpha(\sigma)}{L(z,\sigma(z))}\right)^{\frac{2L(z,\sigma(z))L(\sigma(z),\ttz)}{L^2(z,\ttz)}}\nonumber\\
   &\qquad\qquad\qquad\qquad\, \Bigg|\,{\beta}:\Sc_{]z,W[} \to  \mathbb{R}_{+},\sum_{\sigma\in \Sc_{]z,W[}} \beta(\sigma)=1,\alpha(\sigma):=\frac{\sqrt{\beta(\sigma)}}{\sum_{\sigma\in \Sc_{]z,W[}} \sqrt{\beta(\sigma)}}  \Biggr\}\\
   &= \big|S_1(z)\big| \inf_{W\subset S_2(z)}\Big(1-K(z,W)\Big) =\big|S_1(z)|\big(1-K(z,S_2(z))\big).
\end{align*}
This ends the proof of inequality \eqref{rr_3}.

The proof of the lower bound $\widetilde r_2$ of the $\widetilde T_2$-entropic curvature of the space is similar to the one of Theorem \ref{thmprinc} or   Theorem \ref{thmprincbis}.  
Starting again from inequality \eqref{epuise0} and setting   
\[\widetilde A_t^2(z):=\sum_{\sigma\in \Sc,\sigma(z)\in V_{_\leftarrow}(z)}\Big(
A_t(z,\sigma(z)) \,L(z,\sigma(z))\Big)^2,\]
one gets 
\begin{align*}
H_t(z)&\nonumber\geq \widetilde A_t^2(z)+ \sum_{\sigma,\tau\in \Sc,\sigma\neq \tau, \sigma(z),\tau(z)\in V_{_\leftarrow}(z)} A_t(z,\sigma(z)) \,L(z,\sigma(z)) A_t(z,\tau(z)) \,L(z,\tau(z))\\
&\nonumber\qquad+\rho \bigg(\sum_{ \ttz\in \V_{_\leftarrow}(z)} L^2(z,\ttz)\prod_{\tz\in [z,\ttz] \cap S_1(z)} A_t(z,\tz)^{\frac{2L(z,\tz)L(\tz, \ttz)}{ L^2(z,\ttz)}},\overline{\mathbb A}_t(z)\bigg)\\
&= \widetilde A_t^2(z)\Big[1+\sum_{\sigma,\tau\in \Sc_{]z,\V_{_\leftarrow}(z)[},\sigma\neq \tau} \sqrt{\beta(\sigma,z)}\, \sqrt{\beta(\tau,z)}\Big]\\
\nonumber&\qquad+ \rho \bigg(\widetilde A_t^2(z) \sum_{ \ttz\in \V_{_\leftarrow}(z)} L^2(z,\ttz)\prod_{\sigma \in \Sc_{]z,\ttz[}} \bigg(\frac{\beta(\sigma,z)}{\big(L(z,\sigma(z))\big)^2}\bigg)^{\frac{L(z,\sigma(z))L(\sigma(z), \ttz)}{ L^2(z,\ttz)}},\overline{\mathbb A}_t(z)\bigg)
\end{align*}
where for any $z\in\X$ and  $\sigma\in  \Sc_{]z,\V_{_\leftarrow}(z)[}$,  $\beta(\sigma,z):=\frac{\left(
A_t(z,\sigma(z)) \,L(z,\sigma(z))\right)^2}{\widetilde A_t^2(z)}.$
Using the inequality $\rho(a,b)\geq -a$, it follows that 
\begin{align*}
H_t(z)&\geq \widetilde A_t^2(z)
\Biggl[1+  \sum_{\sigma,\tau\in \Sc_{]z,\V_{_\leftarrow}(z)[},\sigma\neq \tau} \sqrt{\beta(\sigma,z)}\, \sqrt{\beta(\tau,z)} \\ 
&\qquad\qquad - \sum_{ \ttz\in \V_{_\leftarrow}(z)} L^2(z,\ttz)\prod_{\sigma \in \Sc_{]z,\ttz[}} \bigg(\frac{\beta(\sigma,z)}{\big(L(z,\sigma(z))\big)^2}\bigg)^{\frac{L(z,\sigma(z))L(\sigma(z), \ttz)}{ L^2(z,\ttz)}}\Biggr]\nonumber,
\end{align*}
and therefore,
from the definition of $\widetilde K(z,\V_{_\leftarrow}(z))$, 
\begin{equation*}
H_t(z)\geq (1-\widetilde K(z,\V_{_\leftarrow}(z)) \widetilde A_t^2(z).
\end{equation*}
According to the definition of constant $\widetilde r_2$, one has 
\begin{align}\label{coudeville}
\int H_t \, d\widehat \nu_t &\geq \widetilde r_2 \int \widetilde A_t^2(z) \, d\widehat \nu_t(z)\nonumber\\
&=\widetilde r_2\sum_{\sigma\in \Sc}\int \sum_{z\in \X}\frac{\big(a_t(z,\sigma(z),y)L(z,\sigma(z))\big)^2}{a_t(z,y)}\,\1_{z\in [y,\sigma(z)[} \,d\nu_1(y)\nonumber\\
&\geq\widetilde r_2\sum_{\sigma\in \Sc}\int \frac{\Big(\sum_{z\in \X,z\in [y,\sigma(z)[ } a_t(z,\sigma(z),y)L(z,\sigma(z))\Big)^2}{\sum_{z\in \X, z\in [y,\sigma(z)[ } a_t(z,y)}\,d\nu_1(y)
\end{align}
where the last inequality holds if $\widetilde r_2\geq 0$ by applying Cauchy-Schwarz inequality.

According to the  definition of $a_t(z,\sigma(z),y)$ given by \eqref{a_t'}, one has 
\begin{align*}
    &\sum_{z\in \X,z\in [y,\sigma(z)[ } a_t(z,\sigma(z),y)L(z,\sigma(z))\\&=\sum_{w\in \X}\sum_{z\in\X}\1_{(z,\sigma(z))\in [y,w]}\,d(y,w) r(y,z,\sigma(z),w)L(z,\sigma(z))\,\rho_t^{d(y,w)-1}(d(z,w)-1)\,\widehat \pi_{\leftarrow}(w|y)\\
    &=\sum_{w\in \X}   \!\!\!\sum_{k=0}^{d(w,y)-1} \!\!\! d(y,w)\rho_t^{d(y,w)-1}(d(y,w)-1-k) \,\widehat \pi_{\leftarrow}(w|y) \!\!\!  \!\!\!\!\!\! \!\!\!  \!\!\! \!\!\!\!\!\!\!\!\!\!\sum_{z\in \X, d(y,z)=k, (z,\sigma(z))\in [y,w]} \!\!\! \!\!\! \!\!\!\!\!\!\!\!\!\! \!\!\! \!\!\! \!\!\! r(y,z,\sigma(z),w)L(z,\sigma(z))
\end{align*}
From Lemma \ref{lemL} \ref{3}, we know that if $(z,\sigma(z))\in [y,w]$ then $\sigma(y)\in]y,w]$,  and $\sigma(z)\in [\sigma(y),w]$, and   Lemma \ref{lemL} \ref{5} implies 
\begin{align*}
&\sum_{z\in \X, d(y,z)=k, (z,\sigma(z))\in [y,w]} r(y,z,\sigma(z),w)L(z,\sigma(z))\\
& =\frac{\sum_{z\in \X, d(y,z)=k, (z,\sigma(z))\in [y,w]} L^{d(y,z)}(y,z)L(z,\sigma(z))L^{d(\sigma(z),w)}(\sigma(z),w)}{L^{d(y,w)}(y,w)}\\
&=\frac{L(y,\sigma(y))L^{d(y,w)-1}(\sigma(y),w)}{L^{d(y,w)}(y,w)}=r(y,\sigma(y),\sigma(y),w)
\end{align*}
It follows that
\begin{align*}
    \sum_{z\in \X,z\in [y,\sigma(z)[ } a_t(z,\sigma(z),y)L(z,\sigma(z))&=\sum_{w\in \X}\1_{\sigma(y)\in]y,w]}d(y,w) r(y,\sigma(y),\sigma(y),w)\,\widehat\pi_{_\leftarrow}(w|y)\\&=\Pi^\sigma_\leftarrow(y).
\end{align*} 
Observing that 
\begin{equation}\label{obvious}
\sum_{z\in \X, z\in [y,\sigma(z)[ } a_t(z,y)\leq 1,
\end{equation}
inequality \eqref{coudeville} therefore provides
\[\int H_t \, d\widehat \nu_t \geq \widetilde r_2\int\sum_{\sigma\in \Sc} \Pi^\sigma_\leftarrow(y)^2d\nu_1(y).\]
We similarly prove that 
$\int K_t \, d\widehat \nu_t \geq \widetilde r_2\int\sum_{\sigma\in \Sc} \Pi^\sigma_\rightarrow(x)^2d\nu_0(x)$,
and thus we get 
\[\int 
(H_t+K_t) \, d\widehat \nu_t \geq \widetilde r_2\, \widetilde{T}_2(\widehat \pi).\]
The proof of the first part of Theorem \ref{Thmstructure}  ends by applying Theorem \ref{thmsam21}.

Let now assume that condition \eqref{restr} also holds, the second part of Theorem \ref{Thmstructure} will follows from improving the trivial bound \eqref{obvious}.
One has 
\begin{align}\label{long}
   &\sum_{z\in \X, z\in [y,\sigma(z)[ } a_t(z,y)=1-\sum_{z\in \X, z\not\in [y,\sigma(z)[ } a_t(z,y)\nonumber\\
   &\leq 1-\sum_{x\in\X}\1_{\sigma(y)\in]y,x]}\sum_{z\in [y,x], z\not\in [y,\sigma(z)[ } r(x,z,z,y)\rho_t^{d(x,y)}(d(x,z)) \,\widehat \pi_{_\leftarrow}(x|y)\\
   \nonumber &=1-\sum_{x\in\X}\frac{\1_{\sigma(y)\in]y,x]}}{L^{d(x,y)}(x,y)} \sum_{k=0}^{d(x,y)}\rho_t^{d(x,y)}(k)\left(\sum_{z\in Z(y,x,k,\sigma) }\quad
   \sum_{\gamma\in G(y,x), z\in \gamma} L(\gamma)\right) \widehat \pi_{_\leftarrow}(x|y),
\end{align}
where  $Z(y,x,k,\sigma):=\big\{z\in [y,x]\,\big|\,d(y,z)=d(x,y)-k, z\not\in [y,\sigma(z)[\big\}$.
 Given $\sigma\in \Sc$ the set of geodesics from $y$ to $x$ contains the set of geodesics using the move $\sigma$, more precisely, setting $d(x,y)=d$
\[G(y,x)\supset \bigcup_{\ell=0}^{d-1}G_{\sigma,\ell}(y,x),\quad \mbox{with}\quad G_{\sigma,\ell}(y,x):=\big\{\gamma=(z_0,\ldots,z_d)\in G(y,x)\,\big|\,z_{\ell+1}=\sigma(z_l)\big\}.\]
Observe that according to Lemma \ref{lemL} \ref{5}, if $\sigma(y)\in ]y,x]$ then 
for any $\ell\in\{0,\ldots, d-1\}$
\[\sum_{\gamma\in G_{\sigma,\ell}(y,x)}L(\gamma)=L(y,\sigma(y))L^{d(\sigma(y),x)}(\sigma(y),x).\]
According to assumption \eqref{restr}, for $\ell\neq \ell'$, $G_{\sigma,\ell}(y,x)$ and $G_{\sigma,\ell'}(y,x)$ are disjoints sets, and therefore
\begin{align*}
\sum_{z\in Z(y,x,k,\sigma) }\;
   \sum_{\gamma\in G(y,x), z\in \gamma} L(\gamma)
   &\geq \sum_{\ell=0}^{d-1}\quad\sum_{z\in [y,x],d(y,z)=d(x,y)-k, z\not\in [y,\sigma(z)[ }\quad
   \sum_{\gamma\in G_{\sigma,\ell}(y,x), z\in \gamma}  L(\gamma) \\
   &=\sum_{\ell=0}^{d-1} \quad
   \sum_{\gamma=(z_0,\ldots,z_d)\in G_{\sigma,\ell}(y,x)} \1_{z_{d-k}\not\in [y,\sigma(z_{d-k})[ } L(\gamma) 
\end{align*}
Assume that $\sigma(y)\in]y,x]$ and let $\gamma=(z_0,\ldots,z_d)\in G_{\sigma,\ell}(y,x)$ with $z_{d-k}\not\in [y,\sigma(z_{d-k})[$. Observe first that $k\neq d$ (otherwise $y=z_0\not\in [y,\sigma(z_0)[=[y,\sigma(y)[$ which is impossible). If $0\leq \ell<d-k$ then $\gamma=(z_0,\ldots,z_\ell,z_{\ell+1},\ldots,z_{d-k},\ldots,z_d)$ with  $z_{\ell+1}=\sigma(z_\ell)$. It follows that necessarily $z_{d-k}\not\in [y,\sigma(z_{d-k})[$, since otherwise $(z_0,\ldots,z_\ell,z_{\ell+1},z_{d-k},\sigma(z_{d-k}))$ is a geodesic from $z_0$ to $\sigma(z_{d-k})$ that uses the move $\sigma$ twice. It follows that 
\[\sum_{\gamma=(z_0,\ldots,z_d)\in G_{\sigma,\ell}(y,x)}\!\!\!\!\!\!\!\!\!\!\!\! \1_{z_{d-k}\not\in [y,\sigma(z_{d-k})[ } L(\gamma)=\!\!\!\!\!\!\sum_{\gamma=(z_0,\ldots,z_d)\in G_{\sigma,\ell}(y,x)} \!\!\!\!\!\!\!\!\! L(\gamma) = L(y,\sigma(y))L^{d(\sigma(y),x)}(\sigma(y),x).\]
Assume now that  $ d-k\leq\ell\leq d-1$, then   $\gamma=(z_0,\ldots,z_{d-k},\ldots, z_\ell,z_{\ell+1},\ldots,z_d)$ with  $z_{\ell+1}=\sigma(z_\ell)$. According to Lemma \ref{lemL} \ref{1}, $(z_0,\ldots,z_{d-k},\sigma(z_{d-k}),\ldots, \sigma(z_\ell),z_{\ell+2},\ldots,z_d)$ is also a geodesic in $G(y,x)$ and therefore  $z_{d-k}\in [y,\sigma(z_{d-k})[$.
As a consequence \[\sum_{\gamma=(z_0,\ldots,z_d)\in G_{\sigma,\ell}(y,x)} \1_{z_{d-k}\not\in [y,\sigma(z_{d-k})[ } L(\gamma)=0.\]
Finally, if $\sigma(y)\in]y,x]$, one gets for any fixed $k\in\{0,\ldots,d-1\}$,
\begin{align*}
\sum_{z\in [y,x],d(y,z)=d(x,y)-k, z\not\in [y,\sigma(z)[ }\quad
   \sum_{\gamma\in G(y,x), z\in \gamma} L(\gamma) 
&\geq\sum_{d=0}^{d-k-1}L(y,\sigma(y))L^{d(\sigma(y),x)}(\sigma(y),x)\\
   &=(d-k)\, L(y,\sigma(y))L^{d(\sigma(y),x)}(\sigma(y),x).
   \end{align*}
Observe that if assumption \eqref{restr} is not fulfilled, using the fact that $G(y,x)\supset G_{\sigma,\ell}(y,x)$, one identically gets 
\[\sum_{z\in [y,x],d(y,z)=d(x,y)-k, z\not\in [y,\sigma(z)[ }\quad
   \sum_{\gamma\in G(y,x), z\in \gamma} L(\gamma)\geq L(y,\sigma(y))L^{d(\sigma(y),x)}(\sigma(y),x).\]
As a consequence, if assumption \eqref{restr} is satisfied,  \eqref{long} provides 
\begin{align*}
\sum_{z\in \X, z\in [y,\sigma(z)[ } a_t(z,y)&\leq 1-\sum_{x\in\X}{\1_{\sigma(y)\in]y,x]}}r(y,\sigma(y),\sigma(y),x) \sum_{k=0}^{d(x,y)-1}\rho_t^{d(x,y)}(k)\big(d(x,y)-k\big) \widehat \pi_{_\leftarrow}(x|y)\\
&\quad=1-(1-t)\Pi^\sigma_\leftarrow(y), 
\end{align*}
and if assumption \eqref{restr} is not fulfilled,  \eqref{long} implies 
\begin{align*}
\sum_{z\in \X, z\in [y,\sigma(z)[ } a_t(z,y)&\leq 1-\sum_{x\in\X}\1_{\sigma(y)\in]y,x]} r(y,\sigma(y),\sigma(y),x) \sum_{k=0}^{d(x,y)-1}\rho_t^{d(x,y)}(k)\, \widehat \pi_{_\leftarrow}(x|y)\\
&=1-\sum_{x\in\X}\1_{\sigma(y)\in]y,x]} r(y,\sigma(y),\sigma(y),x) \big(1-t^{d(x,y)}\big)\, \widehat \pi_{_\leftarrow}(x|y)\\
&\leq 1-(1-t) \sum_{x\in\X}\1_{\sigma(y)\in]y,x]} r(y,\sigma(y),\sigma(y),x) \, \widehat \pi_{_\leftarrow}(x|y)\\
&\leq1-{(1-t)}\frac{\Pi^\sigma_\leftarrow(y)}{\textnormal{Diam}(\mathcal{X})},
\end{align*}
Setting $D=1$ if assumption \eqref{restr} holds and $D=\textnormal{Diam}(\mathcal{X})$ otherwise, inequality \eqref{coudeville} then provides
\[\int H_t \, d\widehat \nu_t \geq \widetilde r_2\int\sum_{\sigma\in \Sc} \frac{\Pi^\sigma_\leftarrow(y)^2}{1-(1-t)\frac{\Pi^\sigma_\leftarrow(y)}D}d\nu_1(y)= \xi_\leftarrow''(t),\]
with
\[ \xi_\leftarrow(t):=\frac{\widetilde r_2 D^2}2 \int \sum_{\sigma\in \Sc} h\Big((1-t) \,\frac{\Pi^\sigma_\leftarrow(y)}D\Big) \,d\nu_1(y).\]
One identically  proves that 
\[\int K_t \, d\widehat \nu_t \geq  \xi_\rightarrow''(t),\]
with
\[ \xi_\rightarrow(t):=\frac{\widetilde r_2 D^2}2 \int \sum_{\sigma\in \Sc} h\Big(t \,\frac{\Pi^\sigma_\rightarrow(x)}D\Big) \,d\nu_0(x).\]
The proof of the second part of Theorem \ref{Thmstructure} ends by applying Theorem \ref{thmsam21}.
\end{proof}

\subsection{Proofs of Theorem \ref{BM}, Theorem \ref{thmstructure}, Theorem \ref{Logsob}, Theorem \ref{Logsobbis} and Theorem \ref{ttens}}

\begin{proof}[Proof of Theorem \ref{BM}]
Let $x,y$ be distinct vertices in $\mathcal{X}$. By definition of entropic curvature we have 
\begin{equation}\label{bmyers}
\HR(\widehat \nu_t|m)\leq (1-t) \HR(\nu_0|m)+t \HR(\nu_1|m)- \kappa\, \frac{t(1-t)}{2}\,T_2(\widehat \pi) \hspace{0.1cm} ,
\end{equation}
for any $\nu_{0},\nu_{1}\in \Pc(\X)$ .
Let $\nu_{0}=\delta_{x}$ and $\nu_{1}=\delta_{y}$ then \eqref{bmyers} becomes 
\begin{align*}
d(x,y)(d(x,y)-1)&\leq \frac{-2}{\kappa t(1-t)}\HR(\widehat \nu_t|m) +\frac2{\kappa t} \log \frac{1}{m(x)} +\frac2{\kappa(1-t)} \log \frac{1}{m(y)}\\
&\leq \frac{2}{\kappa t(1-t)}\left(-\HR(\widehat \nu_t|m) + \log\frac1{\inf_{x\in \X}m(x)}\right)
\end{align*}
Furthermore by Jensen inequality we have
\begin{equation*}
-\HR(\widehat \nu_t|m)\leq \log m(\supp(\widehat \nu_t))\leq \log \Big(|\supp(\widehat \nu_t)| \sup_{x\in\X} m(x)\Big)\hspace{0.1cm},
\end{equation*}
and since  $|supp(\widehat \nu_t)|\leq \Delta(G)^{d(x,y)}$ one finally gets 
\begin{equation*}
d(x,y)\leq \frac{2}{\kappa t(1-t)}\log \Big(\Delta(G)\, \frac{\sup_{x\in\X} m(x)}{\inf_{x\in \X}m(x)}\Big) +1 <\infty \hspace{0.1cm} .
\end{equation*}
Choosing $t=1/2$ and maximizing over $x$ and $y$ ends the proof of Theorem \ref{BM}.
\end{proof}

\begin{proof}[Proof of Theorem \ref{thmstructure}]
According to Theorem \ref{thmprinc}, in order to prove that a structured graph $(\X,d,m_0,L_0)$  with associated finite set of moves $\Sc$ has non negative entropic curvature, it suffices to show that for any $z\in \X$, $K_0(z,S_2(z))\leq 1$.
Let $z\in \X$ be a fixed vertex, and for any $\ttz\in S_2(z)$ let  
\[ U(\ttz):=\big\{(\tau,\sigma)\in \Sc \,\big|\,\tau(\sigma(z))=\ttz\}.\]
Each couple $(\tau,\sigma)$   can be associated to a single geodesic $(z,\sigma(z),\tau(\sigma(z)))$ from $z$ to $\ttz$. Obviously for  $w''\in S_{2}(z)$ with $w''\neq \ttz$, the sets $U(\ttz)$ and $U(w'')$ are disjoints.

According to the definition of structured graphs, if $(\tau,\sigma)\in U(\ttz)$ then $(\psi(\sigma),\tau)\in U(\ttz)$ where $\psi:\Sc_z^{\cdot\rightarrow \tau}\to \Sc_z^{\tau\rightarrow \cdot}$ is a fixed one to one map. As a consequence,  given $(\sigma_2,\sigma_1)\in U(\ttz)$, one may construct by induction a sequence  $(\sigma_{k+1},\sigma_k)$, $k\in \N^*$, of elements in $U(\ttz)$ defined by $\sigma_{k+1}=\psi_{k-1}(\sigma_k)$ for all $k\geq 2$ with $\psi_{k-1}:\Sc_z^{\cdot\rightarrow \sigma_{k-1}}\to \Sc_z^{\sigma_{k-1}\rightarrow \cdot}$. 
Let us define 
\[\overline{(\sigma_2,\sigma_1)}:=\big\{(\sigma_{k+1},\sigma_k)\,\big|\, k\geq 1\big\}.\]
Since $\Sc$ is finite, there exists $k\geq 2$ and $j\leq k$ such that $\sigma_{k+1}=\sigma_j$.
Let 
\[\ell:=\min\big\{k\geq 1\,\big|\,\exists j\in\{1,\ldots,k\}, \sigma_{k+1}=\sigma_j\big\}.\]
 The maps $\sigma_1,\sigma_2,\ldots, \sigma_\ell$ all differs. Let $j\in [\ell]$ such that $\sigma_j=\sigma_{\ell+1}$. Let us prove that $j=1$. If $j\geq 2$ then $\sigma_{\ell+1}\sigma_\ell(z)=\ttz=\sigma_j\sigma_{j-1}(z)=\sigma_{\ell+1}\sigma_{j-1}(z)$.  Lemma \ref{lemL} \ref{6} implies $\sigma_\ell=\sigma_{j-1}$ which contradicts the definition of $\ell$. 
It follows that $j=1$, i. e.  $\sigma_{\ell+1}=\sigma_1$. As a consequence,  one has $\sigma_{\ell+2}=\psi_\ell(\sigma_{\ell+1})=\psi_\ell(\sigma_{1})$ with $\sigma_{\ell+2} \sigma_{1}(z)=\sigma_{\ell+2} \sigma_{\ell+1}(z)=\ttz=\sigma_2\sigma_1(z)$. Therefore  $\sigma_{\ell+2}=\sigma_2$ and $\sigma_2=\psi_\ell(\sigma_1)$. By induction it follows that 
\[\overline{(\sigma_2,\sigma_1)}:=\big\{(\sigma_{2},\sigma_1),(\sigma_{3},\sigma_2), \ldots,(\sigma_{\ell+1},\sigma_\ell)\big\}.\]
Then one easily checks that for any $(\sigma_{k+1},\sigma_k)\in \overline{(\sigma_2,\sigma_1)}$, one has $\overline{(\sigma_{k+1},\sigma_k)}= \overline{(\sigma_2,\sigma_1)}$. It follows that the set 
\[{\mathcal C}(\ttz)=\big\{\overline{(\tau,\sigma)}\,\big|\, (\tau,\sigma)\in U(\ttz)\big\},\]
is a partition of $U(\ttz)$.

For $c=\overline{(\sigma_2,\sigma_1)}$ as above,  one denotes by
\[s(c):=\{\sigma_1(z),\ldots,\sigma_\ell(z)\}\subset ]z,\ttz[.\]
Observe that $c$ and $s(c)$ have same number of elements. 
We claim that all sets $s(c)$ are pairwise disjoints. Indeed, note that if $c'\in {\mathcal C}(\ttz)$ with $s(c)\cap s(c')\neq \emptyset$, then there exist $\sigma,\tau,\tau'\in \Sc$ such that $\sigma(z),\tau(z)\in s(c)$, $\sigma(z),\tau'(z)\in s(c')$ and $\tau\sigma(z)=\ttz=\tau'\sigma(z)$. It follows that  $\tau=\tau'$,  $(\tau,\sigma)\in c\cap c'$ and therefore $c=c'$, $s(c)=s(c')$. 
Since ${\mathcal C}(\ttz)$ is a partition of $U(\ttz)$, we finally get that the collection of sets  
$\{s(c)\,|\,c\in {\mathcal C}(\ttz)\}$ is also a partition of $]z,\ttz[$.

Let $\alpha:S_1(z)\to \R^+$ such that $\sum_{\tz\in S_1(z)} \alpha(\tz)=1$.  Given $\ttz\in S_2(z)$, applying the arithmetic-geometric mean inequality gives
\begin{align*}
    \big|]z,\ttz[\big| \Big(\prod_{\tz\in  ]z,\ttz[} {\alpha(\tz)}\Big)^{\frac{2}{|]z,\ttz[|}}&=\Big(\sum_{c\in {\mathcal C}(\ttz)}|s(c)|\Big)\prod_{c\in {\mathcal C}(\ttz)} \Big(\prod_{\tz\in s(c)} \alpha(z')\Big)^{\frac{2}{\sum_{c\in {\mathcal C}(\ttz)}|s(c)|}}\\
    &\leq \sum_{c\in {\mathcal C}(\ttz)} |s(c)| \Big(\prod_{\tz\in s(c)} \alpha(z')\Big)^{\frac{2}{|s(c)|}}
\end{align*}
Since $|s(c)|=|c|$, observing that 
\[\Big(\prod_{\tz\in s(c)} \alpha(z')\Big)^2= \prod_{(\tau,\sigma)\in c} \alpha(\tau(z))\alpha(\sigma(z)),\]
and applying again the 
arithmetic-geometric mean inequality, one gets 
\[\big|]z,\ttz[\big| \Big(\prod_{\tz\in  ]z,\ttz[} {\alpha(\tz)}\Big)^{\frac{2}{|]z,\ttz[|}} \leq \sum_{c\in {\mathcal C}(\ttz)} \sum_{(\tau,\sigma)\in c} \alpha(\tau(z))\alpha(\sigma(z)).\]

Given $c\in {\mathcal C}(\ttz)$, either $|c|=1=|s(c)|$ and there exists $\tau\in \Sc$ such that $c=\big\{(\tau,\tau)\}$ and $d(z,\tau\tau(z))=2$, either $|c|\geq 2$ and for any $(\tau,\sigma)\in c$, $\tau(z)\neq \sigma(z)$.  
As a consequence, setting 
\[s_1(\ttz):=\bigcup_{c\in {\mathcal C}(\ttz), |c|=1} s(c),\]
one has 
\[\sum_{c\in {\mathcal C}(\ttz)} \sum_{(\tau,\sigma)\in c} \alpha(\tau(z))\alpha(\sigma(z))\leq \sum_{\tz\in s_1(\ttz)} \alpha(\tz)^2  +\sum_{c\in {\mathcal C}(\ttz),|c|\geq 2} \sum_{(\tau,\sigma)\in c} \alpha(\tau(z))\alpha(\sigma(z)).\] 
and therefore 
\begin{align*}
    \sum_{\ttz\in S_2(z)} \big|]z,\ttz[\big| \Big(\prod_{\tz\in  ]z,\ttz[} {\alpha(\tz)}\Big)^{\frac{2}{|]z,\ttz[|}} &\leq \sum_{\ttz\in S_2(z)}\,\, \sum_{\tz\in s_1(\ttz)} \alpha(\tz)^2  \\
    &\qquad +\sum_{\ttz\in S_2(z)}\,\,\sum_{c\in {\mathcal C}(\ttz),|c|\geq 2} \sum_{(\tau,\sigma)\in c} \alpha(\tau(z))\alpha(\sigma(z))\\
    &\leq \sum_{\tz \in S_1(z)} \alpha(\tz)^2 + \sum_{((\tz,w')\in S_1(z),w'\neq \tz} \alpha(\tz)\alpha(w')\\
    &=\Big(\sum_{\tz \in S_1(z)} \alpha(\tz)\Big)^2=1,
\end{align*}
where the last inequality is a consequence of the fact that for $\ttz\neq w''$, one has  $U(\ttz)\cap U(w'')=\emptyset$ and also  $s_1(\ttz)\cap s_1(w'')=\emptyset$. Then according to the definition of $K_0(z,S_2(z))$, one has $K_0(z,S_2(z))\leq 1$.

If for any $\sigma\in \Sc$, $d(z,\sigma\sigma(z))\leq 1$, then for any $\ttz\in S_2(z)$, $s_1(\ttz)=\emptyset$.
Applying Cauchy Schwarz inequality, the above estimates provide
\begin{align*}
\sum_{\ttz\in S_2(z)} \big|]z,\ttz[\big| \Big(\prod_{\tz\in  ]z,\ttz[} {\alpha(\tz)}\Big)^{\frac{2}{|]z,\ttz[|}} &\leq \sum_{((\tz,w')\in S_1(z),w'\neq \tz} \alpha(\tz)\alpha(w')\\
&= \Big(\sum_{\tz \in S_1(z)} \alpha(\tz)\Big)^2-\sum_{\tz \in S_1(z)} \alpha(\tz)^2\\
&\leq 1-\frac1{|S_1(z)|},
\end{align*}
and therefore $K_0(z,S_2(z))\leq 1-1/|S_1(z)|$. The last statement of Theorem \ref{thmstructure} then follows from Theorem \ref{thmprinc}.
\end{proof}

\begin{proof}[Proofs of Theorem \ref{Logsob} and Theorem \ref{Logsobbis}]
Let $\nu_0$ and $\nu_1$ be  probability measures with same convex bounded support $\Cc\subset \X$ and respective densities $f_0$ and $f_1$ with respect to the measure $m$. Modified logarithmic Sobolev inequalities will follow from the convexity property 
\begin{equation}\label{versHWI}
 \HR(\nu_0|m)\leq - \frac{ \HR(\widehat \nu_t|m)-\HR(\nu_0|m)}{t} + \HR(\nu_1|m)- \frac12 (1-t)C_t(\widehat\pi), \qquad t\in(0,1),
\end{equation}
as $t$ goes to zero, with $C_t(\widehat\pi)={\widetilde{\kappa}}  \,\widetilde{T}(\widehat \pi)$ in Theorem \ref{Logsob},  and $C_t(\widehat\pi)={\widetilde{\kappa}_3}  \,\widetilde{T}_3(\widehat \pi)$ or $C_t(\widehat\pi)={\widetilde{\kappa}_3}  \,\widetilde{C}_t(\widehat \pi)$ in Theorem \ref{Logsobbis}.   

Let us start with the proof of Theorem \ref{Logsob}. Observing  that for any $x,y\in \X$ and $z\in[x,y]$,
\begin{multline*}\partial_t \nu_t^{x,y}(z)_{|t=0}=\frac{L^{d(x,z)}(x,z)L^{d(z,y)}(z,y)}{L^{d(x,y)}(x,y)}\,\binom{d(x,y)}{d(x,z)} \left( \1_{[x,y]}(z)\1_{z\sim x}-d(x,y)\1_{x=z}\right)\\
=\sum_{x'\in S_1(x)\cap[x,y]}  d(x,y) 
\frac{L(x,x')L^{d(x',y)}(x',y)}{L^{d(x,y)}(x,y) }\left(\delta_{x'}(z)-\delta_x(z)\right)
\end{multline*}
and since for any $t\in[0,1]$ the finite convex subset $\Cc$ is the support of $\widehat \nu_t$, one gets 
\begin{align}\label{tempete}
&\partial_t \HR(\widehat \nu_t|m)_{|t=0}\nonumber=\sum_{z\in \Cc} \partial_t \widehat\nu_t(z)_{|t=0} \log f_0(z)\\
&=\nonumber\sum_{z\in \Cc}  \sum_{x,y\in \Cc}  \sum_{x'\in S_1(x)\cap[x,y]} d(x,y) \frac{L(x,x')L^{d(x',y)}(x',y)}{L^{d(x,y)}(x,y) }\left(\delta_{x'}(z)-\delta_x(z)\right) \log f_0(z) \,\widehat \pi(x,y)\\
&= \sum_{x,y\in \Cc}\sum_{x'\in S_1(x)\cap[x,y]}\left(\log f_0(x')-\log f_0(x)\right) d(x,y) \frac{L(x,x')L^{d(x',y)}(x',y)}{L^{d(x,y)}(x,y) } \widehat \pi(x,y)\\
&\nonumber\geq -\sum_{x\in\Cc} \max_{x',x'\sim x}[\log f_0(x)-\log f_0(x')]_+ \Big(\sum_{y\in \Cc} d(x,y) \widehat \pi_\rightarrow(y|x) \Big) \,\nu_0(x)\\
&\nonumber\geq - \frac1{2 \widetilde{\kappa}} \sum_{x\in\Cc} \max_{x',x'\sim x}[\log f_0(x)-\log f_0(x')]_+^2 \nu_0(x) - \frac{\widetilde{\kappa}}2 \sum_{x\in\Cc} \Big(\sum_{y\in \Cc} d(x,y) \widehat \pi_\rightarrow(y|x) \Big)^2 \,\nu_0(x),
\end{align}
where for the last inequality, one uses the inequality $ab\leq a^2/2+b^2/2$, $a,b\in \R$. 
Therefore, from the definition of $\widetilde{T}(\widehat \pi)\geq \widetilde{T}_{_\rightarrow}(\widehat \pi)$, \eqref{versHWI} implies as $t$ goes to zero 
\[\HR(\nu_0|m)\leq  \frac1{2 \widetilde{\kappa}} \sum_{x\in\Cc} \max_{x',x'\sim x}[\log f_0(x)-\log f_0(x')]_+^2 \nu_0(x) +\HR(\nu_1|m) .\]
By choosing $\nu_1=\mu_\Cc:=\frac{\1_\Cc m}{m(\Cc)}$ and $f_0:=\frac{f\1_\Cc}{m(f\1_\Cc)}$, it gives 
 \[   {\rm Ent}_{\mu_\Cc}(f) 
    \leq \frac1{2 \widetilde{\kappa}} \int \max_{x',x'\sim x}[\log f(x)-\log f(x')]_+^2 f(x) \,d\mu_\Cc(x).\]
Applying this inequality with  $\Cc=\Cc_n$ where $(\Cc_n)$ is an increasing sequence of convex subsets with $\bigcup_{n}\Cc_n=\X$, the  monotone convergence theorem provides the expected  modified logarithmic Sobolev inequality \eqref{logsob} for $\mu$ since $f$ is bounded and $m(\X)<+\infty$.

Let $g:\X\to \R$ be a bounded function such that $\mu(g)=0$. As usual, applying \eqref{logsob} to the function $f=1+\varepsilon g $ where $\varepsilon$ is a sufficiently small parameter so that $f>0$, a  Taylor expansion as $\varepsilon$ goes to zero gives
\[\frac{\varepsilon^2}{2}\mu(g^2) +\circ(\varepsilon^2)\leq \frac{\varepsilon^2}{2\widetilde{\kappa}} 
\int \max_{x', x'\sim x}\left[g(x)-g(x')\right]_+^2 d\mu(x)+\circ(\varepsilon^2)
.\]
It provides the  Poincaré inequality of Theorem \ref{Logsob} as $\varepsilon$ goes to zero.

The proof of Theorem \ref{Logsobbis} $\textit{(ii)}$  and $\textit{(iii)}$ is similar. Starting again from equality \eqref{tempete}, one has
\begin{align*}
\partial_t \HR(\widehat \nu_t|m)_{|t=0}
&=\sum_{x,y\in \Cc} \sum_{\sigma\in \Sc,  \sigma(x)\in ]x,y]} \partial_\sigma \log f_0(x) \, \frac{L(x,\sigma(x))L^{d(\sigma(x),y)}(\sigma(x),y)}{L^{d(x,y)}(x,y) } \widehat \pi(x,y)\\
&\geq -\sum_{x\in \Cc} \sum_{\sigma\in \Sc} \left[\partial_\sigma \log f_0(x)\right]_-\Pi^\sigma_\rightarrow(x)\,\nu_0(x) \\
&\geq - \frac1{2\widetilde  \kappa_2 } \sum_{x\in \Cc} \sum_{\sigma\in \Sc} 
[\partial_\sigma( \log f_0)(x)]_-^2 \nu_0(x) - \frac{\widetilde  \kappa_2}2 \sum_{x\in\Cc} \sum_{\sigma\in \Sc}\left(\Pi^\sigma_\rightarrow(x)\right)^2\,\nu_0(x)
\end{align*}
The above inequality together with  \eqref{versHWI}  imply as $t$ goes to zero 
\begin{equation*}
\HR(\nu_0|m)\leq \frac1{2\widetilde  \kappa_2 } \int \sum_{\sigma\in \Sc} [\partial_\sigma( \log f)]_-^2 d\nu_0 +\HR(\nu_1|m) - \frac{\widetilde  \kappa_2}2 \int \sum_{\sigma\in \Sc} \left(\Pi^\sigma_\leftarrow(y)\right)^2 d\nu_1(y).
\end{equation*}
Then the end of  the  proof of the first part of Theorem \ref{Logsobbis}  is similar to the one of Theorem \ref{Logsob} with approximation's arguments. For the proof of its second part, one uses the inequality
\begin{align*}
&\partial_t \HR(\widehat \nu_t|m)_{|t=0}
\geq -\sum_{x\in \Cc} \sum_{\sigma\in \Sc} \left[\partial_\sigma \log f_0(x)\right]_-\Pi^\sigma_\rightarrow(x)\,\nu_0(x) \\
&\geq -  \sum_{x\in\Cc} \sum_{\sigma\in \Sc} \frac{\widetilde \kappa_2 D^2}2 \, h^*\left(\frac{2}{D\widetilde\kappa_2} [\partial_\sigma( \log f)(x)]_-\right)  \nu_0(x) -  \sum_{x\in\Cc} \sum_{\sigma\in \Sc}\frac{\widetilde{\kappa}_2 D^2}2h\left(\frac{\Pi^\sigma_\rightarrow(x)}D\right)\,\nu_0(x), 
\end{align*}
where $h^*(v):=\sup_{0\leq u<1}\big\{uv-h(u)\big\}=2\left(e^{-v/2}+v/2-1\right), v\geq 0$.
Since 
\[\lim_{t\to 0} \widetilde C_t^D(\widehat \pi)= \int \sum_{\sigma\in \Sc} D^2 h\left(\frac{\Pi^\sigma_\rightarrow(x)}{D}\right) d\nu_0(x)+\int \sum_{\sigma\in \Sc} D^2 h_{1}\left(\frac{\Pi^\sigma_\leftarrow(y)}D\right) d\nu_1(y), \]
with $h_1(u):=uh'(u)-h(u)=2(-u-\log(1-u))$, $u\in[0,1)$, 
as before inequality  \eqref{versHWI}  implies as $t$ goes to zero 
\begin{equation*}
\HR(\nu_0|m)\leq  \int \sum_{\sigma\in \Sc} \frac{\widetilde \kappa_2 D^2}2 \, h^*\left(\frac{2}{D\widetilde\kappa_2} [\partial_\sigma( \log f)]_-\right)  d\nu_0 +\HR(\nu_1|m) - \int \sum_{\sigma\in \Sc} D^2 h_{1}\left(\frac{\Pi^\sigma_\leftarrow(y)}D\right) d\nu_1(y).
\end{equation*}
The proof of the second part of Theorem \ref{Logsobbis}  ends as the one of Theorem \ref{Logsob}. By applying inequality \eqref{logsobT3} to a function $f=1+\varepsilon g$ with $g$ bounded and $\mu(g)=0$, a Taylor expansion as $\varepsilon$ goes to 0 implies
\[{\rm Var}_\mu(g)\leq \frac{1}{\widetilde{\kappa}_2} 
\int \sum_{\sigma\in \Sc}[\partial_\sigma g]_-^2 d\mu
. \]
Applying this inequality to $-g$ and adding these two inequalities provide the Poincaré inequality of Theorem \ref{Logsobbis}.

It remains to prove  part (i) of Theorem \ref{Logsobbis}. The proof of the first transport-entropy inequality is the same as the one of Corollary \ref{Transport}. In order to get the second one \eqref{transportentC_0^D} assuming that condition \eqref{restr} holds, since $\HR(\widehat \nu_t|\mu)\geq 0$, Theorem \ref{Thmstructure} ensures that  for any probability measure $\nu_0,\nu_1\in \Pc(\{0,1\}^n)$, and for any $t\in(0,1)$,
\[\frac{\widetilde\kappa_2}{2} \inf_{\pi\in\Pi(\nu_0,\nu_1)} \widetilde C_t^D(\pi)\leq \frac{1}{t}\HR(\nu_0|\mu_v)+\frac{1}{1-t} \HR(\nu_1|\mu_v).\]
Then \eqref{transportentC_0^D} easily follows by choosing $\nu_0=\mu_v$, $\nu_1=\nu\in \Pc(\X)$ and  letting  $t$ goes to $0$. 
\end{proof}

\begin{proof}[Proof of Theorem \ref{ttens}] $(\mathcal{X},d,m,L)$ is the Cartesian product 
of the graph spaces \[(\mathcal{X}_{i},d_{\X_{i}},m_{i},L_{i}), \qquad i\in[n].\]
According to the structure of Cartesian product of graph spaces, one has for any $z=(z_1,\ldots,z_n)\in \X$
\begin{multline*}
    K_{L} \big(z,S_{2}(z)\big)=\sup_\alpha \Big\{  \sum_{i=1}^{n}\sum_{z_{i}^{\prime \prime} \in S_{2}(z_{i})} 
L_{i}^{2}(z_{i},z_{i}^{\prime \prime})
\Big(\prod_{z_{i}^{\prime} \in ]z_{i},z_{i}^{\prime \prime}[}  \frac{\alpha_i(z_{i}^{\prime})}{L_{i}(z_{i},z_{i}^{\prime})}\Big)^{\frac{2L_{i}(z_{i},z_{i}^{\prime})L_{i}(z_{i}^{\prime},z_{i}^{\prime \prime})}{L_{i}^{2}(z_{i},z_{i}^{\prime \prime}) }} 
\\+ 2 \sum_{\{i,j\}\subset[n]} \sum_{z_{i}^\prime,z_{i}^\prime\sim z_{i}}\sum_{ z_{j}^\prime,z_{j}^\prime\sim z_{j}}\alpha_i(z_{i}^\prime)\alpha_j(z_{j}^\prime) \Big\},
\end{multline*}
where the supremum runs over all non negative vector $\alpha$ with coordinates $\alpha_i(z_i')$, $i\in[n]$, $z_i'\sim z_i$, such that $\sum_{i\in [n]}\sum_{z_{i}^\prime,z_{i}^\prime\sim z_{i}} \alpha_i(z_i')=1$. Setting $\alpha_i=\sum_{z_{i}^\prime,z_{i}^\prime\sim z_{i}} \alpha_i(z_i')$, and according to the definition of $K_{L_{i}}\big(z_{i},S_{2}(z_{i})\big)$, it follows that
\begin{align*}
K_{L} \big(z,S_{2}(z)\big)&=\sup_\alpha 
\Big\{  \sum_{i=1}^{n} \alpha_i^2  K_{L_{i}}\big(z_{i},S_{2}(z_{i})+ 2 \sum_{\{i,j\}\subset[n]} \alpha_i\alpha_j
\Big\}\\
&=1-\inf_\alpha\Big\{\sum_{i=1}^{n} \alpha_i^2 \big(1- K_{L_{i}}\big(z_{i},S_{2}(z_{i})\big)\Big\}\\
&\leq 1-\inf_\alpha\Big\{\big(1-\max_{i\in[n]}K_{L_{i}}\big(z_{i},S_{2}(z_{i})\big)\sum_{i=1}^{n} \alpha_i^2 \Big\},
\end{align*}
where the supremum is over all $\alpha=(\alpha_1,\ldots,\alpha_n)\in \R_+^n$ with $\sum_{i=1}^{n} \alpha_i=1$. The expected result follows from the sign of $\big(1-\max_{i\in[n]}K_{L_{i}}\big(z_{i},S_{2}(z_{i})\big)$ and since $\inf_\alpha\Big\{\sum_{i=1}^{n} \alpha_i^2 \Big\}=1/n$ and  $\sup_\alpha\Big\{\sum_{i=1}^{n} \alpha_i^2 \Big\}=1$.
This concludes the first part of Theorem \ref{ttens}.

Let us now study the constant $\widetilde r_2(z)$. By easy induction arguments, it suffices to get the result for $n=2$. Let $W\subset S_2(z)$. The set $]z,W[$ is the disjoint union of the two sets 
\[]z,W[^{1}:=\{z^{\prime} \in ]z,W[\hspace{0.06cm}| \hspace{0.06cm} z^{\prime}=(z_{1}^{\prime},z_{2}), z_{1}\sim z_{1}^{\prime} \}\]
and \[ ]z,W[^{2}:=\{ z^{\prime}\in ]z,W[ \hspace{0.06cm}| \hspace{0.06cm} 
z^{\prime}=(z_{1},z_{2}^{\prime}), z_{2}\sim z_{2}^{\prime} \}.\] 
Let $V_1:=\big\{z_1'\,\big|\,(z_1',z_2)\in ]z,W[^{1}\big\}$ and $V_2:=\big\{z_2'\,\big|\,(z_1,z_2')\in ]z,W[^{2}\big\}$.
Similarly, the set $W$ is the disjoint union of the three following  sets
\[\overline{W}_1:=\big\{(z_1'',z_2)\,|\, d_{\X_1}(z_1,z_1'')=2\big \},\quad
 \overline{W}_2:=\big\{(z_1,z_2'')\,|\, d_{\X_2}(z_2,z_2'')=2\big \},\]
 \[\mbox{and}\quad \overline{W}_3:=\big\{(z_1',z_2')\,|\, z_1\sim z_1', z_2\sim z_2'\big \}.\]
 One has
 \[]z,\overline{W}_1[\cup]z,\overline{W}_2[\cup]z,\overline{W}_3[=]z,W[^{1}\cup]z,W[^{2}.\]
 If $W_1:=\big\{z_1''\,\big|\,(z_1'',z_2)\in \overline{W}_1\big\}$ and $W_2:=\big\{z_2''\,\big|\,(z_1,z_2'')\in \overline{W}_{2}\big\}$, then
one has $]z_1,W_1[\subset V_1$ and $]z_2,W_2[\subset V_2$.
With the above notations and according to the definition \eqref{defKtilde} of $\widetilde K_L(z,W)$,  the structure of product of graphs gives 
\begin{align*}
   1-\widetilde K_{L_1\oplus L_2}(z,W)
   &=
   \inf_\beta  \Biggl\{ \bigg(\sum_{z_1'\in V_1} \sqrt{\beta_1(z^{\prime}_1)} +\sum_{z_2' \in V_2} \sqrt{\beta_2(z^{\prime}_2)}\bigg)^{2} -2\sum_{{(z_1',z_2') \in \overline{W}_3}} \sqrt{\beta_1(z^{\prime}_1)} \sqrt{\beta_2(z^{\prime}_2)}\\
   &\qquad\qquad-\sum_{z_1''\in W_1} L^{2}_1(z_1,z''_1) \prod_{z_1',\in ]z,z_1''[}\bigg(\frac{\beta_1(z^{\prime}_1)}{\big(L_1(z_1,z^{\prime}_1)\big)^{2}}\bigg)^{\frac{L_1(z_1,z^{\prime}_1)L_1(z_1^{\prime},z''_1)}{L^{2}_1(z_1,z''_1)}}\nonumber\\
   &\qquad\qquad-\sum_{z_2''\in W_2} L^{2}_2(z_2,z''_2) \prod_{z_2',\in ]z_2,z_2''[}\bigg(\frac{\beta_2(z^{\prime}_2)}{\big(L_2(z,z^{\prime}_2)\big)^{2}}\bigg)^{\frac{L_2(z_2,z^{\prime}_2)L_2(z_2^{\prime},z''_2)}{L^{2}_2(z_2,z''_2)}}\nonumber\Biggr\}, 
\end{align*}
where the supremum runs over all non negative vector $\beta$ with coordinates $\beta_i(z_i')$, $i\in[2]$, $z_i'\in V_i$, such that $\sum_{i\in [n]}\sum_{z_{i}^\prime\in V_i} \beta_i(z_i')=1$.
 Since $\overline W_3\subset V_1\times V_2$, it follows  
\begin{align*}
   1-&\widetilde K_{L_1\oplus L_2}(z,W)\\
   &\geq
   \inf_\beta  \Biggl\{\bigg(\sum_{z_1' \in V_1} \sqrt{\beta_1(z^{\prime}_1)}\bigg)^{2}-\sum_{z_1''\in W_1} L^{2}_1(z_1,z''_1) \prod_{z_1',\in ]z_1,z_1''[}\bigg(\frac{\beta_1(z^{\prime}_1)}{\big(L_1(z_1,z^{\prime}_1)\big)^{2}}\bigg)^{\frac{L_1(z_1,z^{\prime}_1)L_1(z_1^{\prime},z''_1)}{L^{2}_1(z_1,z''_1)}}\\
   &+\bigg(\sum_{z_2' \in V_2} \sqrt{\beta_2(z^{\prime}_2)}\bigg)^{2}-\sum_{z_2''\in W_2} L^{2}_2(z_2,z''_2) \prod_{z_2',\in ]z_2,z_2''[}\bigg(\frac{\beta_2(z^{\prime}_2)}{\big(L_2(z_2,z^{\prime}_2)\big)^{2}}\bigg)^{\frac{L_2(z_2,z^{\prime}_2)L_2(z_2^{\prime},z''_2)}{L^{2}_2(z_2,z''_2)}}\Biggl\}
\end{align*}
For $i\in[2]$, $]z_i,W_i[\subset V_i$. Therefore, setting $\beta_i=\sum_{z_{i}^\prime\in V_i} \beta_i(z_i')$, from the definition of $\widetilde K_{L_i}(z_i,W_i)$ one gets 
\begin{align*}
   1-\widetilde K_{L_1\oplus L_2}(z,W)&\geq
   \inf_{\beta_1+\beta_2=1}  \Big\{\beta_1 \big(1-\widetilde K_{L_1}(z_1,W_1)\big) +\beta_2 \big(1-\widetilde K_{L_2}(z_2,W_2)\big)\Big\}\\ 
   &=\min\big(1-\widetilde K_{L_1}(z_1,W_1), 1-\widetilde K_{L_2}(z_2,W_2)\big) \geq \min\big(\widetilde r_1(z_1), \widetilde r_2(z_2)\big). 
\end{align*}
 The second part of  Theorem \ref{ttens} then follows optimizing over all $W\subset S_2(z)$. 
\end{proof}

\subsubsection{Proofs of Lemma \ref{Dvcube}, Lemma \ref{lemcubeK_v} and Lemma \ref{KLvZ}}

\begin{proof}[Proof of Lemma \ref{Dvcube} ]

We want to bound from below  the quantity defined by \eqref{mulh}, $D_t v(x,y)$, for any $x,y\in \{0,1\}^n$ with $d=d(x,y)\geq 2$. The identity \eqref{Dvijcube} provides 
\begin{align*} 
&D_t v(x,y)\\&=2\sum_{ z\in[x,y]}\;\; \sum_{\{i,j\}\subset[n],(z,\sigma_i\sigma_j(z))\in[x,y]}\!\!\!\!\!\!\!\!\!\!\!\!(2z_i-1)(2z_j-1)\,\partial_{ij}^2 v(z_{\overline{ij}}) \,r(x,z,\sigma_i\sigma_j(z),y)\rho_t^{d-2}(d(x,z)).
\end{align*}
If $\sigma_i\sigma_j(z)\in[x,y]$ then $(2z_i-1)=x_i-y_i$ and $(2z_j-1)=x_j-y_j$. As a consequence, if $\partial_{ij}^2 v(z_{\overline{ij}}):=V_{ij}$ does not depend on $z_{\overline{ij}}$, one has 
\begin{align*}
D_t v(x,y)&=2\sum_{\{i,j\}\subset[n]} (x_i-y_i)(x_j-y_j)V_{ij}\!\!\!\!\!\!\!\sum_{ z\in[x,y],(z,\sigma_i\sigma_j(z))\in[x,y]}\!\!\!\!\!\!\!\!\!\!\!\!r(x,z,\sigma_i\sigma_j(z),y)\rho_t^{d-2}(d(x,z)) \\
&=\frac2{d(d-1)}  \sum_{\{i,j\}\subset[n]} (x_i-y_i)(x_j-y_j)V_{ij},
\end{align*}
which ends the proof of the first part of Lemma \ref{Dvcube}.
In any case, when $\partial_{ij}^2 v(z_{\overline{ij}})$ depends on $z_{\overline{ij}}$, we also have 
\begin{align*}
D_t v(x,y)
&=\sum_{k=0}^{d-2} \ell_{t}^{x,y}(k)\frac{k!(d-2-k)!}{d!} \,\rho_t^{d-2}(k)\\&= \sum_{k=1}^{d-1} \ell_{t}^{x,y}(k-1)\frac{(k-1)!(d-k-1)!}{d!} \,\rho_t^{d-2}(k-1),
\end{align*}
with 
\begin{align*} \ell_{t}^{x,y}(k)&:=2 \sum_{z\in [x,y], d(x,z)=k}\sum_{ \{i,j\}\subset[n],(z,\sigma_i\sigma_j(z))\in[x,y]}(2z_i-1)(2z_j-1)\,\partial_{ij}^2 v(z_{\overline{ij}})\\
&=2 \sum_{z\in [x,y], d(x,z)=k} \sum_{ \{i,j\}\subset[n]} (2z_i-1)\1_{z_i\neq y_i}(2z_j-1)\1_{z_j\neq y_j} \partial_{ij}^2 v(z_{\overline{ij}}) ,
\end{align*}
or by symmetry,
\[\ell_{t}^{x,y}(k)=2 \sum_{z\in [x,y], d(x,z)=k+2}\sum_{ \{i,j\}\subset[n]}(2z_i-1)\1_{z_i\neq x_i}(2z_j-1)\1_{z_j\neq x_j} \partial_{ij}^2 v(z_{\overline{ij}}).\]

It follows that  for $k\in\{1,\ldots , d-1\}$,
\begin{align*}
\ell_{t}^{x,y}(k-1)&\geq \sum_{z\in [x,y], d(x,z)=k-1} \lambda_{\min}(Hv(z))  \sum_{i\in [n]} \1_{z_i\neq y_i}(2z_i-1)^2\\
&\geq\lambda_{\min}^\infty(Hv)\, (d-k+1) \,\frac{d!}{(k-1)!(d-k+1)!} \,
\end{align*}
and by symmetry
\[\ell_{t}^{x,y}(k-1)\geq \lambda_{\min}^\infty(Hv)\, (k+1) \,\frac{d!}{(k+1)!(d-k-1)!},\]
Since $Hv(z)$ has  off-diagonal entries,  $\lambda_{\min}^\infty(Hv)\leq \lambda_{\min}(Hv(z))\leq 0$ and therefore we get   
\begin{align*}
D_t v(x,y)&\geq \lambda_{\min}^\infty(Hv) \sum_{k=1}^{d-1} \min\{1/k,1/(d-k)\}\,\rho_t^{d-2}(k-1)\\
&=\lambda_{\min}^\infty(Hv)  \sum_{k=1}^{d-1} \min\{d-k,k\}\,\frac{\rho_t^{d}(k)}{t(1-t)d(d-1)}\\
&\geq \frac{\lambda_{\min}^\infty(Hv)}{2(d-1)}\, \frac{1-\rho_t^{d}(0)-\rho_t^{d}(n)}{t(1-t)}=\frac{\lambda_{\min}^\infty(Hv)}{2(d-1)}\,
\gamma_t(d). 
\end{align*}
Then inequality \eqref{subtile} provides the expected result, 
\[\int_0^1 D_s v(x,y)\,q_t(s)\,ds\geq \frac{\lambda_{\min}^\infty(Hv)}{2(d-1)}\int_0^1 \gamma_s(d)\,ds\geq  \frac{\lambda_{\min}^\infty(Hv)}{d-1}\,{\sum_{k=1}^{d-1}\frac1k}\geq \lambda_{\min}^\infty(Hv).\]
\end{proof}

\begin{proof}[Proof of Lemma \ref{lemcubeK_v}] Let $z\in \{0,1\}^n$ and $W\subset S_2(z)$.
Using the definition of the subset of indices $A^1$ associated to the set $W$ and given by
\eqref{hard}, the quantity $K^v(s,W)$ given by \eqref{defKvtilde} can be written as
\[K^v(s,W)=\sup_{\alpha}\Big\{ 2\sum_{\{i,j\}\subset A^1}  e^{-(2z_i-1)(2z_j-1)\,\partial_{ij}^2v(z_{\overline{ij}})/2} \alpha_i\alpha_j\Big\},\]
where the infimum runs over all $\alpha=(\alpha_i)_{i\in A^1}$ with positive coordinates $\alpha_i$ satisfying $\sum_{i\in A^1}\alpha_i=1$.
 The upper bound on $K^v(s,W)$ is a consequence of the inequality $e^s\leq 1+s+|s|k(|s|)$, $s\in \R$. It provides 
\begin{align*}
&2\sum_{\{i,j\}\subset A^1}  e^{-(2z_i-1)(2z_j-1)\,\partial_{ij}^2v(z_{\overline{ij}})/2} \alpha_i\alpha_j\\&\leq2 \!\!\!\!\sum_{\{i,j\}\subset A^1} \!\!\!\!\!\!\Big(1-\frac12 (2z_i-1)(2z_j-1)\,\partial_{ij}^2v(z_{\overline{ij}})\Big) \alpha_i\alpha_j  +k\Big(\frac{|Hv|_{\max,\infty}}{2}\Big)\sum_{\{i,j\}\subset A^1} \!\!\!\!\big|\partial_{ij}^2v(z_{\overline{ij}})\big| \alpha_i\alpha_j
\\
&\leq 1- \Big(\sum_{i\in A^1} \alpha_i^2\Big) \Big[1+\frac{\lambda_{\min}(Hv(z))}2-k\Big(\frac{|Hv|_{\max,\infty}}{2}\Big) \frac{\lambda_{\max}(|Hv|(z))}2\Big]&\\
&\leq 1- \Big(\sum_{i\in A^1} \alpha_i^2\Big) r(v)
\end{align*}
which ends the proof of the first part of Lemma \ref{lemcubeK_v}.

For the second part of Lemma \ref{lemcubeK_v}. The proof is similar for the upper bound of $\widetilde  K^v(z)=\sup_{W\in S_2(z)}\widetilde  K^v(z,W)$ with according to \eqref{defKvtilde}
\[\widetilde  K^v(s,W)=\sup_{\beta}\Big\{ 2\sum_{\{i,j\}\subset A^1}  \Big(e^{-(2z_i-1)(2z_j-1)\,\partial_{ij}^2v(z_{\overline{ij}})/2}-1\Big) \sqrt{\beta_i}\sqrt{\beta_j}\Big\},\]
where the infimum runs over all $\beta=(\beta_i)_{i\in A^1}$ with positive coordinates $\beta_i$ satisfying $\sum_{i\in A^1}\beta_i=1$.
As above it follows that 
\begin{align*}
&2\sum_{\{i,j\}\subset A^1} \Big(e^{-(2z_i-1)(2z_j-1)\,\partial_{ij}^2v(z_{\overline{ij}})/2}-1\Big) \sqrt{\beta_i}\sqrt{\beta_j} 
\\&\leq-\!\!\!\!\!\sum_{\{i,j\}\subset A^1} \!\!\!\!\! (2z_i-1)(2z_j-1)\,\partial_{ij}^2v(z_{\overline{ij}})\sqrt{\beta_i}\sqrt{\beta_j}   +k\Big(\frac{|Hv|_{\max,\infty}}{2}\Big)\!\!\!\!\!\sum_{\{i,j\}\subset A^1} \!\!\!\!\!\big|\partial_{ij}^2v(z_{\overline{ij}})\big| \sqrt{\beta_i}\sqrt{\beta_j} \\
&\leq -\frac{\lambda_{\min}(Hv(z))}2+k\Big(\frac{|Hv|_{\max,\infty}}{2}\Big) \frac{\lambda_{\max}(|Hv|(z))}2\leq 1-r(v),
\end{align*}
which implies the expected upper bound on $\widetilde K^v=\sup_{z\in \{0,1\}^n}\widetilde K^v(z)$.
For the lower bound on $\widetilde K^v$, since for any $z\in\{0,1\}^n$,  $A^1=[n]$ for $W=S_2(z)$, the inequality $e^s-1\geq s$ gives 
 \begin{align*}
\widetilde K^v(z)&\geq \widetilde K^v(s,S_2(z))=\sup_{\beta}\Big\{ 2\sum_{\{i,j\}\subset [n]}  \Big(e^{-(2z_i-1)(2z_j-1)\,\partial_{ij}^2v(z_{\overline{ij}})/2}-1\Big) \sqrt{\beta_i}\sqrt{\beta_j}\Big\}\\
&\geq \sup_{\beta}\Big\{-2\sum_{\{i,j\}\subset [n]} (2z_i-1)(2z_j-1)\,\partial_{ij}^2v(z_{\overline{ij}})\sqrt{\beta_i}\sqrt{\beta_j}\\
&\geq -\lambda_{\min}(Hv(z)).
\end{align*}
It follows that $\widetilde K^v\geq -\lambda_{\min}^\infty(Hv)$. 
\end{proof}

\begin{proof}[Proof of Lemma \ref{KLvZ}]
We want to upper bound the quantity $K^v(z, S_2(z))$ for any $z\in \Z^n$ whose expression is given by \eqref{defKLv}. According to the structure of the lattice $\Z^n$ and from the identity \eqref{DvZn}, one has for any $z\in \Z^n$,
\begin{multline*}
   K^v(z, S_2(z)):=\sup_\alpha\Big\{ 2 \sum_{\{i,j\}\subset[n]} \Big(\alpha_{i+}\alpha_{j+}e^{-\partial_{ij}v(z)/2}+\alpha_{i-}\alpha_{j-}e^{-\partial_{ij}v(z-e_i-e_j)/2} \\
   +\alpha_{i+}\alpha_{j-}e^{\partial_{ij}v(z-e_j)/2}
   +\alpha_{i-}\alpha_{j+}e^{\partial_{ij}v(z-e_i)/2}\Big) 
   +\sum_{i\in[n]} \Big(\alpha_{i+}^2 e^{-\partial_{ii}v(z)/2}+ \alpha_{i-}^2 e^{-\partial_{ii}v(z-2e_i)/2}\Big)\Big\},  
\end{multline*}
where the supremum runs over all vectors $\alpha$ with non-negative coordinates $\alpha_{i+},\alpha_{i-}$ satisfying $\sum_{i\in[n]} (\alpha_{i+}+\alpha_{i-})=1$.
According to the definition of the matrix $Av(z)$ in Lemma \ref{KLvZ}, one has 
\begin{align*}
   &K^v(z, S_2(z))\\&\leq \sup_\alpha\Big\{ 2\!\!\!\!\!\sum_{\{i,j\}\subset[n]} \!\!\!\!\! \big(\alpha_{i+}+\alpha_{i-}\big)\big(\alpha_{j+}+\alpha_{j-}\big)\big((Av(z))_{ij} +1\big)  +\!\sum_{i\in[n]} \!\!\big(\alpha_{i+}+ \alpha_{i-}\big)^2\big((Av(z))_{ii} +1\big)\Big\}\\
   &\leq 1 +\lambda_{\max}\big(Av(z)\big)\sum_{i\in[n]} \big(\alpha_{i+}+ \alpha_{i-}\big)^2\\
   &\leq 1+\frac{\lambda_{\max}\big(Av(z)\big)}n,
\end{align*}
where the last inequality is a consequence of Cauchy-Schwarz inequality if $\lambda_{\max}\big(A(z)\big)\leq 0$.

We want now to upper bound  $\widetilde K^v(z, W)$ for any $W\subset S_2(z)$. According to \eqref{defKtilde} it can be expressed as follows 
\begin{align*}
\widetilde{K}^v(z,W)&=\sup  \Biggl\{  e^{-D v(z,\ttz)/2}\sum_{\ttz\in W} |]z,\ttz[| \Big(\prod_{\sigma\in \Sc_{]z,\ttz[}} {\beta(\sigma)}\Big)^{\frac1{|]z,\ttz[|}}-\!\!\!\!\!\!\!\!\!\!\!\!\!\sum_{(\sigma,\tau)\in \Sc_{]z,W[}^2, \sigma\neq \tau}\!\!\!\!\!\!\!\!\!\!\!\!\! \sqrt{\beta(\sigma)}\sqrt{\beta(\tau)}
 \\&\qquad\qquad\qquad\qquad\qquad\qquad\qquad\qquad\qquad\Bigg|\,{\beta}: \Sc_{]z,W[}\to \mathbb{R}_{+},\sum_{\sigma\in \Sc_{]z,W[}} \beta(\sigma)=1 \Biggr\}.
\end{align*}
Given $W\subset S_2(z)$, there exist  subsets  $I_+,I_-\subset[n]$, $J_+, J_-\subset\{(i,j)\in[n]\times [n]\,|\, i< j\}$ and $K\subset\{(i,j)\in[n]\times [n]\,|\, i\neq j\}$, such that 
\begin{multline*}
    W=\big\{z+2e_i\,\big|\,i \in I_+\big\}\cup\big\{z-2e_i\,\big|\,i \in I_-\big\}\cup \big\{z+e_i+e_j\,\big|\, (i,j)\in J_+\big\}\\
    \cup\big\{z-e_i-e_j\,\big|\, (i,j)\in J_-\big\}\cup\big\{z+e_i-e_j\,\big|\, (i,j)\in K\big\}.
\end{multline*}
Setting $L_+=I_+\cup \{i\,|\, \exists k\in [n], (i,k)\in J_+ \cup  K\mbox{ or } (k,i)\in J_+\}$ and $L_-=I_-\cup \{j\,|\, \exists k\in [n], (k,j)\in J_- \cup  K\mbox{ or } (j,k)\in J_-\}$, the identity \eqref{DvZn} gives 
\begin{align*}
\widetilde{K}^v(z,W)&=\sup_\beta 
\Big\{ 2 \sum_{(i,j)\in J_+} \sqrt{\beta_{i+}}\sqrt{\beta_{j+}}e^{-\partial_{ij}v(z)/2} +2\sum_{(i,j)\in J_-} \sqrt{\beta_{i-}}\sqrt{\beta_{j-}}e^{-\partial_{ij}v(z-e_i-e_j)/2} \\&\quad+2 \sum_{(i,j)\in K} \sqrt{\beta_{i+}}\sqrt{\beta_{j-}}e^{\partial_{ij}v(z-e_j)/2}+ \sum_{i\in I_+} \beta_{i+} e^{-\partial_{ii}v(z)/2} +  \sum_{i\in I_-} \beta_{i-} e^{-\partial_{ii}v(z-2e_i)/2}\\
&\quad-\sum_{(i,j)\in L_+, i\neq j} \sqrt{\beta_{i+}}\sqrt{\beta_{j+}} -\sum_{(i,j)\in L_-, i\neq j} \sqrt{\beta_{i-}}\sqrt{\beta_{j-}} -2\sum_{(i,j)\in L_+\times L_-} \sqrt{\beta_{i+}}\sqrt{\beta_{j-}}\Big\},
\end{align*}
where the supremum runs over all vectors $\beta$ with non-negative coordinates $\beta_{i+}, i\in L_+,\beta_{j-}, j\in L_-$ satisfying $\sum_{i\in L_+}\beta_{i+}+\sum_{j\in L_-}\beta_{j-}=1$.
The definition of the matrix $Av(z)$ then provides
\begin{align*}
\widetilde{K}^v(z,W)
&\leq 1+\sup_\beta 
\Big\{ 2 \sum_{\{i,j\}\subset L_+} \sqrt{\beta_{i+}}\sqrt{\beta_{j+}}((Av(z))_{ij}+1) +2\sum_{\{i,j\}\subset L_-} \sqrt{\beta_{i-}}\sqrt{\beta_{j-}}((Av(z))_{ij}+1) \\&\quad+2 \sum_{(i,j)\in L_+\times L_-, i\neq j} \sqrt{\beta_{i+}}\sqrt{\beta_{j-}}((Av(z))_{ij}+1)+ \sum_{i\in L_+} \beta_{i+} (Av(z))_{ii} +  \sum_{i\in L_-} \beta_{i-} (Av(z))_{ii}\\
&\quad-\sum_{(i,j)\in L_+, i\neq j} \sqrt{\beta_{i+}}\sqrt{\beta_{j+}} -\sum_{(i,j)\in L_-, i\neq j} \sqrt{\beta_{i-}}\sqrt{\beta_{j-}} -2\sum_{(i,j)\in L_+\times L_-} \sqrt{\beta_{i+}}\sqrt{\beta_{j-}}\Big\}\\
&=1+ \sup_\beta 
\Big\{ 2 \sum_{\{i,j\}\subset L_+} \sqrt{\beta_{i+}}\sqrt{\beta_{j+}}(Av(z))_{ij} +2\sum_{\{i,j\}\subset L_-} \sqrt{\beta_{i-}}\sqrt{\beta_{j-}}(Av(z))_{ij}\\&\qquad\qquad\qquad+2 \sum_{(i,j)\in L_+\times L_-,i\neq j} \sqrt{\beta_{i+}}\sqrt{\beta_{j-}}(Av(z))_{ij}+ \sum_{i\in L_+} \beta_{i+} (Av(z))_{ii} \\
&\qquad\qquad\qquad\qquad+  \sum_{i\in L_-} \beta_{i-} (Av(z))_{ii}-2\sum_{i\in L_+\cap L_-} \sqrt{\beta_{i+}}\sqrt{\beta_{i-}}\Big\}
\end{align*}
Since $-1\leq (Av(z))_{ii}$ it follows that 
\begin{align*}
&\widetilde{K}^v(z,W)\\&\leq 1+\sup_\beta 
\Big\{2 \sum_{\{i,j\}\subset [n], i\neq j} \big(\sqrt{\beta_{i+}}\1_{i\in L_+}+\sqrt{\beta_{i-}}\1_{i\in L_-}\big)\big(\sqrt{\beta_{j+}}\1_{j\in L_+}+\sqrt{\beta_{j-}}\1_{j\in L_-}\big)(Av(z))_{ij}\\
&\qquad\qquad\qquad\qquad\qquad\qquad +\sum_{i\in [n]} \big(\sqrt{\beta_{i+}}\1_{i\in L_+}+\sqrt{\beta_{i-}}\1_{i\in L_-}\big)^2 (Av(z))_{ii}\Big\}\\
&\leq 1+\sup_\beta\Big\{\lambda_{\max}\big(Av(z)\big) \sum_{i\in [n]} \big(\sqrt{\beta_{i+}}\1_{i\in L_+}+\sqrt{\beta_{i-}}\1_{i\in L_-}\big)^2\Big\}\leq 1+ \lambda_{\max}\big(Av(z)\big),
\end{align*}
where the last inequality holds since $\lambda_{\max}\big(Av(z)\big)\leq 0$.
Finally we get   \[\widetilde{K}^v(z)=\sup_{W\subset S_2(z)}\widetilde{K}^v(z,W)\leq 1+ \lambda_{\max}\big(Av(z)\big).\] This ends the proof of Lemma \ref{KLvZ}.
\end{proof}

\subsection{Proof of Proposition \ref{Dieppe} and   Proposition \ref{propcycles}}

\begin{proof}[Proof of Proposition \ref{Dieppe}]
Let $z\in \X$. According to hypothesis \eqref{Hypo1} one may consider the finest partition $\{J_1,\ldots ,J_q\}$, $q\in \N^*$, of $S_2(z)$ such that for all distinct $i,j\in [k]$, for all $\ttz\in J_i$ and all $w''\in J_j$,
$]z,\ttz[\,\cap\,]z,w''[=\emptyset$.
According to  \eqref{defRbis} and from the above property of the $J_i$'s, one has
\begin{align*}
   K_0(z,S_2(z)) =\sup_{\alpha} \sum_{i=1}^q \sum_{\ttz\in J_i} |]z,\ttz[| \prod_{z'\in]z,\ttz[}\alpha(\tz)^{\frac2{|]z,\ttz[|}} 
\end{align*}
where the supremum runs over all  function $\alpha:S_1(z)\to \R_+$ with 
\[\sum_{i=1}^q \sum_{\tz\in ]z,J_i[} \alpha(\tz)=1.\]
For any $i\in[q]$, let 
\[K(J_i):=\sup_\beta \sum_{\ttz\in J_i} |]z,\ttz[| \prod_{z'\in]z,\ttz[}\beta(\tz)^{\frac2{|]z,\ttz[|}},\]
where the supremum runs over all  function $\beta:J_i\to \R_+$ with 
\[ \sum_{\tz\in ]z,J_i[} \beta(\tz)=1.\]
By homogeneity, we also have 
\[K_0\big(z,S_2(z)\big) =\sup_{\delta} \sum_{i=1}^q \delta_i^2 K(J_i)=\sup_{i\in[q]} K(J_i), 
\]
since the supremum runs over all vector $\delta$ with non-negative coordinates $\delta_i$, $i\in[q]$, satisfying $\delta_1+\cdots +\delta_q=1$.
Therefore it remains to show that $K(J_i)\leq 7/8$ for any $i\in[q]$.

Let $J$ denotes an arbitrary set $J_i$ of the partition of $S_2(z)$. 
Let $n:=|J|$ and $J=\{z_1'',\ldots,z_n^{\prime\prime}\}$. For simplicity let us denote $M_j:=]z,z_j^{\prime\prime}[\subset S_1(z)$ the set of midpoints between $z$ and $z_j^{\prime\prime}$ and $m_j:=|M_j|$ its cardinality, 
\[M_j:=\big\{w^j_1,\ldots,w^j_{m_j}\big\}.\]  Observe that since $\{J_1,\ldots, J_q\}$ is the finest partition all subsets $M_j$ are connected by intersection, which means that for any two sets $M_j$ and $M_{j'}$ there exists a sequence of sets $M_{j_0},M_{j_1},\ldots, M_{j_k}$ with $M_{j_0}=M_j$, $M_{j_k}=M_{j'}$ and $M_{j_\ell}\cap M_{j_{\ell-1}}\neq \emptyset$ for any $\ell\in [k]$.

Let $A_2$ denotes the subset of indices  $j\in [n]$ such that $m_j=2$ and $A_3$ its complementary (for any $j\in A_3$, $m_j\geq 3$). From the arithmetic-geometric mean inequality, one has 
\begin{align*}
    K(J)&=\sup_\beta \Big\{\sum_{j\in A_2} 2 \beta(w^j_1)\beta(w^j_2) + \sum_{j\in A_3} m_j \big(\beta(w^j_1)\cdots \beta(w^j_{m_j})\big)^{2/m_j}\Big\} \\
    &\leq \sup_\beta \Big\{\sum_{j\in A_2} 2 \beta(w^j_1)\beta(w^j_2) + \sum_{j\in A_3} \Big[\beta(w^j_1)\beta(w^j_{m_j})+\sum_{k=1}^{m_j-1} \beta(w^j_k)\beta(w^j_{k+1}) \Big]\Big\}.
\end{align*}
Let $G^*:=(V^*,E^*)$ denotes the graph with set of vertices $V^*:=]z,J[=\cup_{j\in [n]} M_j$ and set of edges 
\[E^*:=\bigcup_{j\in[n]}\Big\{ \{w^j_1,w^j_{m_j}\}, \{w^j_1,w^j_{2}\}, \{w^j_2,w^j_{3}\},\ldots, \{w^j_{m_j-1},w^j_{m_j}\}\Big\}.\]
The restriction of $G^*$ to a set of vertices  $M_j$, $j\in A_3,$ is a cycle. For convenience, let $M_j$  also denotes this cycle.
According to hypothesis \eqref{Hypo1} any vertex $w$ of $V^*=\cup_{j\in [n]} M_j$ do not belong to more than two different sets of the collection of $M_j$, $j\in[n]$. Therefore it is also the case for any edge  of $E^*$. 

Observe that if $\{w,w'\}$ is an edge of $E^*$ that belongs to two different sets $M_i$ and $M_j$, then necessarily either $m_i\geq 
3$ or $m_j\geq 3$. Indeed if $m_i=m_j=2$ then $M_i=M_j$ and for $W=\{z_i^{\prime\prime},z_j^{\prime\prime}\}$, one has $|]z,W[|=|M_1|=2$, that contradicts hypothesis \eqref{Hypo2}.

As a consequence, one has 
\[\sum_{j\in A_2} 2 \beta(w^j_1)\beta(w^j_2) + \sum_{j\in A_3} \Big[\beta(w^j_1)\beta(w^j_{m_j})+\sum_{k=1}^{m_j-1} \beta(w^j_k)\beta(w^j_{k+1}) =\!\!\!\!\!\!\!\! \sum_{\{w,w'\}\in E^*} \!\!\!\!\!\!\!\!b(w,w') \beta(w)\beta(w'),\]
where  coefficients  $b(w,w')$ belong to $\{1,2,3\}$.
Note that 
\begin{itemize}
\item $ b(w,w')=1 $  if and only if $\{w,w'\}\in E^*$ is a subset of a single set $M_j$ with $j\in A_3$. Let $E_1^*$ denotes the set of these edges.
\item $  b(w,w')=2 $ if and only if $\{w,w'\}\in E^*$ is a subset of a single set $M_i$ with $i\in A_2$, or $\{w,w'\}\in E^*$ is a subset of two different  sets $M_j$ and $M_{j'}$ with $j,j'\in A_3$. Let $E_2^*$ denotes the set of these edges.
\item $ b(w,w')=3 $  if and only if $\{w,w'\}$ is a subset of a set  $M_j$ with $j\in A_3$ and another set $M_i$, $i\in A_2$. Let $E_3^*$ denotes the set of these edges.
\end{itemize}
It follows that 
\begin{eqnarray}\label{dieppe}
\sum_{\{w,w'\}\in E^*} b(w,w') \beta(w)\beta(w') =2 \sum_{\{w,w'\}\in E_2^*\cup E_3^*} \beta(w)\beta(w')+ \sum_{\{w,w'\}\in E_1^*\cup E_3^*} \beta(w)\beta(w').
\end{eqnarray}
The proof  of Proposition \ref{Dieppe} ends by bounding the maximum clique number of the graph $G_1^*$ generated by the set of edges $E_1^*\cup E_3^*$ and the graph $G_2^*$ generated by the set of edges  $E_2^*\cup E_3^*$. 

Recall that according to hypothesis \eqref{Hypo1}, three distinct sets $M_i$, $M_j$ and $M_k$ do not intersect.
It follows that two edges of $E_2^*\cup E_3^*$ share a same vertex if and only if one of the following holds
\begin{itemize}
    \item either the two edges belongs to the same two different cycles $M_j$ and $M_{j'}$ with $j,j'\in A_3$,
    \item either each of the edges belongs to a different set  $M_i$ with $i\in A_2$.
\end{itemize}
This observation implies that  if $\{w_1,w_2,w_3\}$ is the set of vertices of a triangle in $G_2^*$, then 
\begin{itemize}
    \item either $\{w_1,w_2\},\{w_2,w_3\}, \{w_3,w_1\}$ belongs to two different $M_j$ and $M_{j'}$ with $j,j'\in A_3$. But since $M_j$ and $M_{j'}$ are cycles, this implies $M_j=M_{j'}=\{w_1,w_2,w_3\}$ which is nonsense, 
    \item either each of the edges $\{w_1,w_2\},\{w_2,w_3\}, \{w_3,w_1\}$ belongs to a different set  $M_i$ with $i\in A_2$. Without loss of generality for example $M_1=\{w_1,w_2\}$, $M_2=\{w_2,w_3\}$ and $M_3=\{w_3,w_1\}$. Due to hypothesis \eqref{Hypo1}, these sets $M_i$ can not be connected by intersection to other sets $M_j$. It follows that $J:=\{z_1^{\prime\prime},z_2^{\prime\prime},z_3^{\prime\prime}\}$
    and 
    \[K(J):=\sup_{\beta} \Big\{2 \beta(w_1,w_2)+2 \beta(w_2,w_3) +2 \beta(w_3,w_1)\Big\}=2/3.\]
\end{itemize}
We conclude that either the maximum clique number of $G_2^*$ is less than 2 either $K(J)= 2/3$.

Let us now consider the maximum click number of $G_1^*$.
Two edges of $E_1^*\cup E_3^*$ share a same vertex if an only if one of the following holds 
\begin{itemize}
    \item the two edges belongs to the same $M_j$, $j\in A_3$,
    \item one of the edges belongs to some $M_j$, $j\in A_3$ and the second to another $M_{j'}$ $j'\in A_3$,
    \item one of the edges belongs to some $M_j$, $j\in A_3$ and the second also belongs to $M_j$ and to some $M_i$ with $i\in A_2$.
\end{itemize}

Therefore, if $\{w_1,w_2,w_3\}$ is the set of vertices of a triangle in $G_1^*$ then either the $\{w_1,w_2,w_3\}$ corresponds to a cycle $M_j$ of cardinality 3, either the three edges  $\{w_1,w_2\},\{w_2,w_3\}$ and  $\{w_3,w_1\}$ belongs to three different cycles $M_j$, either two of same are in the same cycle $M_j$ and the third  is a different one. 
Actually, assume that each cycle $M_j$ as a different color, one may easily check that with the constrain that any vertex can be share with at most 2 colors, one may construct 4 types of complete graphs $K_4$ with 2 or 3 colors, and only the two following  types of complete graph $K_5$ with vertex set $\{w_1,w_2,w_3,w_4,w_5\}$, namely  
\begin{itemize}
    \item either with $M_1:=\{w_1,w_2,w_3,w_4\}$, $M_2=\{w_2,w_4,w_5\}$ and $M_3=\{w_1,w_3,w_5\}$,
    \item or $M_1:=\{w_1,w_2,w_3,w_4,w_5\}$ and $M_2:=\{w_1,w_4,w_2,w_5,w_3\}$
\end{itemize}
In that cases since any vertex is already shared by 2 cycles $M_i$ it can not be shared by another one. This implies that $G_1^*=G^*$
and according to Theorem \ref{cliquethm}, $K(J)= \frac12\Big(1-\frac15\Big)=\frac25$. 
In the above both cases of  configurations of the $M_i$'s for constructing a complete graph $K_5$, the size of each cycle $M_i$ is  fixed, therefore  $G^*_1$ can not contain any $K_6$. 
We conclude that either the maximum clique number of $G_1^*$ is less than 4 or $K(J)= 2/5$.

According to \eqref{dieppe} and from the last observation, applying Theorem \ref{cliquethm}, we finally get
\begin{align*}
    K(J)&\leq\max\Big\{\frac23,\frac25,2 \sup_\beta  \sum_{\{w,w'\}\in E_2^*\cup E_3^*} \beta(w)\beta(w')+ \sup_\beta \sum_{\{w,w'\}\in E_1^*\cup E_3^*} \beta(w)\beta(w')\Big\}\\
    &\leq \max\Big\{\frac23,\frac25, \Big(1-\frac12\Big)+\frac12\Big(1-\frac14\Big)\Big\}=\frac78,
\end{align*}
which ends the proof of Proposition \ref{Dieppe}.
\end{proof}

\begin{proof}[Proof of Proposition \ref{propcycles}]

Let $x,y$ be two vertices of a graph $G$ such that $x\sim y$. Let suppose that $g(G)<5$ then there exist $x^{\prime}\sim x$ and 
$y^{\prime}\sim y$ such that $d(x^{\prime},y^{\prime})\leq 1$. The measure 
$m_{x}^{\alpha}$  can be expressed as follows
\begin{equation*}
    m_{x}^{\alpha}=\Big(1-\frac{\text{deg}(x)}{\Delta(G)}\alpha \Big)\delta_{x}+\frac{\alpha}{\Delta(G)}\delta_{S_{1}(x)},
\end{equation*}
and $m_{y}^{\alpha}$ in an analogous way.
Without loss of generality, let suppose that $\text{deg}(y)\geq \text{deg}(x)$ so that  
    $S_{1}(x)\setminus \{x^{\prime},y\}=\{z_{1},\ldots , z_{\text{deg}(x)-2}\}$ and 
    $S_{1}(y)\setminus \{x,y^{\prime}\}=\{z_{1}^{\prime},\ldots , z^{\prime}_{\text{deg}(y)-2}\}$.
For $\alpha\leq \frac{\Delta(G)}{\text{deg}(y)+1}$, let  $\pi$ be the coupling probability measure with first marginal $m_{x}^{\alpha}$ and second marginal $m_{y}^{\alpha}$ given by
    \[\pi(x,x)=\pi(y,y)=\frac{\alpha}{\Delta(G)},\quad
    \pi(x,y)=1-\frac{\text{deg}(y)}{\Delta(G)}\alpha-\frac{\alpha}{\Delta(G)}, \quad
    \pi(x^{\prime},y^{\prime})=\frac{\alpha}{\Delta(G)},\]
    \vspace{-0,15 cm}
    \[\pi(z_{i},z^{\prime}_{i})=\frac{\alpha}{\Delta(G)} \mbox{ for } 1\leq i\leq \text{deg}(x)-2,\quad\pi(x,z^{\prime}_{j})=\frac{\alpha}{\Delta(G)} \mbox{ for } \text{deg}(x)-1\leq j\leq  \text{deg}(y)-2.\]
Since  $d(x,y)=1$, $d(x^{\prime},y^{\prime})\leq 1$,   $d(z_{i},z_{i}^{\prime})\leq 3$ and $d(x,z^{\prime}_{j})\leq 2$, the definition of the $W_1$-Wasserstein distance ensures that 
\begin{align*}
1&-\kappa_{\alpha}(x,y)=W_{1}(m_{x}^{\alpha},m_{y}^{\alpha})\\
&\leq 3\frac{\alpha}{\Delta(G)}(\text{deg}(x)-2)+2\frac{\alpha}{\Delta(G)}(\text{deg}(y)-\text{deg}(x))+\frac{\alpha}{\Delta(G)}+1-\frac{\text{deg}(y)}{\Delta(G)}\alpha-\frac{\alpha}{\Delta(G)}.
\end{align*}
It follows  that $\kappa_{LLY}(x,y)\geq \frac{6-\text{deg}(y)-\text{deg}(x)}{\Delta(G)}$ which is a contradiction.
\end{proof}
\bibliographystyle{plain}
\bibliography{criteria-curvature}

\begin{thebibliography}{10}

\bibitem{RJSS23}
Rados{\l}aw Adamczak, Joscha Prochno, Marta Strzelecka, and Micha{\l}
  Strzelecki.
\newblock Norms of structured random matrices.
\newblock {\em Mathematische Annalen}, pages 1--65, 2023.

\bibitem{AKPS18}
Radosław Adamczak, Michal Kotowski, Bartlomiej Polaczyk, and Michał
  Strzelecki.
\newblock A note on concentration for polynomials in the ising model.
\newblock {\em Electronic Journal of Probability}, 2018.

\bibitem{AGZ10}
Greg~W Anderson, Alice Guionnet, and Ofer Zeitouni.
\newblock {\em An introduction to random matrices}.
\newblock Number 118. Cambridge university press, 2010.

\bibitem{BE85}
Dominique Bakry and Michel {\'E}mery.
\newblock Diffusions hypercontractives.
\newblock In {\em Seminaire de probabilit{\'e}s XIX 1983/84}, pages 177--206.
  Springer, 1985.

\bibitem{BBH23}
Afonso~S Bandeira, March~T Boedihardjo, and Ramon van Handel.
\newblock Matrix concentration inequalities and free probability.
\newblock {\em Inventiones mathematicae}, pages 1--69, 2023.

\bibitem{BB19}
Roland Bauerschmidt and Thierry Bodineau.
\newblock A very simple proof of the lsi for high temperature spin systems.
\newblock {\em Journal of Functional Analysis}, 276(8):2582--2588, 2019.

\bibitem{BG99}
Sergey Bobkov and G{\"o}tze Friedrich.
\newblock Exponential integrability and transportation cost related to
  logarithmic sobolev inequalities.
\newblock {\em Journal of Functional Analysis}, 163:1--28, 1999.

\bibitem{BHT00}
Sergey Bobkov, Christian Houdr{\'e}, and Prasad Tetali.
\newblock $\lambda_\infty$, vertex isoperimetry and concentration.
\newblock {\em Combinatorica}, 20(2):153--172, 2000.

\bibitem{YR22}
Yuansi Chen and Ronen Eldan.
\newblock Localization schemes: A framework for proving mixing bounds for
  markov chains.
\newblock In {\em 2022 IEEE 63rd Annual Symposium on Foundations of Computer
  Science (FOCS)}, pages 110--122. IEEE, 2022.

\bibitem{HS13}
Hee~Je Cho and Seong-Hun Paeng.
\newblock Ollivier’s ricci curvature and the coloring of graphs.
\newblock {\em European Journal of Combinatorics}, 34(5):916--922, 2013.

\bibitem{Chu07}
Fan Chung.
\newblock Four proofs for the cheeger inequality and graph partition
  algorithms.
\newblock In {\em Proceedings of ICCM}, volume~2, page 378. Citeseer, 2007.

\bibitem{CY96}
Fan~RK Chung and S-T Yau.
\newblock Logarithmic harnack inequalities.
\newblock {\em Mathematical Research Letters}, 3(6):793--812, 1996.

\bibitem{CKKLP21}
David Cushing, Supanat Kamtue, Riikka Kangaslampi, Shiping Liu, and Norbert
  Peyerimhoff.
\newblock Curvatures, graph products and ricci flatness.
\newblock {\em Journal of Graph Theory}, 96(4):522--553, 2021.

\bibitem{CKL21}
David Cushing, Supanat Kamtue, Shiping Liu, and Norbert Peyerimhoff.
\newblock Bakry-{\'e}mery curvature on graphs as an eigenvalue problem.
\newblock {\em Calculus of Variations and Partial Differential Equations},
  61(2):62, 2022.

\bibitem{CLP20}
David Cushing, Shiping Liu, and Norbert Peyerimhoff.
\newblock Bakry--{\'e}mery curvature functions on graphs.
\newblock {\em Canadian Journal of Mathematics}, 72(1):89--143, 2020.

\bibitem{DLP09}
Jian Ding, Eyal Lubetzky, and Yuval Peres.
\newblock The mixing time evolution of glauber dynamics for the mean-field
  ising model.
\newblock {\em Communications in Mathematical Physics}, 289(2):725--764, 2009.

\bibitem{EKZ22}
Ronen Eldan, Frederic Koehler, and Ofer Zeitouni.
\newblock A spectral condition for spectral gap: fast mixing in
  high-temperature ising models.
\newblock {\em Probability theory and related fields}, 182(3-4):1035--1051,
  2022.

\bibitem{EHMT17}
Matthias Erbar, Christopher Henderson, Georg Menz, and Prasad Tetali.
\newblock Ricci curvature bounds for weakly interacting markov chains.
\newblock {\em Electron. J. Probab}, 22(40):1--23, 2017.

\bibitem{EM12}
Matthias Erbar and Jan Maas.
\newblock Ricci curvature of finite markov chains via convexity of the entropy.
\newblock {\em Archive for Rational Mechanics and Analysis}, 206(3):997--1038,
  2012.

\bibitem{EM14}
Matthias Erbar and Jan Maas.
\newblock Gradient flow structures for discrete porous medium equations.
\newblock {\em Discrete \& Continuous Dynamical Systems-A}, 34(4):1355, 2014.

\bibitem{FS18}
Max Fathi and Yan Shu.
\newblock Curvature and transport inequalities for markov chains in discrete
  spaces.
\newblock {\em Bernoulli}, 24(1):672--698, 2018.

\bibitem{GRST14}
Nathael Gozlan, Cyril Roberto, Paul-Marie Samson, and Prasad Tetali.
\newblock Displacement convexity of entropy and related inequalities on graphs.
\newblock {\em Probability Theory and Related Fields}, 160(1-2):47--94, 2014.

\bibitem{GRST14bis}
Nathael Gozlan, Cyril Roberto, Paul-Marie Samson, and Prasad Tetali.
\newblock Kantorovich duality for general transport costs and applications.
\newblock {\em Journal of Functional Analysis}, 273(11):3327--3405, 2017.

\bibitem{JL14}
J{\"u}rgen Jost and Shiping Liu.
\newblock Ollivier’s ricci curvature, local clustering and
  curvature-dimension inequalities on graphs.
\newblock {\em Discrete \& Computational Geometry}, 51(2):300--322, 2014.

\bibitem{Kam20}
Supanat Kamtue.
\newblock A note on a bonnet-myers type diameter bound for graphs with positive
  entropic ricci curvature.
\newblock {\em arXiv preprint arXiv:2003.01160}, 2020.

\bibitem{Karp}
Richard~M Karp.
\newblock {\em Reducibility among combinatorial problems}.
\newblock Springer, 2010.

\bibitem{KKRT16}
Bo\'az Klartag, Gady Kozma, Peter Ralli, and Prasad Tetali.
\newblock Discrete curvature and abelian groups.
\newblock {\em Canadian Journal of Mathematics}, 68(3):655–674, 2016.

\bibitem{HL82}
Jean-Marie Laborde and Surya Prakash~Rao Hebbare.
\newblock Another characterization of hypercubes.
\newblock {\em Discrete Mathematics}, 39(2):161--166, 1982.

\bibitem{LHY18}
Rafa{\l} Lata{\l}a, Ramon van Handel, and Pierre Youssef.
\newblock The dimension-free structure of nonhomogeneous random matrices.
\newblock {\em Inventiones mathematicae}, 214:1031--1080, 2018.

\bibitem{Leo16}
Christian L{\'e}onard.
\newblock Lazy random walks and optimal transport on graphs.
\newblock {\em The annals of Probability}, 44(3):1864--1915, 2016.

\bibitem{Leo17}
Christian L{\'e}onard.
\newblock On the convexity of the entropy along entropic interpolations.
\newblock In {\em Measure Theory in Non-Smooth Spaces}, pages 194--242. De
  Gruyter Open Poland, 2017.

\bibitem{LLY11}
Yong Lin, Linyuan Lu, and Shing-Tung Yau.
\newblock Ricci curvature of graphs.
\newblock {\em Tohoku Mathematical Journal, Second Series}, 63(4):605--627,
  2011.

\bibitem{YS10}
Yong Lin and Shing-Tung Yau.
\newblock Ricci curvature and eigenvalue estimate on locally finite graphs.
\newblock {\em Mathematical research letters}, 17(2):343--356, 2010.

\bibitem{LMP18}
Shiping Liu, Florentin M{\"u}nch, and Norbert Peyerimhoff.
\newblock Bakry--{\'e}mery curvature and diameter bounds on graphs.
\newblock {\em Calculus of Variations and Partial Differential Equations},
  57(2):1--9, 2018.

\bibitem{LMP17}
Shiping Liu, Florentin M{\"u}nch, and Norbert Peyerimhoff.
\newblock Rigidity properties of the hypercube via bakry--{\'e}mery curvature.
\newblock {\em Mathematische Annalen}, pages 1--35, 2022.

\bibitem{LV09}
John Lott and C{\'e}dric Villani.
\newblock Ricci curvature for metric-measure spaces via optimal transport.
\newblock {\em Annals of Mathematics}, pages 903--991, 2009.

\bibitem{Maa11}
Jan Maas.
\newblock Gradient flows of the entropy for finite markov chains.
\newblock {\em Journal of Functional Analysis}, 261(8):2250--2292, 2011.

\bibitem{Mar96}
Katalin Marton.
\newblock {Bounding $\bar{d}$-distance by informational divergence: a method to
  prove measure concentration}.
\newblock {\em The Annals of Probability}, 24(2):857 -- 866, 1996.

\bibitem{MT06}
Ravi Montenegro and Prasad Tetali.
\newblock Mathematical aspects of mixing times in markov chains.
\newblock {\em Theoretical Computer Science}, 1(3):237--354, 2006.

\bibitem{MS65}
Theodore~S Motzkin and Ernst~G Straus.
\newblock Maxima for graphs and a new proof of a theorem of tur{\'a}n.
\newblock {\em Canadian Journal of Mathematics}, 17:533--540, 1965.

\bibitem{MW19}
Florentin M{\"u}nch and Rados{\l}aw~K Wojciechowski.
\newblock Ollivier ricci curvature for general graph laplacians: Heat equation,
  laplacian comparison, non-explosion and diameter bounds.
\newblock {\em Advances in Mathematics}, 356:106759, 2019.

\bibitem{Oll09}
Yann Ollivier.
\newblock Ricci curvature of markov chains on metric spaces.
\newblock {\em Journal of Functional Analysis}, 256(3):810--864, 2009.

\bibitem{Oll11}
Yann Ollivier.
\newblock A visual introduction to riemannian curvatures and some discrete
  generalizations.
\newblock {\em Analysis and Geometry of Metric Measure Spaces: Lecture Notes of
  the 50th S{\'e}minaire de Math{\'e}matiques Sup{\'e}rieures (SMS),
  Montr{\'e}al}, 56:197--219, 2011.

\bibitem{OllVill12}
Yann Ollivier and C{\'e}dric Villani.
\newblock A curved brunn--minkowski inequality on the discrete hypercube, or:
  What is the ricci curvature of the discrete hypercube?
\newblock {\em Siam Journal on Discrete Mathematics}, 26(3):983--996, 2012.

\bibitem{Oys87}
Oystein Ore.
\newblock {\em Theory of graphs}, volume~38.
\newblock American Mathematical Soc., 1987.

\bibitem{Sam07}
Paul-Marie Samson.
\newblock Infimum-convolution description of concentration properties of
  product probability measures, with applications.
\newblock {\em Annales de l'I.H.P. Probabilités et statistiques},
  43(3):321--338, 2007.

\bibitem{Sam21}
Paul-Marie Samson.
\newblock Entropic curvature on graphs along schr{\"o}dinger bridges at zero
  temperature.
\newblock {\em Probability Theory and Related Fields}, 184(3-4):859--937, 2022.

\bibitem{Sch99}
Michael Schmuckenschl{\"a}ger.
\newblock Curvature of nonlocal markov generators.
\newblock {\em Convex geometric analysis (Berkeley, CA, 1996)}, 34:189--197,
  1999.

\bibitem{Stu06a}
Karl-Theodor Sturm.
\newblock On the geometry of metric measure spaces.
\newblock {\em Acta mathematica}, 196(1):65--131, 2006.

\bibitem{Tal95}
Michel Talagrand.
\newblock Concentration of measure and isoperimetric inequalities in product
  spaces.
\newblock {\em Publications Math{\'e}matiques de l'Institut des Hautes Etudes
  Scientifiques}, 81:73--205, 1995.

\bibitem{Tal96}
Michel Talagrand.
\newblock New concentration inequalities in product spaces.
\newblock {\em Inventiones mathematicae}, 126(3):505--563, 1996.

\bibitem{Tal10}
Michel Talagrand.
\newblock {\em Mean field models for spin glasses: Volume I: Basic examples},
  volume~54.
\newblock Springer Science \& Business Media, 2010.

\bibitem{TT21}
Akihisa Tamura and Kazuya Tsurumi.
\newblock Directed discrete midpoint convexity.
\newblock {\em Japan Journal of Industrial and Applied Mathematics}, 38:1--37,
  2021.

\bibitem{Vil09}
C{\'e}dric Villani.
\newblock Optimal transport, volume 338 of.
\newblock {\em Grundlehren der Mathematischen Wissenschaften [Fundamental
  Principles of Mathematical Sciences]}, 2009.

\end{thebibliography}
\end{document}